\documentclass{amsart}
\usepackage{amsmath,amssymb, amsbsy}
\usepackage{subfigure}
\usepackage{graphpap,latexsym,epsf}
\usepackage{color,psfrag}
\usepackage[dvips]{graphicx}
\usepackage[a4paper,width={15cm},left=3cm,bottom=2.5cm, top=2.5cm,
includeheadfoot]{geometry}
\usepackage{graphpap,latexsym,epsf}
\usepackage{color,graphicx,psfrag}

\newcommand{\R}{{\mathbb R}}
\newcommand{\N}{{\mathbb N}}

\newcommand{\C}{{\mathbb C}}
\newcommand{\SN}{{\mathbb S}^{N-1}}
\newcommand{\weakly}{\rightharpoonup}
\newcommand{\e }{\varepsilon}

\newcommand{\alchi}{\raisebox{1.7pt}{$\chi$}}
\newcommand{\dive }{\mathop{\rm div}}

\newcommand{\bi}{\mathbf b}
\newcommand{\tildO}{\widetilde \Omega}
\newcommand{\tildA}{\widetilde A}
\newcommand{\tildb}{\widetilde\bi}
\newcommand{\boldbeta}{\boldsymbol{\beta}}
\newcommand{\boldalfa}{\boldsymbol{\alpha}}
\renewcommand{\ge }{\geqslant}
\renewcommand{\geq }{\geqslant}
\renewcommand{\le }{\leqslant}
\renewcommand{\leq }{\leqslant}

\newenvironment{pf}{\noindent{\sc Proof}.\enspace}{\hfill\qed\medskip}
\newenvironment{pfn}[1]{\noindent{\bf Proof of
    {#1}.\enspace}}{\hfill\qed\medskip}
\newtheorem{Theorem}{Theorem}[section]
\newtheorem{Corollary}[Theorem]{Corollary}
\newtheorem{Lemma}[Theorem]{Lemma}
\newtheorem{Proposition}[Theorem]{Proposition}
\theoremstyle{definition}

\newtheorem{remark}[Theorem]{Remark}
\begin{document}

\title[Elliptic equations in domains with corners] {Almgren-type
  monotonicity methods
 for the classification\\ of behavior at corners
  of solutions to\\ semilinear elliptic equations }

\author[Veronica Felli \and Alberto Ferrero]{Veronica Felli \and
  Alberto Ferrero}
\address{\hbox{\parbox{5.7in}{\medskip\noindent{Universit\`a di Milano
        Bicocca,\\
        Dipartimento di Ma\-t\-ema\-ti\-ca e Applicazioni, \\
        Via Cozzi
        53, 20125 Milano, Italy. \\[3pt]
        \em{E-mail addresses: }{\tt veronica.felli@unimib.it,
          alberto.ferrero@unimib.it}.}}}}

\date{July 19, 2011}


 \begin{abstract}
   \noindent
   A monotonicity approach to the study of the asymptotic behavior
   near corners of solutions to semilinear elliptic equations in
   domains with a conical boundary point is discussed.  The presence
   of logarithms in the first term of the asymptotic expansion
   is excluded for boundary profiles sufficiently close to straight
   conical surfaces.
 \end{abstract}

\maketitle

\section{Introduction}\label{sec:intro}

This paper presents a monotonicity approach to the study of the
asymptotic behavior near corners of solutions to semilinear elliptic
equations
\begin{equation}\label{eq:92}
  -\dive(A(x)\nabla u(x))+\bi(x)\cdot\nabla
  u(x)-\dfrac{V\big(\frac{x}{|x|}\big)}{|x|^2}u(x)=
  h(x)u(x)+f(x,u(x))
\end{equation}
in a domain $\Omega\subset\R^N$, $N\geq 2$, having the origin as a
conical boundary point.  The coefficients $\bi:\Omega\to \R^N$ and
$h:\Omega\to \R$ are possibly singular at $0$ but satisfy suitable
decaying conditions (see assumptions (\ref{eq:b_assumption}) and
(\ref{eq:h_assumption}) below) which make the corresponding terms
negligible with respect to the homogeneity of the operator, while the
nonlinearity $f$ has at most critical growth in the Sobolev sense (see
assumptions (\ref{eq:F_assumption1}--\ref{eq:F_assumption2})).

Due to their own theoretical interest and their numerical application
to convergence analysis of finite element approximations, regularity
and asymptotics near corners of solutions to linear elliptic equations
in domains with piecewise boundary have been intensively studied and a
large literature has been devoted to this subject (see \cite{nistor},
the monographs \cite{dauge} and \cite[Chapter 3]{MNP1}, the surveys
\cite{Grisvard,ko}, and the references therein). Some early
contributions in this field date back to papers \cite{lehman,wasow}
which use methods based on conformal maps and integral representation
to derive asymptotic expansions for harmonic functions at a common
endpoint of two analytic arcs delimiting the $2$-dimensional simply
connected domain; such asymptotic development excludes the presence of
logarithmic terms for irrational values of $\alpha$, where $\alpha\pi$
is the opening of the corner. On the other hand, a simple example
shows that, if $\alpha=\frac{n}{m}$, $n,m\in \N\setminus\{0\}$, is a
rational number, then there exist harmonic functions with smooth trace
on the boundary of the domain but having a logarithmic term in the
leading part of the asymptotic expansion: it is sufficient to consider
the classical example $u(x,y)=\Im (z^m \log z)$, $z=x+iy$, in the domain
$$
\{(x,y)=(r\cos\theta,r\sin\theta) \in \R^2: r>0, \theta\in(0,\pi/m)\} .
$$
In \cite[Theorem 3.4]{lehman} logarithmic terms are excluded in the
leading expansion term in the case of homogeneous boundary conditions
also for rational values of $\alpha$.

Related results for semilinear Dirichlet
problems on plane domains with corners were obtained in
\cite{kawohl,wigley}; see also \cite{llm} for the study of existence
and nonexistence of solutions to singular semilinear elliptic
equations on cone-like domains.  We mention that edge asymptotics
(which is naturally related to corner asymptotics) is investigated in
\cite{cd} (see also the references therein).

In the spirit of the paper \cite{tolksdorf}, which provides
asymptotics of positive solutions to $p$-Laplace equations with
forcing terms and non-homogeneous boundary conditions on straight
$N$-dimensional cones, we mean to describe the rate and the shape of
solutions to (\ref{eq:92}) near corners of domains which are
perturbations of cones, by relating them to the eigenvalues and the
eigenfunctions of a limit operator on the spherical cap measuring the
opening of the vertex.  The method this paper is proposing for
valuating the asymptotic behavior of solutions to (\ref{eq:92}) is
based on the monotonicity method introduced by Almgren \cite{almgren}
in 1979 and then extended by Garofalo and Lin \cite{GL} to elliptic
operators with variable coefficients in order to prove unique
continuation properties. Monotonicity methods were recently used in
\cite{FFT,FFT2,FFT3} to prove not only unique continuation but also
precise asymptotics near singularities of solutions to linear and
semilinear elliptic equations with singular potentials, by extracting
such precious information from the behavior of the quotient associated
with the Lagrangian energy.  Almgren type formulas were also used in
\cite{ae} to prove unique continuation at the boundary; the
diffeomorphic deformation of the domain performed in \cite{ae} (see
also \cite{tz}) to get rid of the boundary contributions inspires our
construction of the equivalent problem (\ref{eq:tildequation}) in
section \ref{sec:an-equiv-probl}, for which a monotonicity formula is
derived in section \ref{sec:monotonicity-formula}.

As a byproduct of our asymptotic analysis we also obtain a unique
continuation principle for solutions of \eqref{eq:92} vanishing with
infinite order at the conical point of the boundary.

The strengths of the monotonicity formula approach are described in
the note \cite{FFT3}: they essentially rely in the sharpness of the
asymptotics derived, in the possibility of allowing quite general
perturbing potentials, and in the unified approach to linear and
nonlinear equations.

In subsection \ref{sec:assumpt-main-results} we introduce notation and
assumptions needed to state our main result Theorem \ref{t:main-u}.

\subsection{Assumptions and main results}\label{sec:assumpt-main-results}

For $N\geq 2$, let $\varphi:\R^{N-1}\to\R$ and $g:{\mathbb
  S}^{N-2}\to\R$ such that, for some $\delta>0$,
\begin{align}
\label{eq:phi1}
&\varphi(0)=0,\quad
\varphi\in C^2(\R^{N-1}\setminus\{0\}),\\[3pt]
\label{eq:g_continuous}
&g \in C^1({\mathbb S}^{N-2})\quad\text{if }N\geq3,\\[3pt]
&\label{eq:phi2} \sup_{\nu\in {\mathbb
    S}^{N-2}}\Big|\frac{\varphi(t\nu)}{t}-g(\nu)\Big|=O(t^\delta)
\quad\text{as }t\to0^+,\\[3pt]
&\label{eq:63}
\begin{cases}
  \sup_{\nu\in {\mathbb S}^{N-2}}\big|\nabla \varphi(t\nu)-g(\nu)\nu-
  \nabla_{\mathbb S^{N-2}}g(\nu)\big|=O(t^\delta)
  ,&\text{if }N\geq 3,\\
  \sup_{\nu\in \{-1,1\}}\big|\varphi'(t\nu)-g(\nu)\nu\big|=O(t^\delta)
,&\text{if
  }N=2
\end{cases}
\quad\text{as }t\to0^+,\\[3pt]
\label{eq:phi4}
&|D^2\varphi(x')|=O(|x'|^{-1})\quad \text{as }|x'|\to0.
\end{align}
As we will show in Lemma \ref{l:phi3}, 
 assumptions (\ref{eq:phi1}--\ref{eq:63}) imply that there exists
 $C_0>0$ such that
\begin{equation}\label{eq:phi3}
|\varphi(x')-\nabla\varphi(x')\cdot x'|\leq C_0|x'|^{1+\delta}
\text{ for all $x'$ in a neighborhood of $x'=0$.}
\end{equation}
Furthermore, from (\ref{eq:g_continuous}) it follows that
 the function $\varphi_0:\R^{N-1}\to\R$,
\begin{equation*}
\varphi_0(x'):=
\begin{cases}
|x'|g\big(\frac{x'}{|x'|}\big),&\text{if }x'\in\R^{N-1}\setminus\{0\},\\
0,&\text{if }x'=0,
\end{cases}
\end{equation*}
satisfies
$$
\varphi_0\in C^0(\R^{N-1})\quad \text{and}\quad
\varphi_0\in C^1(\R^{N-1}\setminus\{0\}).
$$
Hence the cone in $\R^N$ with vertex in $0$ defined as
\begin{equation}\label{eq:cilindro}
{\mathcal C}:=\big\{(x',x_N)\in \R^{N-1}\times \R:
x_N>\varphi_0(x')\big\}
\end{equation}
is open. In particular,
\begin{equation*}
C={\mathcal C}\cap {\mathbb S}^{N-1}
\end{equation*}
is an open connected subset of ${\mathbb S}^{N-1}$.

\begin{figure}[h]
 \centering
   \begin{psfrags}
     \psfrag{a}{$\varphi_0(x')$}
     \psfrag{b}{$\varphi(x')$}
\psfrag{c}{$C$}
\psfrag{d}{$\Omega$}
\psfrag{e}{$x'$}
\psfrag{f}{$x_N$}
     \includegraphics[width=7cm]{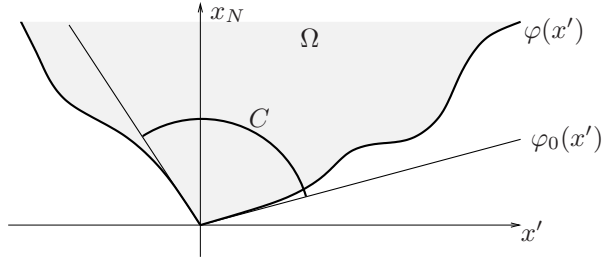}
   \end{psfrags}
 \caption{An example of domain $\Omega$.}\label{fig:dd}
\end{figure}

Let $\Omega$ be an open subset of $\R^N$ such that, for some  $R>0$,
\begin{equation}\label{eq:omega}
\Omega\cap B_R=\{x=(x',x_N)\in B_R: x_N>\varphi(x')\},
\end{equation}
where $B_R$ denotes the ball $\{x\in\R^N: |x|<R\}$
  in $\R^N$ with center at $0$ and radius $R$,
see figure \ref{fig:dd}.

Let $A:\Omega\to {\mathcal M}_{N\times N}$ (with ${\mathcal
  M}_{N\times N}$ denoting the space of $N\times N$  real matrices) satisfying
\begin{align}\label{eq:matrix1}
\begin{cases}
  a_{ij}=(A)_{ij}\in W^{1,\infty}(\Omega) \text{ for all
  }i,j=1,\dots,N, \quad a_{ij}=a_{ji}, \\
\text{there exists }C_A>0\text{ such that }A(x)\xi\cdot\xi
\geq C_A|\xi|^2\text{ for all }\xi\in\R^N\text{ and }
x\in\Omega.
\end{cases}
\end{align}
We observe that under assumption (\ref{eq:matrix1}), the functions
$a_{ij}$ are actually Lipschitz continuous functions on $\overline
\Omega$; moreover, due to symmetry and positive definiteness of $A$,
up to some change of variable, it is not restrictive to assume that
\begin{align}\label{eq:matrix2}
A(0)={\rm Id}_N,
\end{align}
where  ${\rm Id}_N$ denotes the identity $N\times N$ matrix.
Let us assume
\begin{align}
\label{eq:b_assumption}
& {\bi}\in L^\infty_{\rm loc}(\Omega,\R^N), \quad
|\bi(x)|=O(|x|^{-1+\delta})\quad\text{as }|x|\to0,\\
\label{eq:h_assumption}
& h\in
L^\infty_{\rm loc}(\Omega), \quad
h(x)=O(|x|^{-2+\delta})\quad\text{as }|x|\to0.
\end{align}
It is not restrictive to assume that the positive constants $\delta$'s
of formulas (\ref{eq:phi2}), (\ref{eq:63}), (\ref{eq:b_assumption}), and
(\ref{eq:h_assumption}) are the same and that $\delta\in(0,1)$.  Let
$V:\SN\to\R$ such that
\begin{align}\label{eq:V_assumption}
\begin{cases}
V\equiv 0,&\text{if }N=2,\\[5pt]
{{L_V:=\mathop{\sup}\limits_{\substack{\theta,\tau\in\SN\\\theta\neq\tau}}
\frac{|V(\theta)-V(\tau)|}{|\theta-\tau|}<+\infty}}\quad\text{and}\quad
\Lambda(V)<1,&\text{if }N\geq 3,
\end{cases}
\end{align}
where, for $N\geq3$,
\begin{equation}\label{eq:LambdaV}
\Lambda(V):=\sup_{v\in {\mathcal D}^{1,2}(\mathcal C)\setminus\{0\}}\dfrac{{\displaystyle
{\int_{\mathcal C}{{|x|^{-2}}{V(x/|x|)}\,v^2(x)\,dx}}}}
{{\displaystyle{\int_{\mathcal C}{|\nabla v(x)|^2\,dx}}}}
\end{equation}
and ${\mathcal D}^{1,2}(\mathcal C)$ denotes
the completion of $C^\infty_{\rm c}({\mathcal C})$ with respect to the
norm
$$
\|u\|_{{\mathcal D}^{1,2}(\mathcal C)}:=\bigg(\int_{\mathcal C}|\nabla
u(x)|^2\,dx\bigg)^{\!\!1/2}.
$$
Let $f:\Omega\times\R\to\R$ such that
\begin{align}\label{eq:F_assumption1}
f\in C^0(\Omega\times \R),\quad F\in C^1(\Omega\times \R),
\quad s\mapsto f(x,s)\in C^1(\R)\text{ for a.e. }x\in\Omega,
\end{align}
\begin{align}\label{eq:F_assumption2}
|f(x,s)s|+|f'_s(x,s)s^2|+|\nabla_x F(x,s)||x|\leq
\begin{cases}
C_f(|s|^2+|s|^{2^*}),&\text{if }N\geq 3,\\
C_f(|s|^2+|s|^{p}),\quad \text{for some $p>2$},&\text{if }N=2,
\end{cases}
\end{align}
for a.e. $x\in\Omega$ and all $s\in\R$,
where
$F(x,s)=\int_0^s f(x,t)\,dt$, $2^*={2N}/{(N-2)}$ is the
critical Sobolev exponent, $C_f>0$ is a constant independent of
$x\in\Omega$ and $s\in\R$, $\nabla_x F$ denotes the gradient of $F$
with respect to the $x$ variable, and $f'_s(x,s)=\frac{\partial
  f}{\partial s}(x,s)$.

Let $\mu_1(V)$ be the first eigenvalue of the operator
${\mathcal L}_{V}:=-\Delta_{\mathbb S^{N-1}}- V$ on the spherical cap
$C\subset\mathbb S^{N-1}$ under null Dirichlet boundary conditions.
  By classical spectral theory, the spectrum of the  operator
$L_{V}$ is discrete and consists in a nondecreasing diverging sequence of
eigenvalues
\[
\mu_1(V)\leq\mu_2(V)\leq\cdots\leq\mu_k(V)\leq\cdots
\]
with finite multiplicity
the first of which admits the variational characterization
\begin{equation}\label{eq:firsteig}
  \mu_1(V)=\min_{\psi\in H^1_0(C)\setminus\{0\}}\frac{\int_{C}
    \big[\big|\nabla_{\mathbb S^{N-1}}\psi(\theta)\big|^2
    - V(\theta)|\psi(\theta)|^2\big]\,d\sigma(\theta)}{\int_{C}
    |\psi(\theta)|^2\,d\sigma(\theta)}.
\end{equation}
Moreover $\mu_1(V)$ is simple and its associated  eigenfunctions
do not change sign in $C$.

The main result of the present paper provides an evaluation of
 the behavior at the corner $0$ of weak solutions
$u\in H^1(\Omega)$ to
\begin{equation}\label{eq:equation}
\begin{cases}
  -\dive(A(x)\nabla u(x))+\bi(x)\cdot\nabla
  u(x)-\dfrac{V\big(\frac{x}{|x|}\big)}{|x|^2}u(x)=
  h(x)u(x)+f(x,u(x)),&\text{in }\Omega,\\[5pt]
  u=0,&\text{on }\partial\Omega\cap B_R.
\end{cases}
\end{equation}

\begin{Theorem}\label{t:main-u}
  Let $A,\bi,f,h,V$ as in assumptions
  (\ref{eq:matrix1}--\ref{eq:F_assumption2}) and let $\Omega$ satisfying
(\ref{eq:omega}) and
  (\ref{eq:phi1}--\ref{eq:phi4}).  Let $u\in
  H^1(\Omega)\setminus\{0\}$ be a non-trivial weak solution to
  (\ref{eq:equation}).  Then,  there exist $k_0\in
  \N$, $k_0\geq 1$, and
an eigenfunction of the operator ${\mathcal L}_V=-\Delta_{\SN}-V$
  associated to the eigenvalue $\mu_{k_0}(V)$ such that
  $\|\psi\|_{L^2(\SN)}=1$ and
  \begin{equation} \label{convergence-u}
    \lambda^{\frac{N-2}2-\sqrt{\left(\frac{N-2}2\right)^{2}+\mu_{k_0}(V)}}\,
    u(\lambda x)\to
    |x|^{-\frac{N-2}2+\sqrt{\left(\frac{N-2}2\right)^{2}+\mu_{k_0}(V)}}
    \psi\bigg(\frac{x}{|x|}\bigg) \quad \text{as } \lambda\to 0^+
\end{equation}
in $H^1(B_1)$, in $C^{1,\alpha}_{{\rm loc}}(\mathcal C\cap B_1)$
and in $C^{0,\alpha}_{{\rm loc}}(B_1\setminus\{0\})$
for any $\alpha\in (0,1)$, with $u$ being trivially extended
outside~$\Omega$.
\end{Theorem}

As a direct consequence of Theorem \ref{t:main-u}, the following
point-wise upper bound holds.
\begin{Corollary}\label{c:upper_bound}
Under the same assumptions as in Theorem \ref{t:main-u}, let $u\in
  H^1(\Omega)\setminus\{0\}$ be a non-trivial weak solution to
  (\ref{eq:equation}).  Then,  there exists $k_0\in
  \N$, $k_0\geq 1$,  such that
   \begin{equation*}
u(x)=O\Big(|x|^{-\frac{N-2}2+\sqrt{\left(\frac{N-2}2\right)^{2}+\mu_{k_0}(V)}}\Big)
\quad\text{as }|x|\to0^+.
\end{equation*}
\end{Corollary}

A further relevant consequence of our asymptotic analysis is the
following unique continuation principle, whose proof follows straightforwardly
from Theorem \ref{t:main-u}.

\begin{Corollary}\label{c:unique_continuation}
 Under the same assumptions as in Theorem \ref{t:main-u},
let $u\in H^1(\Omega)\setminus\{0\}$ be a weak solution to \eqref{eq:equation}
such that $u(x)=O(|x|^k)$ as $|x|\to 0$, for any $k\in \N$. Then $u\equiv 0$ in
$\Omega$.
\end{Corollary}

Theorem \ref{t:main-u} will be proved by introducing an auxiliary
equivalent problem  obtained as a
diffeomorphic deformation of the original problem (\ref{eq:equation}).
More precisely, letting $C_0$ be as in~(\ref{eq:phi3}), we define
the local diffeomorphism
\begin{equation}\label{eq:def_psi}
\Psi:\R^N\to\R^N,\quad \Psi(y)=\Psi(y',y_N):=(y',y_N+2C_0|y|^{1+\delta}).
\end{equation}
If $u\in H^1(\Omega)$ is a weak solution to (\ref{eq:equation}), then
$w=u\circ \Psi$ weakly solves (in the intersection of
$\Psi^{-1}(\Omega)$ with a sufficiently small neighborhood of $0$)
\begin{equation}\label{eq:36}
  -\dive(\tildA(y)\nabla w(y))+\tildb(y)\cdot\nabla
  w(y)-
\dfrac{V\big(\frac{y}{|y|}\big)}{|y|^2}w(y)
=
  \widetilde h(y)w(y)+\tilde f(y,w(y))
\end{equation}
where
\begin{align}
\label{eq:tildA}  &\tildA(y)=
|\mathop{\rm det}\mathop{\rm Jac}\Psi(y)|
(\mathop{\rm Jac}\Psi(y))^{-1}A(\Psi(y)) ((\mathop{\rm
    Jac}\Psi(y))^T)^{-1}, \\[2pt]
  \label{eq:tildb} &\tildb(y)=|\mathop{\rm det}\mathop{\rm
    Jac}\Psi(y)|
  \bi(\Psi(y))((\mathop{\rm Jac}\Psi(y))^T)^{-1},\\[2pt]
  &\label{eq:tildf}
  \tilde f(y,s)=|\mathop{\rm det}\mathop{\rm Jac}\Psi(y)|f(\Psi(y),s),\\[-5pt]
  & \label{eq:tildh}\widetilde h(y)=|\mathop{\rm det}\mathop{\rm
    Jac}\Psi(y)| h(\Psi(y))+ |\mathop{\rm det}\mathop{\rm Jac}\Psi(y)|
  \Bigg(\dfrac{V\big(\frac{\Psi(y)}{|\Psi(y)|}\big)}{|\Psi(y)|^2}-
  \dfrac{V\big(\frac{y}{|y|}\big)}{|y|^2}\Bigg)\\
\notag&\hskip6cm+ \big( |\mathop{\rm
    det}\mathop{\rm
    Jac}\Psi(y)|-1\big)\dfrac{V\big(\frac{y}{|y|}\big)}{|y|^2}.
\end{align}
For the auxiliary problem (\ref{eq:36}) an Almgren
monotonicity formula is used to describe the rate and the shape of the
singularity of solutions, by relating
them to the eigenvalues and the eigenfunctions of the angular operator
$\mathcal L_V$ on the spherical cap $C$.  The behavior of solutions of
the auxiliary problem \eqref{eq:36}
 (and then of the original one \eqref{eq:92}) near the corner
is indeed classified on the basis of the limit of the following \emph{Almgren
  type frequency} function 
\begin{align}\label{eq:6}
&\mathcal N(r)=\frac{r\int_{\Psi^{-1}(\Omega)\cap B_r} \Big(
\tildA\nabla w\cdot \nabla w+\tildb\cdot\nabla
  w\,w-
\frac{V(\frac{y}{|y|})}{|y|^2}|w|^2-
  \widetilde h w^2-\tilde f(y,w)w\Big)dy}{\int_{\Psi^{-1}(\Omega)\cap \partial B_r}
\frac{\tildA(y)y\cdot y}{|y|^2}
w^2(y)\,d\sigma(y)},
\end{align}
which is defined for $r>0$ sufficiently small (see \eqref{eq:alm_fun}
and \eqref{Almgren}).

\begin{Theorem}\label{t:main-w2}
  Let $A,\bi,f,h,V$ as in assumptions
  (\ref{eq:matrix1}--\ref{eq:F_assumption2}) and let $\Omega$
  satisfying (\ref{eq:omega}) and (\ref{eq:phi1}--\ref{eq:phi4}).  Let
  $u\in H^1(\Omega)\setminus\{0\}$ be a non-trivial weak solution to
  (\ref{eq:equation})
and $w=u\circ\Psi$ with  $\Psi$ as in (\ref{eq:def_psi}). 
Letting  $\mathcal N$
  as in (\ref{eq:6}),  there exists $k_0\in \N$, $k_0\geq 1$, such that
\begin{equation*}
  \lim_{r\to 0^+} \mathcal
  N(r)=-\frac{N-2}2+\sqrt{\left(\frac{N-2}2\right)^{\!\!2}+\mu_{k_0}(V)}.
\end{equation*}
Furthermore there exists $\psi\in H^1_0(C)\subset H^1(\SN)$
  eigenfunction of the operator ${\mathcal L}_V=-\Delta_{\SN}-V$
  associated to the eigenvalue $\mu_{k_0}(V)$ such that
 \begin{equation}\label{eq:8}
\lambda^{\frac{N-2}2-\sqrt{\left(\frac{N-2}2\right)^{\!2}+\mu_{k_0}(V)}} w(\lambda x)\to
  |x|^{-\frac{N-2}2+\sqrt{\left(\frac{N-2}2\right)^{\!2}+\mu_{k_0}(V)}}\psi
\bigg(\frac{x}{|x|}\bigg) \quad \text{as } \lambda\to 0^+
\end{equation}
in $H^1(B_1)$ and  in $C^{1,\alpha}_{{\rm loc}}(\mathcal
C\cap B_1)$ for any $\alpha\in (0,1)$.
\end{Theorem}

Furthermore, Theorem \ref{t:main-w} will provide more precise
informations on the limit angular profile~$\psi$: if $m\geq 1$ is the
multiplicity of the eigenvalue $\mu_{k_0}(V)$ and $\{\psi_i:j_0\leq
i\leq j_0+m-1\}$ is an $L^2(C)$-orthonormal basis for the eigenspace
associated to $\mu_{k_0}(V)$, then the eigenfunction $\psi$ in
(\ref{eq:8}) (which coincides with the one appearing in
(\ref{convergence-u}), as clarified in the proof of Theorem
\ref{t:main-u}, see section \ref{sec:straightening-domain})
can be written as
$$
\psi(\theta)=
\sum_{i=j_0}^{j_0+m-1}
    \beta_i \psi_i(\theta),
$$
where the coefficients $\beta_i$ can be represented in terms of the
\emph{Cauchy's integral type formula} (\ref{eq:88}).

We
emphasize that our monotonicity approach allows excluding the presence
of logarithmic factors in the leading term of the asymptotic
expansion; we refer to \cite{FFT3} for a detailed comparison between
the monotonicity approach to asymptotic analysis and the results
obtained in earlier literature (see
e.g. \cite{cd,dauge,Grisvard,ko,lehman,MNP1,wasow,wigley}) by integral
representation and Mellin transform methods.

In section \ref{sec:examples} we produce an example in dimension $N=2$
of a harmonic function on a domain with a corner of any amplitude and
delimited by arcs violating assumptions (\ref{eq:phi2}--\ref{eq:63}),
satisfying null Dirichlet boundary conditions but exhibiting dominant
logarithmic terms in its asymptotic expansion. Hence assumptions
(\ref{eq:phi2}--\ref{eq:63}) are crucial for excluding the presence of
logarithms, even under null boundary conditions.
Besides the failure of conditions (\ref{eq:phi2}--\ref{eq:63}), other
possible reasons of occurring of logarithms in the expansion could be
boundary conditions (even if very regular when the amplitude is
resonant, see \cite{lehman,wasow}) or lack of linearity with respect
to the first derivatives of $u$, see \cite{wigley}.

\medskip
\noindent
{\bf Notation. } We list below some notation used throughout the
paper.\par
\begin{itemize}
\item[-] For all $r>0$, $B_r$ denotes the ball $\{x\in\R^N: |x|<r\}$
  in $\R^N$ with center at $0$ and radius $r$.
\item[-] ${\mathcal M}_{N\times N}$ denotes the space of $N\times N$
  real matrices.
\item[-] ${\rm Id}_N$ denotes the identity $N\times N$ matrix .
\item[-] For every vector field $\Psi\in C^1(\R^N,\R^N)$, $\mathop{\rm Jac}\Psi$
denotes the Jacobian matrix.

\end{itemize}

\section{An equivalent problem}\label{sec:an-equiv-probl}

In this section we construct  an auxiliary
equivalent problem by a
diffeomorphic deformation of the domain.
\begin{Lemma}\label{l:phi3}
Under assumptions (\ref{eq:phi1}--\ref{eq:63}), there exists $C_0>0$ such that
(\ref{eq:phi3}) holds.
\end{Lemma}
\begin{pf}
From (\ref{eq:phi2}) and (\ref{eq:63}), we can estimate, for $N\geq3$,
\begin{align*}
  & |\varphi(x')-\nabla\varphi(x')\cdot x'|\leq
  {\textstyle{|x'|\Big|\frac{\varphi(x')}{|x'|}-g\big(\frac{x'}{|x'|}\big)\Big|
      +|x'|\Big|\Big(\nabla\varphi(x')-g\big(\frac{x'}{|x'|}\big)
      \frac{x'}{|x'|}\Big)\cdot
      \frac{x'}{|x'|}\Big|}}\\
  &={\textstyle{|x'|\Big|\frac{\varphi(x')}{|x'|}-g
      \big(\frac{x'}{|x'|}\big)\Big|
      +|x'|\Big|\Big(\nabla\varphi(x')-g\big(\frac{x'}{|x'|}\big)
      \frac{x'}{|x'|}-\nabla_{{\mathbb
          S}^{N-2}}g\big(\frac{x'}{|x'|}\big)\Big)\cdot
      \frac{x'}{|x'|}\Big|}}=O(|x'|^{1+\delta})
\end{align*}
as $|x'|\to 0^+$ thus proving (\ref{eq:phi3}). The proof for $N=2$ is similar.
\end{pf}

\noindent We notice that the function $\Psi$ defined in (\ref{eq:def_psi}) 
satisfies $\Psi\in C^1(\R^N,\R^N)$,
\begin{equation*}
\mathop{\rm Jac}\Psi(y',y_N)=
\begin{pmatrix}
  1&0&\cdots&0&0\\
  0&1&\cdots&0&0\\
  \vdots&\vdots&&\vdots&\vdots\\
  0&0&\cdots&1&0\\[3pt]
  2C_0\frac{1+\delta}{|y|^{1-\delta}}y_1&
2C_0\frac{1+\delta}{|y|^{1-\delta}}y_2&\cdots&
  2C_0\frac{1+\delta}{|y|^{1-\delta}}y_{N-1}&1+2C_0\frac{1+\delta}{|y|^{1-\delta}}y_{N}
\end{pmatrix}
\end{equation*}
for all $(y',y_N)\neq 0$, and $\mathop{\rm Jac}\Psi(0)={\rm
  Id}_N$. Hence there exists a bounded neighborhood $U\subset\R^N$ of
$0$ such that the restriction $\Psi\big|_U:U\to\Psi(U)$ is a
$C^1$-diffeomorphism.  Let us denote as
\begin{equation}\label{eq:tildeomega}
\tildO:=\Psi^{-1}(\Omega\cap \Psi(U))
\end{equation}
and let us consider the function
\begin{equation}\label{eq:21}
\widetilde\varphi(y')=\Psi^{-1}(y',\varphi(y'))\cdot e_N,\quad
e_N=(0,0,\cdots,0,1),
\end{equation}
which is well defined in a sufficiently small neighborhood of $0$ in
$\R^{N-1}$.
\begin{Lemma}\label{l:profile}
There exists $\widetilde R>0$ such that
\begin{equation}\label{eq:1}
  \widetilde\varphi(y')+2C_0\big(|y'|^2+
  |\widetilde\varphi(y')|^2\big)^{\frac{1+\delta}2}=\varphi(y')\quad
\text{for all }y'\in\R^{N-1}, \ |y'|<\widetilde R,
\end{equation}
and
\begin{equation}\label{eq:9}
  \tildO\cap B_{\widetilde R}=
  \{(y',y_N)\in B_{\widetilde R}:y_N>\widetilde \varphi(y')\}.
\end{equation}
\end{Lemma}
\begin{pf}
From the definition of $\widetilde \varphi$ we have that
$$
\Psi^{-1}(y', \varphi(y'))=(y',\widetilde\varphi(y')),\quad\text{i.e.}\quad
(y', \varphi(y'))=\Psi(y',\widetilde\varphi(y'))
$$
for all $y'\in\R^{N-1}$ such that $(y', \varphi(y'))\in\Psi(U)$, which
implies (\ref{eq:1}) for some $\widetilde R>0$ sufficiently small. To
prove (\ref{eq:9}) we observe that there exists $R_0>0$ such that for
every fixed $x'\in\R^{N-1},|x'|<R_0$, the function
$$
t\in(-R_0,R_0)\mapsto  \Psi^{-1}(x',t)\cdot e_N
$$
is strictly increasing with respect to $t$, since its derivative
$$
\frac{d}{dt}\Big(\Psi^{-1}(x',t)\cdot
e_N\Big)=\frac1{1+2C_0(1+\delta)|\Psi^{-1}(x',t)|^{-1+\delta}
\Psi^{-1}(x',t)\cdot e_N
}
$$
is strictly positive provided $R_0$ is sufficiently small. In
particular, letting $x=\Psi(y)$ in a sufficiently small neighborhood
of $0$, $x_N>\varphi(x')$ if and only if $\Psi^{-1}(x',x_N)\cdot
e_N>\Psi^{-1}(x',\varphi(x'))\cdot e_N$ and hence if and only if
$y_N>\widetilde\varphi(y')$, which, in view of (\ref{eq:omega}) yields
the conclusion.
\end{pf}

\begin{remark}\label{r:limtildephi}
From assumption (\ref{eq:phi2}) and (\ref{eq:1}), it follows that
\begin{equation}\label{eq:23}
  \sup_{\nu\in {\mathbb S}^{N-2}}\Big|
  \frac{\widetilde\varphi(t\nu)}{t}-g(\nu)\Big|=O(t^\delta)
  \quad\text{as }t\to0^+,
\end{equation}
which implies
\begin{equation}\label{eq:68}
|\widetilde \varphi(y')-\varphi_0(y')|
=O(|y'|^{1+\delta})\quad\text{as }|y'|\to0^+.
\end{equation}
Furthermore, from assumption (\ref{eq:63}) and (\ref{eq:1}), there also holds
\begin{equation}\label{eq:69}
|\nabla\widetilde \varphi(y')-\nabla\varphi_0(y')|
=O(|y'|^{\delta})\quad\text{as }|y'|\to0^+,
\end{equation}
whereas assumption
(\ref{eq:phi4}) implies that
\begin{equation}\label{eq:58}
|D^2\widetilde \varphi(y')|=O(|y'|^{-1})\quad\text{as }|y'|\to 0.
\end{equation}

\end{remark}

If $u\in H^1(\Omega)$ is a weak solution to (\ref{eq:equation}), then
$w=u\circ \Psi\in H^1(\tildO)$ is, up to shrinking $\widetilde R>0$, a
weak solution to
\begin{equation}\label{eq:tildequation}
\begin{cases}
  -\dive(\tildA(y)\nabla w(y))+\tildb(y)\cdot\nabla
  w(y)-
\dfrac{V\big(\frac{y}{|y|}\big)}{|y|^2}w(y)
=
  \widetilde h(y)w(y)+\tilde f(y,w(y)),&\text{in }\tildO,\\[5pt]
  w=0,&\hskip-1.2cm\text{on }\partial\tildO\cap B_{\widetilde R},
\end{cases}
\end{equation}
where $\tildA,\tildb,\widetilde h, \tilde f$ are as in
(\ref{eq:tildA}--\ref{eq:tildh}).

\begin{Lemma}\label{l:new_coeff}
  Let $A,\bi,\Psi,f,h,V$ as in assumptions
  (\ref{eq:matrix1}--\ref{eq:F_assumption2}), (\ref{eq:def_psi}), and
  $\tildA,\tildb,\tilde f,\widetilde h$ as
  (\ref{eq:tildA}--\ref{eq:tildh}). Then
\begin{align}
  \label{eq:2}&
\left\{\!\!
\begin{array}{l}
  \tildA\in C^0(\overline{\tildO},{\mathcal M}_{N\times N}),\quad
  \tildA\in W^{1,\infty}_{\rm loc}(\overline{\tildO}\setminus\{0\},
{\mathcal M}_{N\times N}),\quad
  \tildA(0)={\rm Id}_N,\quad (\tildA)_{ij}=(\tildA)_{ji},\\[5pt]
  \tildA(y)\xi\cdot\xi
  \geq C_{\tildA}|\xi|^2\text{ for all }\xi\in\R^N,
  y\in\tildO,\text{ and some }C_{\tildA}>0,\\[5pt]
  \|d\tildA(y)\|_{{\mathcal
      L}(\R^N,{\mathcal M}_{N\times N})}=O(|y|^{-1+\delta})\quad\text{and}\quad
  \tildA(y)-{\rm Id}_N=O(|y|^{\delta})\quad\text{as }|y|\to0,
\end{array}
\right.\\
\label{eq:3}&\tildb\in L^\infty_{\rm loc}(\tildO,\R^N), \quad
|\tildb(y)|=O(|y|^{-1+\delta})\quad\text{as }|y|\to0,\\
\label{eq:5_1}
&\tilde f\in C^0(\tildO\times \R),\quad \widetilde F\in C^1(\tildO\times \R),
\quad s\mapsto \tilde f(y,s)\in C^1(\R)\text{ for a.e. }y\in\tildO, \\
\label{eq:5_2}
&|\tilde f(y,s)s|+|\tilde f'_s(y,s)s^2|+|\nabla_y \widetilde F(y,s)||y|
\leq
\begin{cases}
C_{\tilde f}(|s|^2+|s|^{2^*}),&\text{if }N\geq 3,\\
C_{\tilde f}(|s|^2+|s|^{p}),&\text{if }N=2,
\end{cases}\\
\label{eq:4}&\widetilde h \in L^\infty_{\rm loc}(\tildO),\quad
\widetilde h(y)=O(|y|^{-2+\delta})\quad\text{as }|y|\to0,
\end{align}
where $\widetilde F(y,s)=\int_0^s \tilde f(y,t)\,dt=
|\mathop{\rm det}\mathop{\rm Jac}\Psi(y)|F(\Psi(y),s)$.
\end{Lemma}
\begin{pf}
  Estimates (\ref{eq:2}--\ref{eq:5_2}) follow  from
  (\ref{eq:tildA}--\ref{eq:tildf}), (\ref{eq:6}--\ref{eq:8}), and
  assumptions (\ref{eq:matrix1}--\ref{eq:b_assumption}),
  (\ref{eq:F_assumption1}--\ref{eq:F_assumption2}). To prove estimate
  (\ref{eq:4}), we first observe that (\ref{eq:h_assumption}) implies
  $|\mathop{\rm det}\mathop{\rm Jac}\Psi(y)|
  h(\Psi(y))=O(|y|^{-2+\delta})$ as $|y|\to0$.
From  (\ref{eq:V_assumption}) and
$$
|\Psi(y)|=|y|(1+O(|y|^{\delta})), \quad
\Psi(y)=y+O(|y|^{1+\delta})\quad \text{as }|y|\to0,
$$
it follows that
\begin{align*}
  &\Bigg|\dfrac{V\big(\frac{\Psi(y)}{|\Psi(y)|}\big)}{|\Psi(y)|^2}-
  \dfrac{V\big(\frac{y}{|y|}\big)}{|y|^2}\Bigg|\leq
  \Bigg|\dfrac{V\big(\frac{\Psi(y)}{|\Psi(y)|}\big)-V\big(\frac{y}{|y|}\big)}
  {|\Psi(y)|^2}\Bigg|+ \Bigg|V\big({\textstyle{\frac{y}{|y|}}}\big)
  \bigg(\frac1{|\Psi(y)|^2}-\frac1{|y|^2}\bigg)\Bigg|\\
  &\leq
  \frac{L_V}{|y|^2(1+O(|y|^{\delta}))}\bigg|\frac{\Psi(y)}{|y|(1+O(|y|^{\delta}))}-
  \frac{y}{|y|} \bigg|+\|V\|_{L^{\infty}(\SN)}
  \bigg|\frac{1}{|y|^2(1+O(|y|^{\delta}))}-
  \frac{1}{|y|^2} \bigg|\\
  &= \frac{L_V}{|y|^2(1+O(|y|^{\delta}))}
  \frac{|\Psi(y)-y+O(|y|^{1+\delta})|}{|y|(1+O(|y|^{\delta}))}+
  \|V\|_{L^{\infty}(\SN)}\frac{O(|y|^{\delta})}{|y|^2(1+O(|y|^{\delta}))}
  =O(|y|^{-2+\delta}),
\end{align*}
which, taking into account that $|\mathop{\rm det}\mathop{\rm
  Jac}\Psi(y)|=1+O(|y|^\delta)$ as $|y|\to0$, yields (\ref{eq:4}).
\end{pf}

\begin{Lemma}\label{l:Ay}
  Let $A$ as in assumptions
  (\ref{eq:matrix1}--\ref{eq:matrix2}) and
  $\tildA$ as in
  (\ref{eq:tildA}). Then
$$
(\tildA(y)y)_i=
\begin{cases}
  y_i+O(|y|^{2}),&\text{if }1\leq i\leq N-1,\\
  \dfrac{y_N-2C_0(1+\delta)|y|^{-1+\delta}|y'|^2}
  {1+2C_0(1+\delta)|y|^{-1+\delta}y_N} +O(|y|^{2}),&\text{if }i=N,
\end{cases}
$$
as $y\to 0$.
\end{Lemma}
\begin{pf}
The proof follows from (\ref{eq:6}--\ref{eq:8}), direct calculations and
the estimate
$$
|a_{ij}(\Psi(y))-\delta_{ij}|=O(|y|)\quad\text{as }|y|\to0,
$$
which is a consequence of (\ref{eq:matrix1}) and  (\ref{eq:matrix2}).
\end{pf}

Let us consider the exterior unit normal $\tilde\nu$ to
$\partial\tildO\cap B_{\widetilde R}$. From (\ref{eq:1}) and
(\ref{eq:9}), it follows that
\begin{align}\label{eq:normaletilde}
  \tilde\nu(y',y_N)&=\frac{(\nabla\widetilde \varphi(y'),-1)}
  {\sqrt{|\nabla\widetilde \varphi(y')|^2+1}}
  =\frac1{\sqrt{|\nabla\widetilde \varphi(y')|^2+1}}\bigg(
  \dfrac{\nabla\varphi(y')-2C_0(1+\delta)|y|^{-1+\delta}y'}
  {1+2C_0(1+\delta)|y|^{-1+\delta}y_N},-1\bigg).
\end{align}
\begin{Lemma}\label{l:tildAytildnu}
  Let $A$ as in assumptions
  (\ref{eq:matrix1}--\ref{eq:matrix2}),
  $\tildA$ as in
  (\ref{eq:tildA}),
$\tildO$ as in (\ref{eq:tildeomega}) with $\Omega$ satisfying (\ref{eq:omega}),
(\ref{eq:phi1}--\ref{eq:63}), and
 $\tilde \nu$ as in (\ref{eq:normaletilde}). Then
$\tildA(y)y\cdot\tilde\nu(y)\geq 0$ for all
$y\in (\partial\tildO\cap B_{r})\setminus\{0\}$
provided $r$ is sufficiently small.
\end{Lemma}
\begin{pf}
  Taking into account that $y_N=\widetilde\varphi(y')$ and $|y'|^2+
  |\widetilde\varphi(y')|^2=|y|^2$ on $\partial\tildO\cap
  B_{\widetilde R}$, from Lemma \ref{l:Ay}, (\ref{eq:normaletilde}),
  and (\ref{eq:1}), we deduce that
\begin{align*}
  (\sqrt{|\nabla\widetilde
    \varphi(y')|^2+1})\tildA(y)y\cdot\tilde\nu(y)&= \frac{y'\cdot
    \nabla\varphi(y')-y_N}{1+2C_0(1+\delta)|y|^{-1+\delta}y_N}
  +O(|y|^{2})\\
  &=\frac{y'\cdot
    \nabla\varphi(y')-\varphi(y')+2C_0|y|^{1+\delta}}{1+2C_0(1+\delta)|y|^{-1+\delta}y_N}
  +O(|y|^{2}).
\end{align*}
Hence Lemma \ref{l:phi3} yields
$$
  (\sqrt{|\nabla\widetilde
    \varphi(y')|^2+1})\tildA(y)y\cdot\tilde\nu(y)\geq
\frac{C_0|y|^{1+\delta}}{1+2C_0(1+\delta)|y|^{-1+\delta}y_N}
  +O(|y|^{2})\geq 0
$$
provided $|y|$ is sufficiently small.
\end{pf}

\noindent The above lemma ensures that, under assumptions
(\ref{eq:phi1}--\ref{eq:63}),   (\ref{eq:omega}--\ref{eq:matrix2}),
(\ref{eq:tildeomega}), and    (\ref{eq:tildA}),
up to shrinking $\widetilde R>0$ there holds
 \begin{equation}\label{eq:pos_bound}
\tildA(y)y\cdot\tilde\nu(y)\geq 0\quad \text{for all
}y\in (\partial\tildO\cap B_{\widetilde R})\setminus\{0\}.
\end{equation}

\section{Hardy type inequalities ($N\geq 3$)}\label{sec:hardy-type-ineq}

Throughout this section we assume $N\geq 3$.
The following lemma establishes the relation between the values
$\Lambda(V)$ defined in (\ref{eq:LambdaV}) and $\mu_1(V)$ defined in
(\ref{eq:firsteig}) and the positivity of the quadratic form
associated with the principal part of the elliptic operator on the
limit domain $\mathcal C$ defined in (\ref{eq:cilindro}).
\begin{Lemma}\label{l:pos}
  If $N\geq 3$ and $V\in L^{\infty}(\SN)$, then the following
  conditions are equivalent:
\begin{align*} {\rm i)}\quad&\inf_{u\in {\mathcal D}^{1,2}(\mathcal C)\setminus\{0\}}
\dfrac{\int_{\mathcal C} |\nabla  v(x)|^2dx-
\int_{\mathcal C}\frac{V(x/|x|)}{|x|^2}\,v^2(x)\,dx}
  {\int_{\mathcal C}|\nabla v(x)|^2\,dx}>0; \\[7pt]
  {\rm ii)} \quad&\Lambda(V)<1;\\
  {\rm iii)}\quad&
  \mu_1(V)>-\big({\textstyle{\frac{N-2}2}}\big)^{\!2}.
\end{align*}
\end{Lemma}
\begin{pf}
The equivalence between i) and ii) follows from the definition of
  $\Lambda(V)$, see (\ref{eq:LambdaV}). The equivalence between
  i) and iii) can be proved arguing as in
  \cite[Proposition 1.3 and Lemma 1.1]{terracini96}.\end{pf}

Let $\tildO$ be as in (\ref{eq:tildeomega}) with $\Omega$ satisfying
(\ref{eq:omega}) and (\ref{eq:phi1}--\ref{eq:63}), and
$\widetilde\varphi$ be as (\ref{eq:21}).  For every $r\in(0,\tilde R)$
let us denote
\begin{align}\label{eq:Cr}
  C_r=\SN\cap
  \big({\textstyle{\frac1r}}{\tildO}\big)=\{(y',y_N)\in\SN:
  y_N>r^{-1}\widetilde\varphi(ry')\}
\end{align}
and, for $V\in L^{\infty}(\SN)$,
 let us consider the first eigenvalue $\mu_1(V,r)$ of the operator
$-\Delta_{\mathbb S^{N-1}}- V$ on the spherical cap
$C_r$ under null Dirichlet boundary conditions, i.e.
\begin{equation}\label{eq:firsteig_r}
  \mu_1(V,r)=\min_{\psi\in H^1_0(C_r)\setminus\{0\}}\frac{\int_{C_r}
    \big[\big|\nabla_{\mathbb S^{N-1}}\psi(\theta)\big|^2
    - V(\theta)|\psi(\theta)|^2\big]\,d\sigma(\theta)}{\int_{C_r}
    |\psi(\theta)|^2\,d\sigma(\theta)}.
\end{equation}
We also define
\begin{equation}\label{eq:Lambda_r}
  \Lambda(V,r)=\max_{\psi\in H^1_0(C_r)\setminus\{0\}}\frac{\int_{C_r}
V(\theta)|\psi(\theta)|^2\,d\sigma(\theta)}
{\int_{C_r}    \big[\big|\nabla_{\mathbb S^{N-1}}\psi(\theta)\big|^2+
\big(\frac{N-2}{2}\big)^2|\psi(\theta)|^2\big]\,d\sigma(\theta)}.
\end{equation}

\begin{Lemma}\label{l:con_eig}
  Let  $N\geq 3$, $V\in L^{\infty}(\SN)$,
$\mu_1(V,r)$ be defined in (\ref{eq:firsteig_r}), $\mu_1(V)$ in
  (\ref{eq:firsteig}), $\Lambda(V,r)$ in (\ref{eq:Lambda_r}), and
  $\Lambda(V)$ in (\ref{eq:LambdaV}). Then
\begin{equation}\label{eq:22}
\lim_{r\to0^+} \mu_1(V,r)=\mu_1(V)
\end{equation}
and
\begin{equation}\label{eq:28}
\lim_{r\to0^+} \Lambda(V,r)=\Lambda(V).
\end{equation}
\end{Lemma}
\begin{pf}
We first claim that
\begin{equation}\label{eq:claim1}
\text{for every }\psi\in C^{\infty}_{\rm c}(C)
\text{ there exists }r_0>0\text{ such that }\mathop{\rm supp} \psi\subseteq C_r
\text{ for all }r\in(0,r_0).
\end{equation}
To prove the claim, let us consider $\psi\in C^{\infty}_{\rm c}(C)$ and
denote $K=\mathop{\rm supp} \psi$. Since $K$ is compact, we have that
$$
\delta=\min_{(y',y_N)\in K}(y_N-\varphi_0(y'))>0.
$$
From (\ref{eq:23}), there exists $r_0$ such that
\begin{equation*}
\Big|
  \frac{\widetilde\varphi(t\nu)}{t}-g(\nu)\Big|<\delta
\quad\text{for all $t\in
(0,r_0)$ and for all $\nu\in {\mathbb S}^{N-2}$}.
\end{equation*}
Then for all $r\in (0,r_0)$ and $(y',y_N)\in K$  we have that
$$
\bigg|\varphi_0(y')-  \frac{\widetilde\varphi(ry')}{r}\bigg|
=|y'|\bigg|
  \frac{\widetilde\varphi(ry')}{r|y'|}-g\Big(\frac{y'}{|y'|}\Big)\Big|<\delta
$$
and hence
$$
y_N-\frac{\widetilde\varphi(ry')}{r}\geq
(y_N-\varphi_0(y'))-\bigg|\varphi_0(y')-  \frac{\widetilde\varphi(ry')}{r}\bigg|
\geq \delta-\bigg|\varphi_0(y')-  \frac{\widetilde\varphi(ry')}{r}\bigg|>0
$$
which implies that $K\subseteq C_r$ for all $r\in(0,r_0)$, thus
proving claim (\ref{eq:claim1}).

From (\ref{eq:claim1}) it follows that for every
$\psi\in C^{\infty}_{\rm c}(C)$
 there exists $r_0>0$ such that, for all $r\in(0,r_0)$,
$$
\mu_1(V,r)\leq
\frac{\int_{C_r}
    \big[\big|\nabla_{\mathbb S^{N-1}}\psi(\theta)\big|^2
    - V(\theta)|\psi(\theta)|^2\big]\,d\sigma(\theta)}{\int_{C_r}
    |\psi(\theta)|^2\,d\sigma(\theta)}=\frac{\int_{C}
    \big[\big|\nabla_{\mathbb S^{N-1}}\psi(\theta)\big|^2
    - V(\theta)|\psi(\theta)|^2\big]\,d\sigma(\theta)}{\int_{C}
    |\psi(\theta)|^2\,d\sigma(\theta)}.
$$
Hence
$$
\limsup_{r\to0^+}\mu_1(V,r)\leq \frac{\int_{C}
    \big[\big|\nabla_{\mathbb S^{N-1}}\psi(\theta)\big|^2
    - V(\theta)|\psi(\theta)|^2\big]\,d\sigma(\theta)}{\int_{C}
    |\psi(\theta)|^2\,d\sigma(\theta)}
$$
for all $\psi\in C^{\infty}_{\rm c}(C)$. By density of
$C^{\infty}_{\rm c}(C)$ in $H^1_0(C)$, we conclude that
$$
\limsup_{r\to0^+}\mu_1(V,r)\leq \mu_1(V).
$$
To prove (\ref{eq:22}), it remains to show that
\begin{equation}\label{eq:24}
\liminf_{r\to0^+}\mu_1(V,r)\geq \mu_1(V).
\end{equation}
Arguing by contradiction, let us assume that (\ref{eq:24}) fails, then
there exists  $\{r_n\}_{n\in\N}\subset(0,\widetilde R)$ such
that $\lim_{n\to+\infty}r_n=0$ and
$\lim_{n\to+\infty}\mu_1(V,r_n)<\mu_1(V)$. For all $n$, let $\psi_n\in
H^1_0(C_{r_n})$ such that
$$
  \mu_1(V,r_n)=\int_{C_{r_n}}
    \big[\big|\nabla_{\mathbb S^{N-1}}\psi_n(\theta)\big|^2
    - V(\theta)|\psi_n(\theta)|^2\big]\,d\sigma(\theta)\quad\text{and}\quad
\int_{C_{r_n}}
    |\psi_n(\theta)|^2\,d\sigma(\theta)=1.
$$
Let us identify $\psi_n$ with its trivial extension in $\SN$ which belongs
to $H^1(\SN)$.  It is easy to verify that $\{\psi_n\}_{n\in\N}$ is bounded in
$H^1(\SN)$ so that there exists a subsequence $\psi_{n_k}$ weakly and a.e.
converging to some $\psi$ in $H^1(\SN)$. By compactness of the
embedding $H^1(\SN)\hookrightarrow L^2(\SN)$, we have that $\int_{\SN}
\psi^2=1$ and by weakly lower semicontinuity
\begin{multline}\label{eq:25}
\int_{\SN}
    \big[\big|\nabla_{\mathbb S^{N-1}}\psi(\theta)\big|^2
    - V(\theta)|\psi(\theta)|^2\big]\,d\sigma(\theta)\\
\leq\liminf_{k\to\infty}
\int_{\SN}
    \big[\big|\nabla_{\mathbb S^{N-1}}\psi_{n_k}(\theta)\big|^2
    - V(\theta)|\psi_{n_k}(\theta)|^2\big]\,d\sigma(\theta)
=\liminf_{k\to\infty}\mu_1(V,r_{n_k})<\mu_1(V).
  \end{multline}
  By a.e. convergence of $\psi_{n_k}$ to $\psi$, it is easy to verify
  that $\psi\in H^1_0(C)$ thus implying that
$$
\mu_1(V)\leq \int_{\SN} \big[\big|\nabla_{\mathbb
  S^{N-1}}\psi(\theta)\big|^2 -
V(\theta)|\psi(\theta)|^2\big]\,d\sigma(\theta)
$$ giving rise to a contradiction with (\ref{eq:25}).
(\ref{eq:22}) is thereby proved.  The proof of (\ref{eq:28}) can be
derived in similar way after observing that
\begin{equation*}
\Lambda(V)=  \max_{\psi\in H^1_0(C)\setminus\{0\}}\frac{\int_{C}
V(\theta)|\psi(\theta)|^2\,d\sigma(\theta)}
{\int_{C}    \big[\big|\nabla_{\mathbb S^{N-1}}\psi(\theta)\big|^2+
\big(\frac{N-2}{2}\big)^2|\psi(\theta)|^2\big]\,d\sigma(\theta)},
\end{equation*}
see
\cite[Lemma 1.1]{terracini96}.
\end{pf}

We extend to singular potentials on corner sets  the
Hardy type inequality with boundary terms proved by Wang and Zhu
in \cite{wz}.
 For every $r\in(0,\tilde R)$ let us denote
\begin{align}\label{eq:omerar}
  \Omega_r=\tildO\cap B_r,\quad S_r=(\partial B_r)\cap\tildO,\quad
  \Gamma_r=(\partial\tildO)\cap B_r,
\end{align}
so that $\partial\Omega_r=\overline{S_r}\cup\Gamma_r$.

\begin{Lemma}\label{l:hardyboundary}
  Let $N\geq 3$ and $V\in L^{\infty}(\SN)$. For every $r\in
  (0,\widetilde R)$ and $v\in H^1(\Omega_r)$ such that $v=0$ on
  $\Gamma_r$, the following inequality holds
  \begin{multline}\label{eq:26}
    \int_{\Omega_r} \bigg( |\nabla
    v(y)|^2-\frac{V(\frac{y}{|y|})}{|y|^2} v^2(y)\bigg)\,dy+
    \frac{N-2}{2r}\int_{S_r}v^2(y)\,d\sigma\\
    \geq
    \int_{\Omega_r}\bigg(\bigg(\frac{N-2}{2}\bigg)^{\!\!2}+\mu_1(V,|y|)\bigg)
    \frac{v^2(y)}{|y|^2}\,dy.
  \end{multline}
  \end{Lemma}
\begin{pf}
  Let $v\in C^\infty_{\rm c}(\tildO\cap \overline{B}_r)$ for some
  $r\in(0,\widetilde R)$.  Passing to polar coordinates and denoting
  as $\tilde v$ the trivial extension of $v$ in $\overline{B}_r$, we
  have that $\tilde v\in C^\infty(\overline{B}_r)$ and
\begin{align}\label{eq:27}
  \int_{\Omega_r} \bigg( |\nabla
  v(y)|^2&-\frac{V(\frac{y}{|y|})}{|y|^2} v^2(y)\bigg)\,dy+
  \frac{N-2}{2r}\int_{S_r}v^2(y)\,d\sigma\\
  \notag=&
  \int_0^r\bigg(s^{N-1}\int_{C_s}|\partial_sv(s,\theta)|^2\,d\sigma(\theta)\bigg)ds
  +\frac{N-2}{2r}r^{N-1}\int_{C_r}v^2(r,\theta)\,d\sigma(\theta)\\
  \notag\quad&+\int_0^r\bigg(s^{N-3}\int_{C_s} \left[|\nabla_{{\mathbb
        S}^{N-1}}v(s,\theta)|^2-V(\theta)
    v^2(s,\theta)\right]\,d\sigma(\theta)\bigg)\,ds\\
  &\notag= \int_{{\mathbb S}^{N-1}}\bigg(
  \int_0^{r}s^{N-1}|\partial_s\tilde
  v(s,\theta)|^2\,ds\bigg)\,d\sigma(\theta)
  +\frac{N-2}{2r}r^{N-1}\int_{{\mathbb S}^{N-1}}\tilde v^2(r,\theta)
  \,d\sigma(\theta)\\
  \notag& \quad +\int_0^r\bigg(s^{N-3}\int_{C_s}
  \left[|\nabla_{{\mathbb S}^{N-1}}v(s,\theta)|^2-V(\theta)
    v^2(s,\theta)\right]\,d\sigma(\theta)\bigg)\,ds.
\end{align}
For all $\theta\in{\mathbb S}^{N-1}$, let $\varphi_{\theta}\in
C^\infty(0,r)$ be defined by $\varphi_{\theta}(r)=\tilde v(r,\theta)$, and
$\widetilde\varphi_{\theta}\in C^\infty(B_r)$ be the radially
symmetric function given by
$\widetilde\varphi_{\theta}(x)=\varphi_{\theta}(|x|)$.  We notice that
$0\not\in \mathop{\rm supp}\widetilde\varphi_{\theta}$.
The Hardy
inequality with boundary term proved in \cite{wz} yields
\begin{align}\label{eq:5ang}
\int_{{\mathbb S}^{N-1}}&\bigg(
  \int_0^{r}s^{N-1}|\partial_s\tilde v(s,\theta)|^2\,ds\bigg)\,d\sigma(\theta)
  +\frac{N-2}{2r}r^{N-1}\int_{{\mathbb S}^{N-1}}\tilde v^2(r,\theta)\,d\sigma(\theta)
\\
  \notag& =\frac1{\omega_{N-1}}  \int_{{\mathbb S}^{N-1}}\bigg(
\int_{B_r}|\nabla\widetilde \varphi_{\theta}(x)|^2\,dx
+\frac{N-2}{2r}\int_{\partial B_r}|\widetilde \varphi_{\theta}(x)|^2\,d\sigma
\bigg)\,d\sigma(\theta)
\\&\notag\geq
\frac1{\omega_{N-1}}  \bigg(\frac{N-2}2\bigg)^{\!\!2}
  \int_{{\mathbb S}^{N-1}}\bigg(\int_{B_r}\frac{|\widetilde
    \varphi_{\theta}(x)|^2}{|x|^2}\,dx\bigg)\,d\sigma(\theta)\\
\notag&=\bigg(\frac{N-2}2\bigg)^{\!\!2}
  \int_{{\mathbb S}^{N-1}}\bigg(
  \int_0^{r}\frac{s^{N-1}}{s^2}\tilde v^2(s,\theta)\,ds\bigg)\,d\sigma(\theta)
=\bigg(\frac{N-2}2\bigg)^{\!\!2}\int_{\Omega_r}\frac{v^2(x)}{|x|^2}\,dx,
\end{align}
where $\omega_{N-1}$ denotes the volume of the unit sphere ${\mathbb
  S}^{N-1}$, i.e. $\omega_{N-1}=\int_{{\mathbb
    S}^{N-1}}d\sigma(\theta)$.  On the other hand, from the definition
of $\mu_1(V,s)$, see (\ref{eq:firsteig_r}), it follows that, for every
$s\in(0,r)$,
\begin{equation}\label{eq:3ang}
\int_{C_s}
\left(|\nabla_{{\mathbb S}^{N-1}}v(s,\theta)|^2-V(\theta)
    v^2(s,\theta)\right)\,d\sigma(\theta)
  \geq \mu_1(s,V)\int_{C_s}v^2(s,\theta)\,d\sigma(\theta).
\end{equation}
From (\ref{eq:27}), (\ref{eq:5ang}), and (\ref{eq:3ang}), we deduce that
\begin{align*}
    \int_{\Omega_r} \bigg( |\nabla
    v(y)|^2-\frac{V(\frac{y}{|y|})}{|y|^2} v^2(y)\bigg)\,dy+
    \frac{N-2}{2r}\int_{S_r}v^2(y)\,d\sigma
\geq \int_{\Omega_r}
 \left[\bigg(\frac{N-2}2\bigg)^{\!\!2}\!\!+\mu_1(|x|,V)\right]
  \frac{v^2(x)}{|x|^2}\,dx
\end{align*}
for all $v\in C^\infty_{\rm c}(\tildO\cap \overline{B}_r)$,
which, by density, yields the stated inequality for all
$H^1(\Omega_r)$-functions vanishing on $\Gamma_r$.
\end{pf}

\begin{Corollary}\label{c:hardyboundary}
  Let $N\geq 3$ and $V\in L^{\infty}(\SN)$ such that $\Lambda(V)<1$,
  where $\Lambda(V)$ is defined in (\ref{eq:LambdaV}). Then, there
  exist $R_0\in (0,\widetilde R)$ and $C_{N,V}>0$ such that, for every
  $r\in (0,R_0)$ and $v\in H^1(\Omega_r)$ such that $v=0$ on
  $\Gamma_r$, the following inequalities hold
  \begin{equation}\label{eq:29}
    \int_{\Omega_r} \!\!\bigg( |\nabla
    v(y)|^2-\frac{V(\frac{y}{|y|})}{|y|^2} v^2(y)\bigg)dy+
    \frac{N-2}{2r}\!\!\int_{S_r}\!\!v^2(y)d\sigma
    \geq \frac12\bigg(\!\bigg(\frac{N-2}{2}\bigg)^{\!\!2}+\mu_1(V)\!\bigg)\!\!
    \int_{\Omega_r}\!\! \frac{v^2(y)}{|y|^2}\,dy,
  \end{equation}
  \begin{equation}\label{eq:30}
    \int_{\Omega_r} \!\!\bigg( |\nabla
    v(y)|^2-\frac{V(\frac{y}{|y|})}{|y|^2} v^2(y)\bigg)dy+
    \frac{1+\Lambda(V)}{2}\frac{N-2}{2r}\int_{S_r}\!\!v^2(y)d\sigma
    \geq  \frac{1-\Lambda(V)}{2}\!\!
    \int_{\Omega_r}\!\! |\nabla
    v(y)|^2\,dy,
  \end{equation}
and
  \begin{gather}\label{eq:31}
    \int_{\Omega_r} \bigg( |\nabla
    v(y)|^2-\frac{V(\frac{y}{|y|})}{|y|^2} v^2(y)\bigg)\,dy+
    \frac{N-2}{2r}
\frac{\Lambda(V)+3}4
\int_{S_r}v^2(y)\,d\sigma\\
    \notag\geq
C_{N,V}
\bigg(    \int_{\Omega_r}\!\! \bigg(|\nabla
    v(y)|^2+ \frac{v^2(y)}{|y|^2}\bigg)\,dy+
\|v\|_{L^{2^*}(\Omega_r)}^2\bigg).
  \end{gather}
\end{Corollary}
\begin{pf}
  Inequality (\ref{eq:29}) follows from Lemmas \ref{l:con_eig} and
  \ref{l:hardyboundary}.  To prove (\ref{eq:30}) we observe that if
  $R_0$ is sufficiently small, then, by (\ref{eq:28}) and assumption
  (\ref{eq:V_assumption}), $\Lambda(V,r)<\frac{\Lambda(V)+1}2$ for all
  $r\in (0,R_0)$.  Consequently for all $v\in C^\infty_{\rm
    c}(\tildO\cap \overline{B}_r)$ with $r\in (0,R_0)$, from
  (\ref{eq:Lambda_r}) and (\ref{eq:5ang}), it follows
\begin{align*}
  \int_{\Omega_r}&\frac{V(\frac{y}{|y|})}{|y|^2} v^2(y)\,dy=
  \int_0^r s^{N-3}\bigg(\int_{C_s}V(\theta)v^2(s,\theta)\,d\sigma(\theta)\bigg)ds\\
  &\leq \int_0^r s^{N-3}\Lambda(V,s)\bigg( \int_{C_s}
  \bigg[\big|\nabla_{\mathbb S^{N-1}}v(s,\theta)\big|^2+
  \Big(\frac{N-2}{2}\Big)^2|v(s,\theta)|^2\bigg]\,d\sigma(\theta)\bigg)ds\\
  &\leq \frac{\Lambda(V)+1}2\int_0^r s^{N-3}\bigg( \int_{C_s}
  \bigg[\big|\nabla_{\mathbb S^{N-1}}v(s,\theta)\big|^2+
  \Big(\frac{N-2}{2}\Big)^2|v(s,\theta)|^2\bigg]\,d\sigma(\theta)\bigg)ds\\
  &\leq \frac{\Lambda(V)+1}2 \int_0^r \bigg( s^{N-3}\bigg( \int_{C_s}
  \big|\nabla_{\mathbb S^{N-1}}v(s,\theta)\big|^2d\sigma(\theta)\bigg)
  +s^{N-1}\bigg(\int_{C_s} |\partial_s\tilde v(s,\theta)|^2
  d\sigma(\theta)\bigg)\bigg)ds\\
  &\qquad+\frac{\Lambda(V)+1}2
  \frac{N-2}{2r}r^{N-1}\int_{C_r} v^2(r,\theta)\,d\sigma(\theta)\\
  &=\frac{\Lambda(V)+1}2\bigg(\int_{\Omega_r}|\nabla v(y)|^2\,dy
  +\frac{N-2}{2r}\int_{S_r} v^2(y)\,dy\bigg)
\end{align*}
which yields (\ref{eq:30}) by density.
From summation of (\ref{eq:29}) and (\ref{eq:30}) and
 Sobolev embeddings, it follows that, for every $r\in (0,R_0)$ and
  $v\in H^1(\Omega_r)$ such that $v=0$ on $\Gamma_r$,
  \begin{multline}\label{eq:33}
    \int_{\Omega_r} \bigg( |\nabla
    v(y)|^2-\frac{V(\frac{y}{|y|})}{|y|^2} v^2(y)\bigg)\,dy+
    \frac{N-2}{2r}\frac{3+\Lambda(V)}{4}\int_{S_r}v^2(y)\,d\sigma\\
    \geq \frac{\widetilde S_N}{4}\min\bigg\{
\bigg(\frac{N-2}{2}\bigg)^{\!\!2}+\mu_1(V),1-\Lambda(V)\bigg\}
   \bigg(\int_{\Omega_r} |v(y)|^{2^*}\,dy\bigg)^{\!\!2/2^*},
  \end{multline}
where  $\widetilde S_N>0$ is the best constant of the Sobolev embedding
$H^1(B_1)\subset L^{2^*}(B_1)$. By summing up (\ref{eq:29}),
  (\ref{eq:30}), (\ref{eq:33}), we conclude that
(\ref{eq:31}) holds with
$$
C_{N,V}=\frac{\min\big\{1,{\widetilde S_N}/{2}\big\} \min\big\{
      \big(\frac{N-2}{2}\big)^{2}\!+\mu_1(V),1-\Lambda(V)\big\}
    }6.
$$
\end{pf}

Repeating the same arguments carried out in this section for the
family of domains $\Omega_r$, we can prove analogous estimates
on the domains  $\Omega_r\cup\mathcal C$.
\begin{Corollary}\label{c:suppC}
  Under the same assumptions as in Corollary \ref{c:hardyboundary},
  there exists $\widetilde R_0$ such that for every $r\in
  (0,\widetilde R_0)$ and $v\in H^1(\Omega_r\cup \mathcal C)$ such
  that $v=0$ on $\partial(\Omega_r\cup \mathcal C)\cap B_r$, the
  following inequality holds
  \begin{gather*}
    \int_{\Omega_r\cup \mathcal C} \bigg( |\nabla
    v(y)|^2-\frac{V(\frac{y}{|y|})}{|y|^2} v^2(y)\bigg)\,dy+
    \frac{N-2}{2r}
\frac{\Lambda(V)+3}4
\int_{(\tildO\cup\mathcal C)\cap\partial B_r}v^2(y)\,d\sigma\\
    \notag\geq C_{N,V}
\bigg(    \int_{\Omega_r\cup \mathcal C}\!\! \bigg(|\nabla
    v(y)|^2+ \frac{v^2(y)}{|y|^2}\bigg)\,dy+
\|v\|_{L^{2^*}(\Omega_r\cup \mathcal C)}^2\bigg).
  \end{gather*}
\end{Corollary}

\section{A Brezis-Kato type estimate in dimension $N\geq
  3$}\label{sec:brezis-kato-type}

Throughout this section, we assume $\tildA,\tildb,\widetilde h$ as in
(\ref{eq:tildA}), (\ref{eq:tildb}), (\ref{eq:tildh}) with
$A,\bi,\Psi,h,V$ as in assumptions
(\ref{eq:matrix1}--\ref{eq:LambdaV}), (\ref{eq:def_psi}), and let
$\tildO$ as in (\ref{eq:tildeomega}) with $\Omega$ satisfying
(\ref{eq:omega}) and (\ref{eq:phi1}--\ref{eq:63}). We also assume
that $W\in L^1_{\rm loc}(\tildO)$ satisfies  the form-bounded condition
$$
\sup_{v\in
  H^1(\tildO)\setminus\{0\}}\frac{\int_{\tildO}|W(y)|v^2(y)\,dy}
{\|v\|^2_{H^1(\tildO)}}<+\infty,
$$
see \cite{Mazja}. The above condition in particular implies that for every
$v\in H^1(\tildO)$, $Wv\in H^{-1}(\tildO)$.
 Let $w\in
H^1(\tildO)\setminus\{0\}$ be a weak solution to
\begin{equation}\label{eq:tildequation_W}
\begin{cases}
  -\dive(\tildA(y)\nabla w(y))+\tildb(y)\cdot\nabla
  w(y)-
\dfrac{V\big(\frac{y}{|y|}\big)}{|y|^2}w(y)
=
  \widetilde h(y)w(y)+W(y)w(y),&\text{in }\tildO,\\[5pt]
  w=0,&\hskip-1.2cm\text{on }\partial\tildO\cap B_{\widetilde R}.
\end{cases}
\end{equation}

\begin{Proposition} \label{SMETS}
Let $N\ge 3$ and let $w$ be a weak solution of (\ref{eq:tildequation_W}).
If $W_+\in
L^{N/2}(\tildO)$, letting
$$
q_{\rm lim}:=\left\{
\begin{array}{ll}
\frac{2^*}{2} \min\left\{\frac{8}{\Lambda(V)+1}-2,2^*\right\}, &
\qquad \text{if }\quad \Lambda(V)>0, \\
\frac{(2^*)^2}2, & \qquad \text{if } \quad \Lambda(V)=0 ,
\end{array}
\right.
$$
then for every $1\leq q<q_{\rm lim}$ there exists $r_q>0$ depending
on $q,N,\tildA, \tildb, V,\widetilde h$ such that
$w\in L^q(\Omega_{r_q})$ with $\Omega_{r_q}$ as in (\ref{eq:omerar}).
\end{Proposition}

\begin{pf}
For any $2<\tau<\frac{2}{2^*} q_{\rm lim}$ define
$C(\tau):=\frac{4}{\tau+2}$ and let $\ell_\tau>0$ be
large enough so that
\begin{equation} \label{elleq}
\bigg(\int\limits_{\tildO\cap\{W_+(y) \geq\ell_\tau\}}\!\!\!\!\!\!
W_+^{\frac{N}{2}}\!(y) \, dy
\bigg)^{\!\!\frac{2}{N}}<\frac{S_N(2C(\tau)-\Lambda(V)-1)}4
\end{equation}
where
$$
S_N=\inf_{\phi\in \mathcal D^{1,2}(\R^N)\setminus \{0\}}
\frac{\int_{\R^N}|\nabla \phi(y)|^2 dy}{\left(\int_{\R^N}
    |\phi(y)|^{2^*}dy\right)^{2/2^*}} .
$$
For any $\phi\in H^1_0(\tildO)$, by H\"older and Sobolev inequalities and
\eqref{elleq}, we have
\begin{align} \label{estW}
\int_{\tildO}
W(y)|\phi(y)|^2\,dy &
\leq \ell_\tau \int_{\tildO} |\phi(y)|^2 \, dy +
\bigg(\int\limits_{\tildO\cap\{W_+(y) \geq\ell_\tau\}}\!\!\!\!\!\!
W_+^{\frac{N}{2}}(y) \, dy \bigg)^{\!\!\frac{2}{N}} \left(\int_{\tildO}
|\phi(y)|^{2^*}\, dy\right)^{\!\!\frac{2}{2^*}}
\\
\notag & \leq \ell_\tau \int_{\tildO} |\phi(y)|^2 \,
dy+\frac{2C(\tau)-\Lambda(V)-1}4 \int_{\tildO} |\nabla \phi(y)|^2 \, dy .
\end{align}
Let $r\in(0,\widetilde R)$ small to chosen later and
$\eta\in C^\infty_c (B_r)$ be such that $\eta\equiv 1$ in
$B_{r/2}$. Let us define $v(y):=\eta(y)w(y)\in H^1_0(\Omega_r)$. Then $v$
is a $H^1(\Omega_r)$-weak solution of the equation
\begin{equation} \label{eq:g} -\dive(\tildA(y)\nabla v(y))+\tildb
  (y)\cdot \nabla v(y)-\frac{V(\frac{y}{|y|})}{|y|^2} v(y)=\widetilde
  h(y)v(y)+W(y)v(y)+g(y), \quad \text{in } \Omega_r,
\end{equation}
where $g(y)=-\dive(\tildA(y)\nabla \eta(y))w(y)-2\tildA(y)\nabla
w(y)\cdot \nabla\eta(y)+(\tildb(y)\cdot \nabla \eta(y))w(y)\in
L^2(\Omega_r)$.  For any $n\in\N$, $n\geq 1$, let us define the
function $v^n:=\min\{|v|,n\}$. Testing \eqref{eq:g} with
$(v^n)^{\tau-2} v\in H^1_0(\Omega_r)$ we obtain
\begin{multline} \label{vn}
 \int_{\Omega_r}  (v^n(y))^{\tau-2} \tildA(y)\nabla v(y)\cdot \nabla v(y) \, dy \\
 +(\tau-2) \!\!\!  \int_{\Omega_r} \!\!\!(v^n(y))^{\tau-3}|v(y)|
 \chi_{\{|v(y)|<n\}}(y) \tildA(y)\nabla v(y)\cdot \!\nabla v(y) \, dy
 -\!\!\!\int_{\Omega_r} \!\!\! \frac{V(\frac{y}{|y|})}{|y|^2}
 (v^n(y))^{\tau-2} v^2(y)\, dy
 \\
 \ =-\int_{\Omega_r} (\tildb(y)\cdot \nabla v(y)) (v^n(y))^{\tau-2}
 v(y)\, dy
 + \int_{\Omega_r} \widetilde h(y)(v^n(y))^{\tau-2} v^2(y)\, dy \\
+\int_{\Omega_r} W(y)(v^n(y))^{\tau-2} v^2(y)\, dy + \int_{\Omega_r}
g(y)(v^n(y))^{\tau-2} v(y)\, dy .
\end{multline}
We observe that the following identities hold true
\begin{equation} \label{id}
\begin{cases}
  \tildA \nabla((v^n)^{\frac{\tau}{2}-1}v)\!\cdot\!
  \nabla((v^n)^{\frac{\tau}{2}-1}v)= (v^n)^{\tau-2} \tildA \nabla v
  \!\cdot\! \nabla v+\tfrac{(\tau-2)(\tau+2)}4 (v^n)^{\tau-2}
  \alchi_{\{|v|<n\}} \tildA \nabla v \!\cdot\! \nabla v,
  \\[5pt]
  |\nabla((v^n)^{\frac{\tau}{2}-1}v)|^2= (v^n)^{\tau-2} |\nabla
  v|^2+\tfrac{(\tau-2)(\tau+2)}4 (v^n)^{\tau-2}
  \alchi_{\{|v|<n\}} |\nabla v|^2 .
\end{cases}
\end{equation}
By \eqref{eq:3}, \eqref{id}, H\"older and Hardy inequalities, we have
\begin{align}\label{estb}
   \bigg|\int_{\Omega_r} (\tildb(y) \cdot \nabla v(y))&
    (v^n(y))^{\tau-2} v(y) \, dy\bigg| \leq C_{\tildb}
  \int_{\Omega_r} |\nabla v(y)| \frac{(v^n(y))^{\tau-2}
    |v(y)|}{|y|^{1-\delta}} \, dy
  \\
  & \notag \leq C_{\tildb} r^{\delta} \left(\int_{\Omega_r}
    (v^n(y))^{\tau-2} |\nabla v(y)|^2 dy\right)^{\!\!1/2}
  \left(\int_{\Omega_r} \frac{((v^n(y))^{\frac{\tau}2 -1}
      v(y))^2}{|y|^2} \, dy\right)^{\!\!1/2} \\ & \notag \leq
  \frac{2C_{\tildb}}{N-2} \, r^{\delta} \int_{\Omega_r} |\nabla
  ((v^n(y))^{\frac{\tau}2 -1} v(y))|^2 dy
\end{align}
for some positive constant $C_{\tildb}$ depending only on $\tildb$.

Then by \eqref{vn}, \eqref{eq:30} applied to the function
$(v^n)^{\frac{\tau}2 -1} v$, \eqref{estW} applied to
$\phi=(v^n)^{\frac{\tau}2 -1} v$, \eqref{estb}, \eqref{eq:2},
\eqref{eq:4} and classical Hardy inequality, we obtain
\begin{gather} \label{vn2} C(\tau) (1-K_{\tildA}\, r^{\delta})
  \int_{\Omega_r}
  |\nabla((v^n(y))^{\frac{\tau}{2}-1}v(y))|^2 \, dy \\
  \notag \leq \int_{\Omega_r} \frac{V(\frac{y}{|y|})}{|y|^2}
  ((v^n(y))^{\frac{\tau}{2}-1}v(y))^2 \, dy -\int_{\Omega_r}
  (\tildb(y)\cdot \nabla v(y)) (v^n(y))^{\tau-2} v(y) \, dy +
  \int_{\Omega_r} \widetilde h(y)((v^n(y))^{\frac{\tau}{2}-1}v(y))^2\,
  dy  \\
  \notag \quad +\int_{\Omega_r}
  W(y)((v^n(y))^{\frac{\tau}{2}-1}v(y))^2\,
  dy + \int_{B_r} g(y)(v^n(y))^{\tau-2} v(y)\, dy \\
  \notag \leq \bigg[\frac{\Lambda(V)+1}{2}+\frac{2C_{\tildb}}{N-2} \,
  r^{\delta}+C_{\widetilde h}
  \bigg(\frac{2}{N-2}\bigg)^{\!\!2}r^\delta+\frac{2C(\tau)-\Lambda(V)-1}4\bigg]
  \int_{\Omega_r} |\nabla((v^n(y))^{\frac{\tau}{2}-1}v(y))|^2 \, dy
  \\
  \notag \quad +\ell_\tau \int_{\Omega_r} (v^n(y))^{\tau-2} (v(y))^2
  \,dy +\int_{\Omega_r} |g(y)|(v^n(y))^{\tau-2} |v(y)|\, dy \ ,
 \end{gather}
 for some positive constants $K_{\tildA}$ depending on $\tildA$ and
$C_{\widetilde h}$ depending on $\widetilde h$.

Arguing as in \cite[Proposition 2.3]{FFT3}, we can easily estimate
\begin{multline}\label{est-g}
\int_{\Omega_r}  |g(y)|(v^n(y))^{\tau-2} |v(y)|\, dy\\
\leq\frac{1}{\tau} \|g\|^\tau_{L^2(\Omega_r)} +\frac{\tau-1}\tau
\bigg(\frac{\omega_{N-1}}N\bigg)^{\!\!\frac{\tau}{2(\tau-1)}-\frac{2}{2^*}}
r^{\frac{N\tau}{2(\tau-1)}-N+2}  S_N^{-1} \int_{\Omega_r}
|\nabla((v^n(y))^{\frac{\tau}2-1} v(y))|^{2} \, dy.
\end{multline}
Inserting \eqref{est-g} into \eqref{vn2} and using Sobolev embedding,
we obtain
\begin{align} \label{lim-q} S_N \bigg[& \frac{2C(\tau)-\Lambda(V)-1}4
  -\left(C(\tau)K_{\tildA}+\frac{2C_{\tildb}}{N-2}+C_h
    \Big(\frac{2}{N-2}\Big)^2\right)r^\delta \\
  & \notag -\frac{\tau-1}\tau
  \bigg(\frac{\omega_{N-1}}N\bigg)^{\!\!\frac{\tau}{2(\tau-1)}-\frac{2}{2^*}}
  r^{\frac{N\tau}{2(\tau-1)}-N+2} S_N^{-1} \bigg] \left(
    \int_{\Omega_r}
    |(v^n(y))^{\frac{\tau}2-1} v(y)|^{2^*} \, dy\right)^{\!\!2/2^*} \\
  \notag & \leq \frac{1}{\tau} \|g\|^\tau_{L^2(\Omega_r)} +\ell_\tau
  \int_{\Omega_r} (v^n(y))^{\tau-2} (v(y))^2 \,dy .
\end{align}
 Since $\tau<\frac{2}{2^*}
q_{\rm lim}$ then $2C(\tau)-\Lambda(V)-1$ is positive and
$\frac{N\tau}{2(\tau-1)}-N+2$ is also positive. Hence we may fix $r$
small enough in such a way that the left hand side of
\eqref{lim-q} becomes positive. Since $v\in L^\tau(B_{r})$, letting
$n\to +\infty$, the right hand side of \eqref{lim-q} remains
bounded and hence by Fatou Lemma we infer that $v\in
L^{\frac{2^*}2 \tau}(B_{r})$. Since $\eta\equiv 1$ in $B_{r/2}$ we
may conclude that $w\in L^{\frac{2^*}2 \tau}(B_{r/2})$. This
completes the proof of the lemma.
\end{pf}

\section{The monotonicity formula} \label{sec:monotonicity-formula}
Let $\tildA,\tildb,\tilde f,\widetilde h$ be as in
(\ref{eq:tildA}--\ref{eq:tildh}) with $A,\bi,\Psi,f,h,V$ as in
assumptions (\ref{eq:matrix1}--\ref{eq:F_assumption2}),
(\ref{eq:def_psi}). Let $\tildO$ be as in (\ref{eq:tildeomega}) with
$\Omega$ satisfying (\ref{eq:omega}) and (\ref{eq:phi1}--\ref{eq:63}).
Let $w\in H^1(\tildO)\setminus\{0\}$ be a non-trivial weak solution
to~(\ref{eq:tildequation}).

For every $r\in(0,\tilde R)$ let us define
\begin{align}
&\label{eq:D(r)}D(r)=\frac1{r^{N-2}}\int_{\Omega_r} \bigg(
\tildA\nabla w\cdot \nabla w+\tildb\cdot\nabla
  w\,w-
\dfrac{V\big(\frac{y}{|y|}\big)}{|y|^2}|w|^2-
  \widetilde h w^2-\tilde f(y,w)w\bigg)dy,\\
&\label{eq:H(r)}H(r)=\frac1{r^{N-1}}\int_{S_r}\mu(y)w^2(y)\,d\sigma(y),
\end{align}
where
\begin{align}\label{eq:mu}
\mu(y)=\mu(y',y_N)=|y|^{-2}\tildA(y)y\cdot y.
\end{align}

\begin{Lemma}\label{l:mu}
  Let $N\ge 2$ and let $\mu$ as in (\ref{eq:mu}) with $\tildA$ as in
  (\ref{eq:tildA}). Then
\begin{align}
  \label{eq:11}&\mu(y)=\frac{1}{1+2C_0(\delta+1)y_N
|y|^{-1+\delta}}+O(|y|)
=1+O(|y|^\delta)
\quad\text{as }|y|\to0,\\[5pt]
\label{eq:12}  &\nabla\mu(y)=O(|y|^{-1+\delta})\quad\text{as }|y|\to0.
\end{align}
\end{Lemma}
\begin{pf}
Estimate (\ref{eq:11}) follows from Lemma \ref{l:Ay} and direct calculations.
Differentiating (\ref{eq:mu}) we obtain
$$
\nabla \mu(y)=-2|y|^{-4} (\tildA(y)y\cdot
y)y+|y|^{-2}(d\tildA(y)y)y+2|y|^{-2}\tildA(y) y.
$$
From (\ref{eq:11}) and (\ref{eq:2}) we then deduce
\begin{align*}
\nabla\mu(y)&=
-\frac{2|y|^{-2}y}{1+2C_0(\delta+1)y_N
|y|^{-1+\delta}}+O(1)+|y|^{-2}O(|y|^{1+\delta})+2|y|^{-2}\big(y+O(|y|^{1+\delta})\big)\\
&=\frac{4C_0(\delta+1)y_N
|y|^{-3+\delta}y}{1+2C_0(\delta+1)y_N
|y|^{-1+\delta}}+O(|y|^{-1+\delta})=O(|y|^{-1+\delta})
\end{align*}
as $|y|\to 0$.
\end{pf}

\begin{Lemma}\label{l:5.2}
  Let $N\ge 2$ and let $\tildA, V$ be as in (\ref{eq:V_assumption}),
  (\ref{eq:tildA}) with $A$ as in (\ref{eq:matrix1}-\ref{eq:matrix2}).
  Define the function
$$
\boldbeta(y):=\frac{\tildA(y)y}{\mu(y)} .
$$
Then we have
\begin{align}
 \label{eq:17} &\boldbeta(y)=y+O(|y|^{1+\delta})=O(|y|) \text{ as }|y|\to 0,\\
  \label{eq:18}&\mathop{\rm Jac}\boldbeta(y)=\tildA(y)+O(|y|^{\delta})
  ={\rm Id}_N +O(|y|^{\delta})\text{ as }|y|\to 0,\\
 \label{eq:19} &\dive\boldbeta(y)=N+O(|y|^{\delta}) \text{ as }|y|\to 0,\\
 \label{eq:20}&\boldbeta(y)\cdot \nabla_{\SN}V(y/|y|)=y\cdot
 \nabla_{\SN}V(y/|y|)+O(|y|^{1+\delta}) =O(|y|^{1+\delta}) \text{ as
 }|y|\to 0,
\end{align}
\end{Lemma}

\begin{pf}
It follows from the definitions of $\boldbeta$ and $\mu$.
\end{pf}

\noindent From  (\ref{eq:2}--\ref{eq:4}) we we derive the following lemma.

\begin{Lemma} \label{l:hprime} Let $N\ge 2$. Then $H\in W^{1,1}_{\rm
    loc}(0,\widetilde R)$ and
\begin{equation}\label{eq:H'}
  H'(r)=\frac{2}{r^{N-1}} \int_{S_r}
\mu(y)w(y)\frac{\partial w}{\partial \nu}(y)d\sigma(y)
+H(r)O(r^{-1+\delta})\quad \text{as }r\to0^+
\end{equation}
in a distributional sense and for a.e. $r\in (0,\widetilde R)$,
where $\nu=\nu(y)$ is the unit outer normal vector to $S_r$, i.e
\begin{equation}\label{eq:nu}
\nu(y)=\frac{y}{|y|}.
\end{equation}
\end{Lemma}
 \begin{pf}
We notice that, for all $r\in(0,\widetilde R)$,
$$
H(r)=\int_{C_r} \mu(r\theta)|w(r\theta)|^2d\sigma
$$
where $C_r$ is defined in (\ref{eq:Cr}).  For every $\phi\in
C^{\infty}_{\rm c}(0,{\widetilde R})$
\begin{align*}
  \int_0^{\widetilde R} H(t)\phi'(t)\,dt&= \int_0^{\widetilde R}
  \bigg(\int_{C_t} \mu(t\theta)|w(t\theta)|^2d\sigma\bigg)\phi'(t)\,dt
=  \int_{\Omega_{\widetilde R}}\frac{\mu(y)w^2(y)\nu(y)}{|y|^{N-1}}
  \cdot\nabla \tilde\phi(y)\,dy\\
  &=-\int_{\Omega_{\widetilde R}}\dive\Big(\mu(y)w^2(y)
  \frac{y}{|y|^{N}}\Big)\tilde\phi(y)\,dy=\\
  &=-\int_{\Omega_{\widetilde R}}\frac{w^2(y)\nabla\mu(y)
    +2w(y)\mu(y)\nabla w(y)}{|y|^{N-1}}\cdot \nu(y)\tilde\phi(y)\,dy\\
  &=-\int_0^{\widetilde R}\bigg(\int_{C_t}
  \big(2\mu(t\theta)w(t\theta)\nabla
  w(t\theta)\cdot\theta+w^2(t\theta)\nabla\mu(t\theta)\cdot\theta\big)
d\sigma\bigg)\phi(t)\,dt,
\end{align*}
where $\tilde\phi(y):=\phi(|y|)$. Hence
\begin{equation}\label{eq:10}
  H'(t)=\int_{C_t} \big(2\mu(t\theta)w(t\theta)\nabla
  w(t\theta)\cdot\theta\big)\,d\sigma+\int_{C_t}
  \big(w^2(t\theta)\nabla\mu(t\theta)
  \cdot\theta\big)d\sigma
\end{equation}
in a distributional sense in $(0,\widetilde R)$. From
$w,\frac{\partial w}{\partial\nu}\in L^2(\Omega_{\widetilde R})$ we
deduce that $H\in W^{1,1}_{\rm loc}(0,\widetilde R)$. Furthermore
(\ref{eq:10}) holds a.e. and can be rewritten as
\begin{align*}
H'(t)&=\frac{2}{t^{N-1}} \int_{S_t}
\mu(y)w(y)\frac{\partial w}{\partial \nu}(y)d\sigma(y)
+\frac{1}{t^{N-1}} \int_{S_t}w^2(y)\nabla\mu(y)\cdot\nu(y)d\sigma(y)\\
&=\frac{2}{t^{N-1}} \int_{S_t}
\mu(y)w(y)\frac{\partial w}{\partial \nu}(y)d\sigma(y)
+H(t)O(t^{-1+\delta})
\end{align*}
as $t\to0^+$ thus proving (\ref{eq:H'}).
\end{pf}

\begin{Lemma}\label{l:HprimeD} Let $N\ge 2$.
Let $D$ and $H$ the functions defined in (\ref{eq:D(r)}--\ref{eq:H(r)}). Then
\begin{align}\label{eq:H'2}
&  H'(r)=\frac{2}{r^{N-1}} \int_{S_r}
(\tildA\nabla w\cdot\nu)w\,d\sigma(y)
+H(r)O(r^{-1+\delta}),\\
\label{eq:35}&H'(r)=\frac2rD(r)+H(r)O(r^{-1+\delta})
\end{align}
as $r\to 0^+$.
\end{Lemma}
\begin{pf}
We have that
\begin{align*}
\int_{S_r}(\tildA\nabla w\cdot\nu)w\,d\sigma
=\int_{S_r}\mu w\frac{\partial w}{\partial\nu}\,d\sigma
+\frac12\int_{S_r}\boldalfa\cdot\nabla(w^2)\,d\sigma
\end{align*}
where
$$
\boldalfa(y)=\frac{\mu(y)(\boldbeta(y)-y)}{|y|}.
$$
Since $\boldalfa(y)\cdot y=0$ and, in view of (\ref{eq:12}), (\ref{eq:17}),
 and (\ref{eq:19}),
\begin{align*}
  \dive\boldalfa=\bigg(\frac{\nabla\mu}{|y|}-\mu\frac{y}{|y|^3}\bigg)
  (\boldbeta-y)
  +\frac{\mu}{|y|}(\dive\boldbeta-N)=O(|y|^{-1+\delta})\quad\text{ as
  }|y|\to 0,
\end{align*}
we deduce that
 \begin{align*}
\int_{S_r}(\tildA\nabla w\cdot\nu)w\,d\sigma=
\int_{S_r}\mu w\frac{\partial w}{\partial\nu}\,d\sigma
  -\frac12\int_{S_r}\dive(\boldalfa)w^2\,d\sigma=\int_{S_r}\mu
  w\frac{\partial w}{\partial\nu}\,d\sigma
  +O(r^{-2+N+\delta})H(r)
\end{align*}
and hence  (\ref{eq:H'2}) follows from
 (\ref{eq:H'}).
  Multiplying equation (\ref{eq:tildequation}) by $w$ and
  integrating on $\Omega_r$, from  (\ref{eq:H'2})  we obtain that
\begin{align*}
  r^{N-2}D(r)=
\int_{S_r}(\tildA\nabla w\cdot\nu)w\,d\sigma
 =\frac{r^{N-1}}{2}H'(r) +O(r^{-2+N+\delta})H(r)
\end{align*}
as $r\to0^+$, thus proving (\ref{eq:35}).
\end{pf}

\noindent We proceed by distinguishing the cases $N\ge 3$ and $N=2$.

\subsection{The case $N\ge 3$}

By (\ref{eq:5_2}) and Sobolev embedding, we infer that the function
$$
W(y):=\left\{
\begin{array}{ll}
\frac{\tilde f(y,w(y))}{w(y)}, & \qquad \text{if } w(y)\neq 0, \\
0, & \qquad \text{if } w(y)=0,
\end{array}
\right.
$$
belongs to $L^{N/2}(\tildO)$ and hence we may apply Proposition
\ref{SMETS} to the function $w$. Therefore, throughout this
section, we may fix
\begin{equation}\label{eq:50}
2^*<q<q_{\rm lim}
\end{equation}
 and $r_q$ as in
Proposition \ref{SMETS} in such a way that $w\in L^q(\Omega_{r_q})$.

\begin{Lemma} \label{welld} There exist $r_0\in(0,\min\{\widetilde R,
  r_q\})$ and a constant
  $\overline{C}=\overline{C}(N,V,\tildA,\tildb,\tilde f,\widetilde
  h,w)>0$~depending on $N$, $V$, $\tildA$, $\tildb$, $\tilde f$,
  $\widetilde h$, $w$ but independent of $r$ such that such that, for
  all $r\in (0,r_0)$,
\begin{align*}
(i)&\quad \mu(y)>1/2 \quad \text{for all }y\in B_{r},\\
(ii)&\quad r^{N-2}\bigg(D(r)+\frac{N-2}{2}H(r)\bigg)\geq
    \overline{C} \bigg(\int_{\Omega_r}\!\! \bigg(|\nabla
    w(y)|^2+ \frac{w^2(y)}{|y|^2}\bigg)\,dy+
\|w\|_{L^{2^*}(\Omega_r)}^2\bigg),\\
(iii)&\quad  H(r)>0,
\end{align*}
where $D$ and $H$ are defined in (\ref{eq:D(r)}) and (\ref{eq:H(r)}).
\end{Lemma}
\begin{pf} Estimate $(i)$ near $0$ follows from the definition of
  $\mu$. To prove $(ii)$, we observe that, from (\ref{eq:D(r)}),
  (\ref{eq:H(r)}), (\ref{eq:11}), and (\ref{eq:2}--\ref{eq:4}), it
  follows that
\begin{multline*}
r^{N-2}\bigg(D(r)+\frac{N-2}{2}H(r)\bigg)\geq
    \int_{\Omega_r} \bigg( |\nabla
    w(y)|^2-\frac{V(\frac{y}{|y|})}{|y|^2} w^2(y)\bigg)\,dy+
    \frac{N-2}{2r}(1+O(r^\delta))\int_{S_r}w^2(y)\,d\sigma\\
+  O(r^\delta)
    \int_{\Omega_r} \bigg( |\nabla
    w(y)|^2+\frac{ w^2(y)}{|y|^2}\bigg)\,dy
-C_{\tilde f}\|w\|_{L^{2^*}(\Omega_r)}^{2^*-2}
\|w\|_{L^{2^*}(\Omega_r)}^{2}\\
\geq
\int_{\Omega_r} \bigg( |\nabla
    w(y)|^2-\frac{V(\frac{y}{|y|})}{|y|^2} w^2(y)\bigg)\,dy+
    \frac{N-2}{2r}
\frac{\Lambda(V)+3}4
\int_{S_r}w^2(y)\,d\sigma\\
+  \Big(O(r^\delta)-C_{\tilde f}\|w\|_{L^{2^*}(\Omega_r)}^{2^*-2}\Big)
\bigg(    \int_{\Omega_r} \bigg( |\nabla
    w(y)|^2+\frac{ w^2(y)}{|y|^2}\bigg)\,dy
+\|w\|_{L^{2^*}(\Omega_r)}^{2}\bigg)
\end{multline*}
as $r\to 0^+$, which, together with (\ref{eq:31}), yields $(ii)$
provided $r$ is sufficiently small.

To prove the positivity of $H$ near $0$,
  suppose by contradiction that there exists a sequence $r_n\to 0^+$
  such that $H(r_n)=0$. Since $\mu(y)>0$ if $|y|$ is sufficiently
  small, then  $w= 0$ a.e. on $S_{r_n}$ for $n$ sufficiently
  large and thus $w\in H^1_0(\Omega_{r_n})$.  Multiplying both sides
  of (\ref{eq:tildequation}) by $w$ and using estimate $(ii)$,
 we obtain, for $n$ sufficiently large,
\begin{align*}
0&=\int_{\Omega_{r_n}}
 \bigg(
\tildA\nabla w\cdot\nabla w+\tildb\cdot\nabla
  w\,w-
\dfrac{V\big(\frac{y}{|y|}\big)}{|y|^2}|w|^2-
  \widetilde h|w|^2-\tilde f(y,w(y))w\bigg)dy\\
&\geq
\overline{C} \bigg(\int_{\Omega_{r_n}}\!\! \bigg(|\nabla
    w(y)|^2+ \frac{w^2(y)}{|y|^2}\bigg)\,dy+
\|w\|_{L^{2^*}(\Omega_{r_n})}^2\bigg)
\end{align*}
which implies $w\equiv 0$ in
$\Omega_{r_n}$ for $n$ large.  Applying away from $0$ classical
unique continuation principles for second order elliptic equations
with locally bounded coefficients (see e.g.  \cite{wolff}), we
conclude that $w=0$ a.e. in $\tildO$, a contradiction.
\end{pf}

\begin{remark}\label{r:reg_boundary}
  If $w\in H^1(\tildO)$ is a weak solution to (\ref{eq:tildequation}),
  with $\tildA,\tildb,\tilde f,\widetilde h$ as in
  (\ref{eq:tildA}--\ref{eq:tildh}), $A,\bi,\Psi,f,h,V$ as in
  assumptions (\ref{eq:matrix1}--\ref{eq:F_assumption2}),
  (\ref{eq:def_psi}), and $\tildO$ as in (\ref{eq:tildeomega}) with
  $\Omega$ satisfying (\ref{eq:omega}) and (\ref{eq:phi1}), then by
  classical elliptic regularity theory and a Brezis-Kato type
  iteration \cite{BrezisKato}, we have that $w\in W^{2,p}_{\rm
    loc}(\tildO)$ for all $1\leq p<\infty$. In particular $w\in
  H^{2}_{\rm loc}(\tildO)\cap C^{1,\alpha}_{\rm loc}(\tildO)$ for all
  $\alpha\in(0,1)$.  Using a local $C^2$-parametrization of the
  boundary away from the origin (see assumption (\ref{eq:phi1})) and
  classical regularity results for elliptic equations with homogeneous
  boundary conditions on half-spaces, we can deduce that $w\in
  C^{1,\alpha}_{\rm loc}(\overline{\Omega}_{\widetilde R}\setminus
  \{0\})\cap H^2(\Omega_{\widetilde R}\setminus
  \overline{\Omega}_{r})$ for all $r\in(0,\widetilde R)$.
\end{remark}

\begin{Proposition}\label{p:poho}
  Let $N\geq3$, $\tildA,\tildb,\tilde f,\widetilde h$ as in
  (\ref{eq:tildA}--\ref{eq:tildh}) with $A,\bi,\Psi,f,h,V$ as in
   (\ref{eq:matrix1}--\ref{eq:F_assumption2}),
  (\ref{eq:def_psi}), and let $\tildO$ as in (\ref{eq:tildeomega})
  with $\Omega$ satisfying (\ref{eq:omega}) and
  (\ref{eq:phi1}--\ref{eq:63}).  If $w\in H^1(\tildO)\setminus\{0\}$
  is a weak solution to~(\ref{eq:tildequation}), then for
  a.e. $r\in(0,\tilde R)$
\begin{multline}\label{eq:poho}
  r\int_{S_r} (\tildA\nabla w\cdot\nabla w)\,d\sigma - 2r\int_{S_r}
  \frac{|\tildA \nabla w\cdot\nu|^2}\mu\,d\sigma
  -\int_{\Gamma_r} \frac{|\nabla
    w|^2}{\mu}(\tildA\tilde\nu\cdot\tilde\nu)(\tildA y\cdot\tilde\nu)\,d\sigma\\
  = \int_{\Omega_r}(\dive\boldbeta)\tildA\nabla w\cdot\nabla w\,dy
  -2\int_{\Omega_r}(\mathop{\rm Jac}\boldbeta)(\tildA\nabla w)
  \cdot\nabla w\,dy\\
  +\int_{\Omega_r}(d\tildA \nabla w)\nabla w\cdot\boldbeta\,dy-
  \int_{\Omega_r}2(\boldbeta\cdot\nabla w)(\tildb\cdot\nabla w)\,dy\\
  -\int_{\Omega_r}
  \frac{V(y/|y|)\dive\boldbeta-2V(y/|y|)+|y|^{-1}\boldbeta\cdot
    \nabla_{\SN}V(y/|y|)}{|y|^2}w^2(y)\,dy +r \int_{S_r}
  \frac{V(y/|y|)}{|y|^2} w^2\,d\sigma\\+2
  \int_{\Omega_r}(\boldbeta\cdot\nabla w) \widetilde h w\,dy
  -2\int_{\Omega_r} (\nabla_y\widetilde
  F(y,w)\cdot\boldbeta+\widetilde F(y,w)\dive\boldbeta)\,dy +
 2 r\int_{S_r}\widetilde F(y,w)\,d\sigma,
\end{multline}
where $\boldbeta(y):=\frac{\tildA(y)y}{\mu(y)}$.
\end{Proposition}
\begin{pf}
  By Remark \ref{r:reg_boundary}, $w\in H^2_{\rm
    loc}(\Omega_{\widetilde R})\cap C^1(\overline{\Omega}_{\widetilde
    R}\setminus \{0\})$ and hence for all $r\in (0,\widetilde R)$
  the following Rellich-Ne\u{c}as identity
\begin{multline}\label{eq:rellich_necas}
  \dive\big((\tildA\nabla w\cdot\nabla w)\boldbeta-
  2(\boldbeta\cdot\nabla w)\tildA \nabla w\big)=
  (\dive\boldbeta)\tildA\nabla w\cdot\nabla w-2(\mathop{\rm
    Jac}\boldbeta)(\tildA\nabla w)
  \cdot\nabla w\\
  +(d\tildA \nabla w)\nabla w\cdot\boldbeta-
  2(\boldbeta\cdot\nabla w)(\tildb\cdot\nabla w)\\
  +2\frac{V(y/|y|)}{|y|^2} (\boldbeta\cdot\nabla
  w)w+2(\boldbeta\cdot\nabla w) \widetilde h w +2(\boldbeta\cdot\nabla
  w)\tilde f(y,w)
\end{multline}
is satisfied in a weak sense in $\Omega_{\widetilde
  R}\setminus\overline \Omega_r$.  By (\ref{eq:5_2}) and Hardy
inequality, we have
\begin{multline*}
  \int_0^{\widetilde R} \left[\int_{S_s} \left(|\nabla
      w(y)|^2+\frac{w^2(y)}{|y|^2}+|\widetilde F(y,w(y))|\right)
    \, d\sigma\right] \, ds \\
  =\int_{\Omega_{\widetilde R}} \left(|\nabla
    w(y)|^2+\frac{w^2(y)}{|y|^2}+|\widetilde F(y,w(y))|\right) \,
  dy<+\infty
\end{multline*}
and hence there exists a decreasing sequence $\{\delta_n\}\subset
(0,{\widetilde R})$ such that $\lim_{n\to +\infty} \delta_n=0$ and
\begin{equation} \label{bordo-delta} \delta_n \int_{S_{\delta_n}}
  \left(|\nabla w(y)|^2+\frac{w^2(y)}{|y|^2}+|\widetilde
    F(y,w(y))|\right) \, d\sigma \longrightarrow 0 \text{\quad as }
  n\to +\infty .
\end{equation}
Let $r\in (0,\widetilde R)$.
Integrating (\ref{eq:rellich_necas}) in $\Omega_r\setminus\Omega_{\delta_n}$
and taking into account Remark \ref{r:reg_boundary},
we obtain
\begin{multline}\label{eq:rellich_necas_int}
  \int_{S_r} (\tildA\nabla w\cdot\nabla w)\boldbeta\cdot\nu\,d\sigma -
  2\int_{S_r} (\boldbeta\cdot\nabla w)\tildA \nabla w\cdot\nu\,d\sigma
  -\int_{S_{\delta_n}} (\tildA\nabla w\cdot\nabla w)\boldbeta\cdot\nu\,d\sigma\\
  + 2\int_{S_{\delta_n}} (\boldbeta\cdot\nabla w)\tildA \nabla
  w\cdot\nu\,d\sigma +\int_{\Gamma_r\setminus\Gamma_{\delta_n}}
  (\tildA\nabla w\cdot\nabla w)\boldbeta\cdot\tilde\nu\,d\sigma-
  2\int_{\Gamma_r\setminus\Gamma_{\delta_n}} (\boldbeta\cdot\nabla
  w)\tildA \nabla w\cdot\tilde\nu\,d\sigma\\
  =
  \int_{\Omega_r\setminus\Omega_{\delta_n}}(\dive\boldbeta)\tildA\nabla
  w\cdot\nabla w\,dy
  -2\int_{\Omega_r\setminus\Omega_{\delta_n}}(\mathop{\rm
    Jac}\boldbeta)(\tildA\nabla w)
  \cdot\nabla w\,dy\\
  +\int_{\Omega_r\setminus\Omega_{\delta_n}}(d\tildA \nabla w)\nabla
  w\cdot\boldbeta\,dy-
  \int_{\Omega_r\setminus\Omega_{\delta_n}}2(\boldbeta\cdot\nabla w)
  (\tildb\cdot\nabla w)\,dy\\
  +2\int_{\Omega_r\setminus\Omega_{\delta_n}}\frac{V(y/|y|)}{|y|^2}
  (\boldbeta\cdot\nabla
  w)w\,dy+2\int_{\Omega_r\setminus\Omega_{\delta_n}}(\boldbeta\cdot\nabla
  w) \widetilde h w\,dy
  +2\int_{\Omega_r\setminus\Omega_{\delta_n}}(\boldbeta\cdot\nabla
  w)\tilde f(y,w)\,dy
\end{multline}
with $\nu$  as in (\ref{eq:nu}) and $\tilde\nu$ as in (\ref{eq:normaletilde}).
Since $\boldbeta\cdot y=|y|^2$,
integration by parts yields
\begin{multline}\label{eq:14}
\int_{\Omega_r\setminus\Omega_{\delta_n}}\frac{V(y/|y|)}{|y|^2}
  (\boldbeta\cdot\nabla
  w)w\,dy\\
=-\frac12\int_{\Omega_r\setminus\Omega_{\delta_n}}
\frac{V(y/|y|)\dive\boldbeta-2V(y/|y|)+|y|^{-1}\boldbeta\cdot
\nabla_{\SN}V(y/|y|)}{|y|^2}w^2(y)\,dy\\
+\frac r2  \int_{S_r} \frac{V(y/|y|)}{|y|^2}
  w^2\,d\sigma-\frac{\delta_n}2  \int_{S_{\delta_n}} \frac{V(y/|y|)}{|y|^2}
  w^2\,d\sigma
\end{multline}
and
\begin{multline}\label{eq:13}
\int_{\Omega_r\setminus\Omega_{\delta_n}}
(\boldbeta\cdot\nabla
  w)\tilde f(y,w)\,dy=-\int_{\Omega_r\setminus\Omega_{\delta_n}}
(\nabla_y\widetilde F(y,w)\cdot\boldbeta+\widetilde F(y,w)\dive\boldbeta)\,dy\\
+r\int_{S_r}\widetilde F(y,w)\,d\sigma-
\delta_n\int_{S_{\delta_n}}\widetilde F(y,w)\,d\sigma.
\end{multline}
By definition of $\boldbeta$
\begin{align}\label{eq:15}
\boldbeta\cdot\nabla w=\frac{r}{\mu}\tildA \nabla w\cdot\nu\quad\text{on }S_r.
\end{align}
Since $w=0$ on $\Gamma_r$
\begin{align}\label{eq:16}
\nabla w=\pm|\nabla w|\tilde\nu\quad\text{a.e. on }\Gamma_r
\end{align}
Taking into account (\ref{eq:14}--\ref{eq:16}),
(\ref{eq:rellich_necas_int}) becomes
\begin{multline}\label{eq:rellich_necas_int2}
  r\int_{S_r} (\tildA\nabla w\cdot\nabla w)\,d\sigma - 2r\int_{S_r}
  \frac{|\tildA \nabla w\cdot\nu|^2}\mu\,d\sigma
  -\delta_n\int_{S_{\delta_n}} (\tildA\nabla w\cdot\nabla w)\,d\sigma\\
  + 2\delta_n\int_{S_{\delta_n}} \frac{|\tildA \nabla
    w\cdot\nu|^2}{\mu}\,d\sigma
  -\int_{\Gamma_r\setminus\Gamma_{\delta_n}} \frac{|\nabla
    w|^2}{\mu}(\tildA\tilde\nu\cdot\tilde\nu)(\tildA y\cdot\tilde\nu)\,d\sigma\\
  =
  \int_{\Omega_r\setminus\Omega_{\delta_n}}(\dive\boldbeta)\tildA\nabla
  w\cdot\nabla w\,dy
  -2\int_{\Omega_r\setminus\Omega_{\delta_n}}(\mathop{\rm
    Jac}\boldbeta)(\tildA\nabla w)
  \cdot\nabla w\,dy\\
  +\int_{\Omega_r\setminus\Omega_{\delta_n}}(d\tildA \nabla w)\nabla
  w\cdot\boldbeta\,dy- \int_{\Omega_r\setminus\Omega_{\delta_n}}
  2(\boldbeta\cdot\nabla w)(\tildb\cdot\nabla w)\,dy\\
  -\int_{\Omega_r\setminus\Omega_{\delta_n}}
  \frac{V(y/|y|)\dive\boldbeta-2V(y/|y|)+|y|^{-1}\boldbeta\cdot
    \nabla_{\SN}V(y/|y|)}{|y|^2}w^2(y)\,dy\\
  +r \int_{S_r} \frac{V(y/|y|)}{|y|^2} w^2\,d\sigma-\delta_n
  \int_{S_{\delta_n}} \frac{V(y/|y|)}{|y|^2} w^2\,d\sigma+2
  \int_{\Omega_r\setminus\Omega_{\delta_n}}(\boldbeta\cdot\nabla
  w) \widetilde h w\,dy\\
  -2\int_{\Omega_r\setminus\Omega_{\delta_n}} (\nabla_y\widetilde
  F(y,w)\cdot\boldbeta+\widetilde F(y,w)\dive\boldbeta)\,dy
  +2r\int_{S_r}\widetilde F(y,w)\,d\sigma-
  2\delta_n\int_{S_{\delta_n}}\widetilde F(y,w)\,d\sigma.
\end{multline}
Letting $n\to\infty$ in (\ref{eq:rellich_necas_int2}) and using
(\ref{bordo-delta}), Lemma \ref{l:5.2}, and
(\ref{eq:2}), we obtain that
$$
\int_{\Gamma_r} \frac{|\nabla
  w|^2}{\mu}(\tildA\tilde\nu\cdot\tilde\nu)(\tildA
y\cdot\tilde\nu)\,d\sigma<+\infty
$$
and (\ref{eq:poho}) holds.
\end{pf}

\begin{Lemma}\label{l:boundary_non}
If $q$ is as in (\ref{eq:50}), the function
\begin{equation*}
g(r):=\frac{r^{\frac{q-2^*}{q}}\int_{S_r}|w|^{2^*} \!
  d\sigma}{\Big(\int_{\Omega_r}|w(y)|^{2^*} \! dy\Big)^{\!1-\frac1N}}
\end{equation*}
is well defined and satisfies
\begin{equation*}
g\in L^1(0,r_0)
\quad\text{and}\quad g\geq 0 \text{ a.e. in }(0,r_0).
\end{equation*}
Furthermore
\begin{equation}\label{eq:37}
\int_{S_r}|w|^{2^*} d\sigma\leq
\Big(\frac{\omega_{N-1}}N\Big)^{\!\!
\frac{q-2^*}{qN}}
\|w\|_{L^q(\Omega_{r_0})}^{{2^*}/N}
  \overline{C}^{-1}
g(r)r^{N-2}\bigg(D(r)+\frac{N-2}{2}H(r)\bigg)
\end{equation}
for a.e. $r\in (0,r_0)$ and
\begin{equation}\label{eq:40}
\int_0^r g(s)\,ds\leq N\|w\|_{L^{2^*}(\Omega_{r_0})}^{2^*/N}
 r^{\frac{q-2^*}{q}}
\end{equation}
for all $r\in (0,r_0)$.
\end{Lemma}

\begin{pf}
  From Lemma \ref{welld}, $\int_{\Omega_r} |w(y)|^{2^*} dy>0$ for any
  $r\in(0,r_0)$ and $g$ is well defined in $(0,r_0)$.  Let us denote
  $\beta=\frac{q-2^*}q>0$.  By a direct calculation, we have that
\begin{align}\label{eq:38}
  g(r)&= \frac{r^{\beta}\int_{S_r}|w|^{2^*} \!
    d\sigma}{\Big(\int_{\Omega_r}|w(y)|^{2^*} \! dy\Big)^{\!1-\frac1N}}\\
  &\notag =N\bigg[\frac{d}{dr}\bigg(r^{\beta}
  \bigg(\int_{\Omega_r}|w(y)|^{2^*} \! dy\bigg)^{\!\!1/N}\bigg) -
  \beta\, r^{-1+\beta} \bigg(\int_{\Omega_r}|w(y)|^{2^*} \!
  dy\bigg)^{\!\!1/N}\bigg]
\end{align}
in the distributional sense and for a.e. $r\in (0,r_0)$.
Since
$$
\lim_{r\to0^+}r^{\beta}
    \bigg(\int_{\Omega_r}|w(y)|^{2^*} \! dy\bigg)^{\!\!1/N}=0
\quad\text{and}\quad
\bigg(\int_{\Omega_r}|w(y)|^{2^*} \!
dy\bigg)^{\!\!1/N} =O(1)
$$
as $r\to 0^+$, we have that $g\in L^1(0,r_0)$. Furthermore,
(\ref{eq:40}) follows from
integration of (\ref{eq:38}).

To prove (\ref{eq:37}) we observe that, by H\"older inequality,
Proposition \ref{SMETS}, and Lemma \ref{welld},
\begin{align*}
  \bigg(\int_{\Omega_r}&|w(y)|^{2^*}\,dy\bigg)^{\!\!1-\frac1N} =
\bigg(\int_{\Omega_r}|w(y)|^{2^*}\,dy\bigg)^{\!\!\frac1N}
  \bigg(\int_{\Omega_r}|w(y)|^{2^*}\,dy\bigg)^{\!\!\frac2{2^*}}\\
  &\notag\leq
\Big(\frac{\omega_{N-1}}N\Big)^{\!\!\frac{\beta}{N}}
\|w\|_{L^q(\Omega_{r_0})}^{\frac{2^*}N}
  \overline{C}^{-1}
  r^{\beta+N-2}
  \bigg(D(r)+{{\frac{N-2}{2}}}H(r)\bigg)
\end{align*}
for all $r\in(0,r_0)$, thus implying (\ref{eq:37}).
\end{pf}

\begin{Lemma}\label{l:dprime}
The  function $D$ defined in (\ref{eq:D(r)}) belongs to  $W^{1,1}_{{\rm\
      loc}}(0,r_0)$ and
\begin{align}\label{eq:94}
  D'(r)=B(r) +\frac{1}{r^{N-2}}\!\int_{S_r}\!\! \tildb\cdot\nabla w\,w
  d\sigma +O\Big(r^{-1+\delta}\!+r^{-1+\frac{2(q-2^*)}{q}}\!+g(r)\Big)
  \Big(D(r)+\tfrac{N-2}{2}H(r) \Big)
\end{align}
as $r\to 0^+$, in a distributional sense and for a.e. $r\in (0,r_0)$,
where
\begin{equation}\label{eq:B(r)}
B(r):=\frac2{r^{N-2}} \int_{S_r} \frac{|\tildA \nabla
    w\cdot\nu|^2}\mu\,d\sigma +\frac1{r^{N-1}}\int_{\Gamma_r}
  \frac{|\nabla w|^2}{\mu}(\tildA\tilde\nu\cdot\tilde\nu)(\tildA
  y\cdot\tilde\nu)\,d\sigma
\end{equation}
and $r_0$ is as in Lemma \ref{welld}.
\end{Lemma}
\begin{pf}
By (\ref{eq:poho}) and Lemmas \ref{l:new_coeff} and \ref{l:5.2}, we have that
\begin{align}\label{eq:poho2}
  \int_{S_r} (\tildA\nabla w\cdot&\nabla w)\,d\sigma - \int_{S_r}
  \frac{V(y/|y|)}{|y|^2} w^2\,d\sigma\\
  = &\notag \,2\int_{S_r} \frac{|\tildA \nabla
    w\cdot\nu|^2}\mu\,d\sigma +\frac1r\int_{\Gamma_r} \frac{|\nabla
    w|^2}{\mu}(\tildA\tilde\nu\cdot\tilde\nu)(\tildA y\cdot\tilde\nu)\,d\sigma\\
  &\notag + \frac{N-2}r \bigg(\int_{\Omega_r}\tildA\nabla w\cdot\nabla
  w\,dy-\int_{\Omega_r} \frac{V(y/|y|)}{|y|^2} w^2\,dy\bigg)
  \\
  &\notag -\frac2r\int_{\Omega_r} \Big(\nabla_y\widetilde F(y,w)\cdot
  y+N\widetilde F(y,w)\Big)\,dy +2\int_{S_r}\widetilde F(y,w)\,d\sigma\\
  &\notag +O(r^{-1+\delta})\bigg( \int_{\Omega_r}\tildA\nabla
  w\cdot\nabla w\,dy +\int_{\Omega_r}\frac{w^2(y)}{|y|^2}
  \,dy+\int_{\Omega_r}\Big(|w|^{2}+|w|^{2^*}\Big)\,dy\bigg) .
\end{align}
From (\ref{eq:D(r)}) and  (\ref{eq:poho2})  we obtain
\begin{align}\label{eq:32}
  D'(r)=\,&-\frac{N-2}{r^{N-1}}\int_{\Omega_r} \bigg( \tildA\nabla w\cdot
  \nabla w+\tildb\cdot\nabla w\,w-
  \dfrac{V\big(\frac{y}{|y|}\big)}{|y|^2}|w|^2- \widetilde
  h|w|^2-\tilde f(y,w(y))w\bigg)dy\\
&\notag+ \frac{1}{r^{N-2}}\int_{S_r}
  \bigg( \tildA\nabla w\cdot \nabla w+\tildb\cdot\nabla w\,w-
  \dfrac{V\big(\frac{y}{|y|}\big)}{|y|^2}|w|^2-
  \widetilde h|w|^2-\tilde f(y,w(y))w\bigg)d\sigma\\
  \notag=\,&B(r)+\frac{1}{r^{N-1}}\int_{\Omega_r}\Big((N-2)  \tilde
  f(y,w(y))w-2\nabla_y\widetilde F(y,w)\cdot
  y-2N\widetilde F(y,w)\Big)\,dy\\
  &\notag+\frac{1}{r^{N-2}}\int_{S_r} \bigg( \tildb\cdot\nabla w\,w-
  \widetilde h|w|^2-\tilde
  f(y,w(y))w+2\widetilde F(y,w)
\bigg)d\sigma\\
&\notag+O(r^{1-N+\delta})\bigg( \int_{\Omega_r}\tildA\nabla w\cdot\nabla
  w\,dy +\int_{\Omega_r}\frac{w^2(y)}{|y|^2}
  +\int_{\Omega_r}\Big(|w|^{2}+|w|^{2^*}\Big)\,dy
  \bigg).
\end{align}
From (\ref{eq:5_2}), H\"older inequality, Proposition \ref{SMETS},
and Lemma \ref{welld} $(ii)$,  we have that, for all $r\in(0,r_0)$,
\begin{multline}\label{eq:34}
  \bigg|\frac{1}{r^{N-1}}\int_{\Omega_r}\Big((N-2) \tilde
  f(y,w(y))w-2\nabla_y\widetilde F(y,w)\cdot
  y-2N\widetilde F(y,w)\Big)\,dy\bigg|\\
  \leq \frac{2NC_{\tilde f}}{r^{N-1}}
  \int_{\Omega_r}\big(w^2(y)+|w(y)|^{2^*}\big)\,dy\\
  \leq \frac{2NC_{\tilde f}}{r^{N-1}}
  \left(\Big(\frac{\omega_{N-1}}N\Big)^{\!\frac2N}r^{2}+
    \|w\|_{L^{2^*}(\Omega_r)}^{2^*-2}\right)
  \bigg(\int_{\Omega_r}|w(y)|^{2^*}dy\bigg)^{\!\!2/2^*}\\
  \leq \frac{2NC_{\tilde f}}{\overline{C}} r^{-1+\frac{2(q-2^*)}{q}}
  \left(\Big(\frac{\omega_{N-1}}N\Big)^{\!\frac2N}r_0^{\frac{22^*}q}+
    \Big(\frac{\omega_{N-1}}N\Big)^{\!\frac{2(q-2^*)}{Nq}}
    \|w\|_{L^{q}(\Omega_{r_0})}^{2^*-2} \right)
  \bigg(D(r)+\frac{N-2}{2}H(r) \bigg).
\end{multline}
On the other hand, from (\ref{eq:4}), Lemma \ref{welld} $(i)$, (\ref{eq:H(r)}),
(\ref{eq:5_2}), (\ref{eq:37}), we can estimate
\begin{align}\label{eq:39}
\frac{1}{r^{N-2}}\int_{S_r} \bigg(
  \widetilde h|w|^2+\tilde
  f(y,w(y))w-2\widetilde F(y,w)
\bigg)d\sigma=O\big(g(r)+r^{-1+\delta}\big)  \bigg(D(r)+\frac{N-2}{2}H(r) \bigg).
\end{align}
In view of (\ref{eq:34}), (\ref{eq:39}), and estimate $(ii)$ in Lemma
\ref{welld}, (\ref{eq:32}) yields (\ref{eq:94}).
\end{pf}

\begin{Lemma}\label{l:if_Nprime_neg}
Let $D$ and $H$ be the functions defined in (\ref{eq:D(r)}--\ref{eq:H(r)}),
 $r_0$ be as in Lemma \ref{welld},  and denote
\begin{equation}\label{eq:Sigma}
\Sigma:=\big\{ r\in (0,r_0):D'(r)H(r)\leq H'(r)D(r)\big\}.
\end{equation}
If $\Sigma\not=\emptyset$ and $0$ is a limit point of $\Sigma$, then
\begin{align*}
  D'(r)=B(r)+
  O\Big(r^{-1+\delta}+r^{-1+\frac{2(q-2^*)}{q}}+g(r)\Big)\bigg(D(r)+\frac{N-2}{2}H(r)
  \bigg)
\end{align*}
as $r\to 0^+$, $r\in\Sigma$.
\end{Lemma}
\begin{pf}
From (\ref{eq:poho}),  (\ref{eq:17}--\ref{eq:20}), (\ref{eq:2}--\ref{eq:4}),
Lemma \ref{welld} $(i)$, and (\ref{eq:H(r)})
 we have that
\begin{multline*}
\frac1{r^{N-3}}\int_{S_r} (\tildA\nabla w\cdot\nabla w)\,d\sigma
=rB(r)+O(r^{-N+2})\bigg(\int_{\Omega_r}\!\! \bigg(|\nabla
    w(y)|^2+ \frac{w^2(y)}{|y|^2}\bigg)\,dy+
\|w\|_{L^{2^*}(\Omega_r)}^2\bigg)\\
+O(1)H(r)+O(r^{-N+3})\int_{S_r}|w|^{2^*}\,d\sigma
\end{multline*}
which, in view of Lemma \ref{welld} and (\ref{eq:37}),
implies
\begin{align}\label{eq:41}
\frac1{r^{N-3}}\int_{S_r} (\tildA\nabla w\cdot\nabla w)\,d\sigma
=rB(r)+O(1+rg(r))\bigg(D(r)+\frac{N-2}{2}H(r)
  \bigg)
\end{align}
as $r\to 0^+$. From (\ref{eq:3}), Schwarz inequality,
Lemma \ref{welld} $(i)$,   (\ref{eq:H(r)}), (\ref{eq:2}),
and (\ref{eq:41}), we have that
\begin{align}\label{eq:42}
  \frac{1}{r^{N-2}}\int_{S_r} \tildb\cdot\nabla w\,w\,
  d\sigma&=O(r^{-1+\delta})\bigg(\frac1{r^{N-3}}\int_{S_r} |\nabla
  w|^2\,d\sigma\bigg)^{\!\!1/2}\sqrt{H(r)}\\
  \notag&=O(r^{-1+\delta})\bigg(\frac1{r^{N-3}}\int_{S_r} (\tildA\nabla
  w\cdot\nabla w)\,d\sigma\bigg)^{\!\!1/2}\sqrt{H(r)}\\
  \notag&=O(r^{-1+\delta})\sqrt{rB(r)H(r)}+O\big(r^{-1+\delta}+g(r)\big)
  \bigg(D(r)+\frac{N-2}{2}H(r) \bigg)
\end{align}
as $r\to 0^+$. From Lemma \ref{l:dprime} and (\ref{eq:42}), it follows that
\begin{align*}
  B(r)&=D'(r)-\frac{1}{r^{N-2}}\int_{S_r} \tildb\cdot\nabla w\,w
  d\sigma
  +O\Big(r^{-1+\delta}+r^{-1+\frac{2(q-2^*)}{q}}+g(r)\Big)\bigg(D(r)+\frac{N-2}{2}H(r)
  \bigg)\\
  &=D'(r)+O(r^{-1+\delta})\sqrt{rB(r)H(r)}
  +O\Big(r^{-1+\delta}+r^{-1+\frac{2(q-2^*)}{q}}+g(r)\Big)\bigg(D(r)+\frac{N-2}{2}H(r)
  \bigg)\\
  &\leq D'(r)+\frac{B(r)}2
  +O\Big(r^{-1+\delta}+r^{-1+\frac{2(q-2^*)}{q}}+g(r)\Big)\bigg(D(r)+\frac{N-2}{2}H(r)
  \bigg)
\end{align*}
thus yielding
\begin{align}\label{eq:43}
  B(r)\leq
  2D'(r)+O\Big(r^{-1+\delta}+r^{-1+\frac{2(q-2^*)}{q}}+g(r)\Big)
  \bigg(D(r)+\frac{N-2}{2}H(r) \bigg)
\end{align}
as $r\to 0^+$. From (\ref{eq:43}), (\ref{eq:35}), and the fact that
$D'H\leq H'D$ a.e. in $\Sigma$, we deduce that, as $r\to 0^+$, $r\in
\Sigma$,
\begin{align*}
rB(r)&H(r)\leq
2rH'(r)D(r)+O\Big(r^{\delta}+r^{\frac{2(q-2^*)}{q}}+rg(r)\Big)
  \bigg(D(r)+\frac{N-2}{2}H(r) \bigg)^{\!\!2}\\
&=
2rD(r)\bigg(\frac2rD(r)+H(r)O(r^{-1+\delta})\bigg)
+O\Big(r^{\delta}+r^{\frac{2(q-2^*)}{q}}+rg(r)\Big)
  \bigg(D(r)+\frac{N-2}{2}H(r) \bigg)^{\!\!2}\\
&=4D^2(r)+O(r^\delta)D(r)H(r)+O\Big(r^{\delta}+r^{\frac{2(q-2^*)}{q}}+rg(r)\Big)
  \bigg(D(r)+\frac{N-2}{2}H(r) \bigg)^{\!\!2}\\
&=O\big(1+rg(r)\big)
  \bigg(D(r)+\frac{N-2}{2}H(r) \bigg)^{\!\!2}
\end{align*}
which implies
\begin{align}\label{eq:44}
\sqrt{rB(r)H(r)}=O\big(1+rg(r)\big)
  \bigg(D(r)+\frac{N-2}{2}H(r) \bigg)
\quad\text{as $r\to 0^+$, $r\in \Sigma$}.
\end{align}
Combining  (\ref{eq:42}) and (\ref{eq:44}), we obtain
\begin{align*}
  \frac{1}{r^{N-2}}\int_{S_r} \tildb\cdot\nabla w\,w\,
  d\sigma=O\big(r^{-1+\delta}+g(r)\big)
  \bigg(D(r)+\frac{N-2}{2}H(r) \bigg)
\quad\text{as $r\to 0^+$, $r\in \Sigma$}
\end{align*}
which, together with Lemma \ref{l:dprime}, yields the conclusion.
\end{pf}

In view of Lemma \ref{welld}, the \emph{Almgren type frequency function}
\begin{equation}\label{eq:alm_fun}
{\mathcal N}(r)=\frac{D(r)}{H(r)}
\end{equation}
is well defined in $(0,r_0)$. Furthermore, by Lemmas \ref{l:hprime}
and \ref{l:dprime}, ${\mathcal N}\in W^{1,1}_{\rm loc}(0,r_0)$. The following lemma
provides the existence of a finite limit of ${\mathcal N}(r)$ as $r\to 0^+$.
\begin{Lemma}\label{l:limitN}
Let $\mathcal N:(0,r_0)\to\R$ be defined in (\ref{eq:alm_fun}). Then the limit
$$
\gamma:=\lim_{r\to 0^+}\mathcal N(r)
$$
exists, is finite and
\begin{equation}\label{eq:59}
\gamma\geq -\frac{N-2}{2}.
\end{equation}
\end{Lemma}
\begin{pf}
By Lemma \ref{welld},
\begin{equation}\label{eq:46}
\mathcal N(r)> -\frac{N-2}2
\quad\text{for all }r\in (0,r_0).
\end{equation}
If the set $\Sigma$ defined in (\ref{eq:Sigma}) is empty or if $0$ is
not a limit point of $\Sigma$, then $\mathcal N'(r)\geq 0$ in a right
neighborhood of $0$ and hence $\mathcal N$ is nondecreasing near $0$
and admits a limit as $r\to0^+$ which is necessarily finite in view of
(\ref{eq:46}). If $\Sigma\not=\emptyset$ and $0$ is a limit point of
$\Sigma$, then from Lemma \ref{l:if_Nprime_neg}, (\ref{eq:35}), and
(\ref{eq:B(r)}), we have that
\begin{align}\label{eq:48}
  \mathcal N'(r)&=\frac{D'(r)H(r)-H'(r)D(r)}{H^2(r)}\\
\notag&=
  \frac{B(r)H(r)}{H^2(r)}
-\frac{H'(r)}{H^2(r)}\bigg(\frac{r}2H'(r)+H(r)O(r^{\delta})\bigg)\\
\notag&\quad +
  O\Big(r^{-1+\delta}+r^{-1+\frac{2(q-2^*)}{q}}+g(r)\Big)\bigg(\mathcal
  N(r)+\frac{N-2}{2} \bigg)\\
\notag&=
  \frac{2r\Big(\int_{S_r} \frac{|\tildA \nabla
    w\cdot\nu|^2}\mu\,d\sigma\Big)\Big(\int_{S_r}\mu w^2 d\sigma\Big)
+\Big(\int_{S_r}\mu w^2 d\sigma\Big)\Big(
\int_{\Gamma_r}
  \frac{|\nabla w|^2(\tildA\tilde\nu\cdot\tilde\nu)(\tildA
  y\cdot\tilde\nu)}{\mu}\,d\sigma\Big)}{
\Big(\int_{S_r}\mu w^2 d\sigma\Big)^2
}\\
\notag&\quad
-\frac r2\frac{(H'(r))^2}{H^2(r)}+\frac{H'(r)}{H(r)}O(r^{\delta})+
  O\Big(r^{-1+\delta}+r^{-1+\frac{2(q-2^*)}{q}}+g(r)\Big)\bigg(\mathcal
  N(r)+\frac{N-2}{2} \bigg)
\end{align}
as $r\to 0^+$, $r\in\Sigma$.
In view of (\ref{eq:H'2}), there holds
\begin{multline*}
(H'(r))^2=\frac{4}{r^{2N-2}}\bigg( \int_{S_r}
(\tildA\nabla w\cdot\nu)w\,d\sigma(y)\bigg)^{\!\!2}\\+H^2(r)O(r^{-2+2\delta})+
2H(r)O(r^{-1+\delta})(H'(r)-H(r)O(r^{-1+\delta}))
\end{multline*}
which yields
\begin{align}\label{eq:47}
  \frac{(H'(r))^2}{H^2(r)}=\frac{4\Big( \int_{S_r} (\tildA\nabla
    w\cdot\nu)w\,d\sigma(y)\Big)^{\!2}}{\Big(\int_{S_r}\mu w^2
    d\sigma\Big)^{\!2}} +O(r^{-2+2\delta})+\frac{H'(r)}{H(r)}
  O(r^{-1+\delta}).
\end{align}
Moreover (\ref{eq:35}) implies
\begin{equation}\label{eq:49}
  \frac{H'(r)}{H(r)}=\frac2r\mathcal N(r)+  O(r^{-1+\delta}).
\end{equation}
From (\ref{eq:48}), (\ref{eq:47}), and (\ref{eq:49}), it follows that
\begin{multline}\label{eq:48-bis}
  \mathcal N'(r)=
  \frac{
\int_{\Gamma_r}
  \frac{|\nabla w|^2(\tildA\tilde\nu\cdot\tilde\nu)(\tildA
  y\cdot\tilde\nu)}{\mu}\,d\sigma}{
\int_{S_r}\mu w^2 d\sigma
}\\
+
  \frac{2r
\Big[
\Big(\int_{S_r} \frac{|\tildA \nabla
    w\cdot\nu|^2}\mu\,d\sigma\Big)\Big(\int_{S_r}\mu w^2 d\sigma\Big)
-\Big( \int_{S_r} (\tildA\nabla
    w\cdot\nu)w\,d\sigma(y)\Big)^{\!2}\Big]}{
\Big(\int_{S_r}\mu w^2 d\sigma\Big)^2
}\\
+\frac{H'(r)}{H(r)}O(r^{\delta})+O(r^{-1+2\delta})+
  O\Big(r^{-1+\delta}+r^{-1+\frac{2(q-2^*)}{q}}+g(r)\Big)\bigg(\mathcal
  N(r)+\frac{N-2}{2} \bigg)\\
=
  \frac{
\int_{\Gamma_r}
  \frac{|\nabla w|^2(\tildA\tilde\nu\cdot\tilde\nu)(\tildA
  y\cdot\tilde\nu)}{\mu}\,d\sigma}{
\int_{S_r}\mu w^2 d\sigma
}+
  \frac{2r
\Big[
\Big(\int_{S_r} \frac{|\tildA \nabla
    w\cdot\nu|^2}\mu\,d\sigma\Big)\Big(\int_{S_r}\mu w^2 d\sigma\Big)
-\Big( \int_{S_r} (\tildA\nabla
    w\cdot\nu)w\,d\sigma(y)\Big)^{\!2}\Big]}{
\Big(\int_{S_r}\mu w^2 d\sigma\Big)^2
}\\
+O(r^{-1+\delta})+
  O\Big(r^{-1+\delta}+r^{-1+\frac{2(q-2^*)}{q}}+g(r)\Big)\bigg(\mathcal
  N(r)+\frac{N-2}{2} \bigg)
\end{multline}
as $r\to 0^+$, $r\in\Sigma$.
From (\ref{eq:48-bis}), Lemma \ref{l:tildAytildnu},
and Schwarz inequality, it follows that
$$
  \mathcal N'(r)\geq
  O\Big(r^{-1+\delta}+r^{-1+\frac{2(q-2^*)}{q}}+g(r)\Big)\bigg(\mathcal
  N(r)+\frac{N}{2} \bigg)
$$
as $r\to 0^+$, $r\in\Sigma$. Since $\mathcal N'(r)\geq 0$ a.e. in
$(0,r_0)\setminus\Sigma$, the above inequality is trivially satisfied
as $r\to 0^+$, $r\in(0,r_0)\setminus\Sigma$. Hence there exists some
$r_1\in(0,r_0)$ and $c_1>0$ such that
\begin{equation}\label{eq:40-bis}
\bigg({\mathcal N}+\frac N2\bigg)'(r)\geq
 -c_1\left({\mathcal
      N}(r)+\frac{N}{2}\right)\Big(r^{-1+\delta}
      +r^{-1+\frac{2(q-2^*)}{q}}+g(r)\Big)
\end{equation}
for a.e. $r\in (0,r_1)$.
After integration over $(r,r_1)$ it follows that
\begin{equation*}
{\mathcal N}(r)\leq -\frac N2+
\left({\mathcal
      N}(r_1)+\frac{N}{2}\right)
\exp\left(
c_1\bigg(
\frac{r_1^\delta}{\delta}+ \frac{q}{2(q-2^*)}
r_1^{\frac{2(q-2^*)}{q}}
+\int_0^{ r_1} g(s)\,ds\bigg)\right)
\end{equation*}
for any $r\in (0,r_1)$, thus proving that there exists $c_2>0$ such that
\begin{equation}\label{eq:45}
\mathcal N(r)\leq c_2\quad\text{for all $r\in (0,r_1)$}.
\end{equation}
Estimates (\ref{eq:40-bis}), (\ref{eq:45}), and the fact that $r\mapsto
r^{-1+\delta} +r^{-1+\frac{2(q-2^*)}{q}}+g(r)\in L^1(0,r_1)$ imply
that $\mathcal N'$ is the sum of a nonnegative function and of a
$L^1$-function on $(0,r_1)$.  Therefore
$$
{\mathcal N}(r)={\mathcal N}(r_1)-\int_r^{r_1} {\mathcal N}'(s)\, ds
$$
admits a limit as $r\rightarrow 0^+$ which is necessarily finite in view of
(\ref{eq:45}) and (\ref{eq:46}).
\end{pf}

\begin{Lemma}\label{l:uppb} 
  There exists $r_1\in (0,r_0)$ and $K_1>0$ such that
\begin{equation} \label{1stest}
H(r)\leq K_1 r^{2\gamma}  \quad \text{for all } r\in (0,r_1)
\end{equation}
and
\begin{equation} \label{eq:doubling}
H(2r)\leq K_1 H(r) \quad \text{for all } r\in (0,r_1/2).
\end{equation}
Furthermore, for any $\sigma>0$ there exists a constant
$K_2(\sigma)>0$ depending on $\sigma$ such that
\begin{equation} \label{2ndest} H(r)\geq K_2(\sigma)\,
  r^{2\gamma+\sigma} \quad \text{for all } r\in (0, r_1).
\end{equation}
\end{Lemma}

\begin{pf}
  By (\ref{eq:46}), (\ref{eq:45}), and Lemma \ref{l:limitN}, there
  exists $r_1\in (0,r_0)$ such that ${\mathcal N}$ is bounded in
  $(0,r_1)$ and ${\mathcal N}'\in L^1(0,r_1 )$.  Then from
  (\ref{eq:40}) and (\ref{eq:40-bis}) it follows that
  \begin{equation} \label{qsopra}
{\mathcal N}(r)-\gamma=\int_0^r
    {\mathcal N}'(s) \, ds\geq -c_3 r^{\tilde\delta}
\end{equation}
for some constant $c_3>0$ and all $r\in(0,  r_1)$, where
\begin{equation}\label{eq:60}
\tilde\delta=\min\bigg\{\delta,
\frac{q-2^*}q
\bigg\}.
\end{equation}
Therefore by (\ref{eq:35})
and (\ref{qsopra}) we deduce that, for $r\in(0,r_1)$,
$$
\frac{H'(r)}{H(r)}=\frac{2\,{\mathcal N}(r)}{r}+O(r^{-1+\delta})\geq
\frac{2\gamma}{r}-2c_3 r^{-1+\tilde\delta}+O(r^{-1+\delta})\quad\text{as }r\to0^+,
$$
which, after integration over the interval $(r, r_1)$ and up to
shrinking $r_1$, yields (\ref{1stest}). On the other hand, from boundedness
of $\mathcal N$ in $(0, r_1)$, we have that
$$
\frac{H'(r)}{H(r)}=\frac{2\,{\mathcal
    N}(r)}{r}+O(r^{-1+\delta})\leq\frac{\rm const}{r},
$$
which, for all $r\in(0,r_1/2)$, after integration over the interval
$(r, 2r)$ yields
$$
\log\frac{H(2r)}{H(r)}\leq {\rm const\,}\log2
$$
thus proving (\ref{eq:doubling}).

Let us prove (\ref{2ndest}). Since $\gamma=\lim_{r\rightarrow 0^+}
{\mathcal N}(r)$ and $\frac{H'(r)}{H(r)}-\frac{2\,{\mathcal
    N}(r)}{r}=O(r^{-1+\delta})$, for any $\sigma>0$ there exists
$r_\sigma>0$ such that ${\mathcal N}(r)<\gamma+\sigma/4$ and
$\frac{H'(r)}{H(r)}-\frac{2\,{\mathcal
    N}(r)}{r}\leq\frac\sigma{2r}$
for any
$r\in (0,r_\sigma)$ and hence
$$
\frac{H'(r)}{H(r)}<\frac{2\gamma+\sigma}{r}
\quad \text{for all } r\in (0,r_\sigma).
$$
Integrating over the interval $(r,r_\sigma)$ and by continuity of $H$
outside $0$, we obtain (\ref{2ndest}) for some constant $K_2(\sigma)$
depending on $\sigma$.
\end{pf}

\subsection{The case $N=2$}
The two-dimensional version of Lemma \ref{welld} we are going to prove
in Lemma~\ref{welld-2} requires the following Sobolev type inequality
with boundary terms.

\begin{Proposition} \label{ine:Sobolev}
Let $N\ge 2$ and let $p\in [1,\infty)$ with $p\leq 2^*=2N/(N-2)$ if $N\ge 3$.
Then there exists a constant $C(N,p)>0$ depending only on $N$ and $p$ such that
for all $r>0$
\begin{equation} \label{Sob-ine}
\|v\|_{L^p(B_r)}^2 \leq C(N,p)\, r^{\frac{2N}p+2-N}
\left(\int_{B_r} |\nabla v(x)|^2 dx+\frac{1}{r} \int_{\partial B_r} v^2(x)\, d\sigma
\right) \quad \text{for all } v\in H^1(B_r)
\end{equation}
\end{Proposition}

\begin{pf}
  The proof in the case $r=1$ follows from the classical Sobolev
  inequality and the fact that the square root of the right hand side
  of \eqref{Sob-ine} is a norm equivalent to the standard norm of
  $H^1(B_1)$.  The proof in the case of a general $r>0$ follows by
  scaling.
\end{pf}

\noindent In the rest of this subsection, we assume $N=2$.
\begin{Lemma} \label{welld-2} Let $N=2$ and let $p>2$ as in
  (\ref{eq:F_assumption2}).  Then for every $\e>0$ there exist $\tilde
  r_\e\in(0,\widetilde R)$ and a constant
  $\overline{C}_\e=\overline{C}_\e(\e,p,\tildA,\tildb,\tilde
  f,\widetilde h,w)>0$ depending on $\e$, $p$, $\tildA$, $\tildb$,
  $\tilde f$, $\widetilde h$, $w$ such that, for all $r\in (0,\tilde
  r_\e)$,
\begin{align*}
(i)&\quad \mu(y)>1/2 \quad \text{for all }y\in B_{r},\\
(ii)&\quad D(r)+\e H(r)\geq
    \overline{C}_\e \bigg(\int_{\Omega_r}\!\! |\nabla
    w(y)|^2 dy+\frac{1}r
\int_{S_r} w^2 d\sigma+\|w\|_{L^p(\Omega_r)}^2\bigg),\\
(iii)&\quad  H(r)>0,
\end{align*}
where $D$ and $H$ are defined in (\ref{eq:D(r)}) and (\ref{eq:H(r)}).
\end{Lemma}

\begin{pf}
  The positivity of $\mu$ near $0$ follows from its definition.  By
  \eqref{eq:4}, H\"older inequality, and Proposition
  \ref{ine:Sobolev}, we have
\begin{gather}\label{sp}
  \left|\int_{\Omega_r} \widetilde h(y) w^2(y) dy\right| \leq O(1)
  \int_{\Omega_r} |y|^{-2+\delta} w^2(y)\, dy \\
  \notag\leq O(1)\left(\int_{\Omega_r}
    |y|^{-\frac{4-\delta}2}dy\right)^{\!\!\frac{2(2-\delta)}{4-\delta}}
  \|w\|_{L^{\frac{2(4-\delta)}{\delta}}(\Omega_r)}^2 =O(r^\delta)
  \left( \int_{\Omega_r} |\nabla w(y)|^2 dy+\frac{1}r \int_{S_r} w^2
    d\sigma \right)
\end{gather}
as $r\to 0^+$.
Similarly by \eqref{eq:3}, H\"older inequality, and \eqref{sp} we also have
\begin{align} \label{sp2}
\left| \int_{\Omega_r} \tildb(y)\cdot \nabla w(y) \ w(y)\, dy\right|&
\leq \left(\int_{\Omega_r} |\nabla w(y)|^2 dy\right)^{1/2}
\left(\int_{\Omega_r} |y|^{-2+2\delta} w^2(y)\, dy \right)^{1/2} \\
& \notag
=O(r^\delta) \left( \int_{\Omega_r} |\nabla w(y)|^2 dy+\frac{1}r
\int_{S_r} w^2 d\sigma \right).
 \end{align}
as $r\to 0^+$.
Finally by \eqref{eq:5_2} and Proposition \ref{ine:Sobolev} we have
\begin{align} \label{sp3}
\left| \int_{\Omega_r} \tilde f (y,w(y))w(y) \, dy \right|
&\leq C_{\tilde f} \int_{\Omega_r} (w^2(y)+|w(y)|^p)\, dy\\
\notag &\leq O(r^{\frac{4}p}) \left( \int_{\Omega_r} |\nabla w(y)|^2 dy+\frac{1}r
\int_{S_r} w^2 d\sigma \right)
\end{align}
as $r\to 0^+$.
Inequality (ii) follows from (\ref{sp}-\ref{sp3}).
Inequality (ii) implies  (iii) by proceeding like in the proof
of Lemma \ref{welld}.
\end{pf}

In view of the previous lemma, we can define the Almgren type
frequency function $\mathcal N$ as in the previous subsection, see
(\ref{eq:alm_fun}).  We now sketch the proof of the existence of a
finite limit of $\mathcal N$ as $r\to 0^+$ in dimension $N=2$. To this
aim, we first notice that, under assumption \eqref{eq:V_assumption},
the Pohozaev-type identity \eqref{eq:poho} proved for $N\ge 3$ admits
the following extension to the two-dimensional case.

\begin{Proposition} \label{Poho-2d} Let $N=2$ and let
  $\tildA,\tildb,\tilde f,\widetilde h$ be as in
  (\ref{eq:tildA}--\ref{eq:tildh}) with $A,\bi,\Psi,f,h,V$ as in
  assumptions (\ref{eq:matrix1}--\ref{eq:F_assumption2}),
  (\ref{eq:def_psi}), and let $\tildO$ as in (\ref{eq:tildeomega})
  with $\Omega$ satisfying (\ref{eq:omega}) and
  (\ref{eq:phi1}--\ref{eq:63}).  If $w\in H^1(\tildO)\setminus\{0\}$
  is a weak solution to (\ref{eq:tildequation}), then for a.e.
  $r\in(0,\tilde R)$
\begin{multline}\label{eq:poho-2d}
  r\int_{S_r} (\tildA\nabla w\cdot\nabla w)\,d\sigma - 2r\int_{S_r}
  \frac{|\tildA \nabla w\cdot\nu|^2}\mu\,d\sigma -\int_{\Gamma_r}
  \frac{|\nabla
    w|^2}{\mu}(\tildA\tilde\nu\cdot\tilde\nu)(\tildA y\cdot\tilde\nu)\,d\sigma\\
  = \int_{\Omega_r}(\dive\boldbeta)\tildA\nabla w\cdot\nabla w\,dy
  -2\int_{\Omega_r}(\mathop{\rm Jac}\boldbeta)(\tildA\nabla w)
  \cdot\nabla w\,dy\\
  +\int_{\Omega_r}(d\tildA \nabla w)\nabla w\cdot\boldbeta\,dy-
  \int_{\Omega_r}2(\boldbeta\cdot\nabla w)(\tildb\cdot\nabla w)\,dy\\
  +2 \int_{\Omega_r}(\boldbeta\cdot\nabla w) \widetilde h w\,dy
  -2\int_{\Omega_r} (\nabla_y\widetilde
  F(y,w)\cdot\boldbeta+\widetilde F(y,w)\dive\boldbeta)\,dy + 2
  r\int_{S_r}\widetilde F(y,w)\,d\sigma.
\end{multline}
\end{Proposition}

\begin{pf} It is enough to follow the proof of Proposition
  \ref{p:poho} recalling that $V\equiv 0$ for $N=2$.
\end{pf}

\noindent The next lemma provides an upper bound for a nonlinear boundary term.
\begin{Lemma}\label{l:boundary_non_2d}
  Under the same assumptions of Proposition \ref{Poho-2d}, let $\tilde
  r_1\in (0,\tilde R)$ as in Lemma \ref{welld-2} with $\e=1$. Let
\begin{equation*}
g(r):=\frac{\int_{S_r}|w|^{p}
  d\sigma}{\Big(\int_{\Omega_r}|w(y)|^{p}  dy\Big)^{\!\frac{p+2}{2p}}} .
\end{equation*}
Then $g\in L^1(0,\tilde r_1)$ and $g\geq 0$ a.e. in $(0,\tilde r_1)$.
Furthermore
\begin{equation}\label{eq:37-2d}
\int_{S_r}|w|^{p} d\sigma\leq
\bigg(\int_{\Omega_{\tilde r_1}} |w(y)|^p dy\bigg)^{\!\!\frac{p-2}{2p}}
 \overline C_1^{-1}  g(r)\Big(D(r)+H(r)\Big)
\end{equation}
for a.e. $r\in (0,\tilde r_1)$. Moreover for any $q>p$ and all
$r\in (0,\tilde r_1)$, we have
\begin{equation}\label{eq:40-2d}
\int_0^r g(s)\,ds\leq
\frac{2p}{p-2} \pi^{\frac{(q-p)(p-2)}{2pq}} r^{\frac{(q-p)(p-2)}{pq}}
\|w\|_{L^q(\Omega_{\tilde r_1})}^{\frac{p-2}2}.
\end{equation}
\end{Lemma}
\begin{pf}
We have
\begin{align}\label{eq:38-2d}
  g(r)&= \frac{\int_{S_r}|w|^{p} d\sigma}
{\Big(\int_{\Omega_r}|w(y)|^{p} dy\Big)^{\!\frac{p+2}{2p}}}
=\frac{d}{dr}\bigg(
\frac{2p}{p-2} \bigg(\int_{\Omega_r}|w(y)|^{p} dy\bigg)^{\!\!\frac{p-2}{2p}}\bigg)
\end{align}
in the distributional sense and for a.e. $r\in (0,\tilde r_1)$
and this clearly implies that $g\in L^1(0,\tilde r_1)$. Furthermore,
(\ref{eq:40-2d}) follows from
integration of (\ref{eq:38-2d}) and H\"older inequality.
The proof of \eqref{eq:37-2d} follows by Lemma \ref{welld-2} and
 the definition of $g$.
 \end{pf}

Next we state the two-dimensional version of Lemma \ref{l:dprime}.

\begin{Lemma}\label{l:dprime-2d}
  The function $D$ defined in (\ref{eq:D(r)}) belongs to
  $W^{1,1}_{{\rm loc}}(0,\tilde R)$.
Moreover
\begin{equation*}
D(r)+H(r)\geq 0 \qquad \text{for all } r\in (0,\tilde r_1)
\end{equation*}
where $\tilde r_1\in (0,\tilde R)$ is as in Lemma \ref{welld-2} with $\e=1$, and
\begin{align*}
D'(r)=B(r)+\int_{S_r} \tildb\cdot\nabla w\,w
  d\sigma
  +O\Big(r^{-1+\tilde \delta}+g(r)\Big)\Big(D(r)+H(r)
  \Big)
\end{align*}
as $r\to 0^+$, in a distributional sense and for a.e. $r\in (0,\tilde r_1)$,
where $\tilde \delta=\min\{\delta,4/p\}$ and
\begin{equation*}
B(r):=2\int_{S_r} \frac{|\tildA \nabla
    w\cdot\nu|^2}\mu\,d\sigma +\frac1{r}\int_{\Gamma_r}
  \frac{|\nabla w|^2}{\mu}(\tildA\tilde\nu\cdot\tilde\nu)(\tildA
  y\cdot\tilde\nu)\,d\sigma.
\end{equation*}
\end{Lemma}

\begin{pf} We give here only a sketch of the proof being essentially
  similar to the proof Lemma \ref{l:dprime}.  By \eqref{eq:poho-2d},
  \eqref{eq:D(r)}, \eqref{sp}, \eqref{sp2}, \eqref{sp3},
  (\ref{eq:5_2}), and Lemma \ref{welld-2}, we obtain
\begin{align}\label{eq:32-2d}
  D'(r)=&B(r)+\int_{S_r} \big( \tildb\cdot\nabla w\,w-
  \widetilde h|w|^2-\tilde
  f(y,w)w+2\widetilde F(y,w)
\big)d\sigma\\
&\notag\quad+O\big(r^{-1+\tilde \delta}\big)\Big(D(r)+H(r)\Big),
\end{align}
where $\tilde \delta=\min\{\delta,4/p\}$.  On the other hand, from
(\ref{eq:5_2}), (\ref{eq:4}), Lemma \ref{welld-2} $(i)$,
(\ref{eq:H(r)}), (\ref{eq:37-2d}), we can estimate
\begin{align}\label{eq:39-2d}
\int_{S_r} \bigg(
  \widetilde h|w|^2+\tilde
  f(y,w(y))w-2\widetilde F(y,w)
\bigg)d\sigma=O\big(g(r)+r^{-1+\delta}\big)  \Big(D(r)+H(r) \Big).
\end{align}
In view of (\ref{eq:39-2d}),
(\ref{eq:32-2d}) yields the conclusion.
\end{pf}

\noindent
We now give the statement of two-dimensional version of Lemma
\ref{l:if_Nprime_neg}.
\begin{Lemma}\label{l:if_Nprime_neg-2d}
  Let $D$ and $H$ be defined in
  (\ref{eq:D(r)}--\ref{eq:H(r)}), $\tilde r_1$ be as in Lemma
  \ref{welld-2}, $\tilde \delta$ as in Lemma \ref{l:dprime-2d}, and
  denote $\Sigma:=\big\{ r\in (0,\tilde r_1):D'(r)H(r)\leq H'(r)D(r)\big\}$.
If $\Sigma\not=\emptyset$ and $0$ is a limit point of $\Sigma$, then
\begin{align*}
  D'(r)=B(r)+
  O\Big(r^{-1+\tilde \delta}+g(r)\Big)\big(D(r)+H(r)
  \big)
\quad\text{as $r\to 0^+$, $r\in\Sigma$}.
\end{align*}
\end{Lemma}

\begin{pf}
The proof follows that of  Lemma \ref{l:if_Nprime_neg}
 and exploits   Lemma \ref{l:dprime-2d}.
\end{pf}

\noindent If $\tilde r_1$ is as in Lemma \ref{welld-2} (with $\e=1$),
then the function
\begin{equation} \label{Almgren}
\mathcal N:(0,\tilde r_1)\to\R,\quad
 \mathcal N(r)=\frac{D(r)}{H(r)}
\end{equation}
is well defined.

\begin{Lemma}\label{l:limitN-2d}
  Let $\mathcal N:(0,\tilde r_1)\to\R$ be defined in (\ref{Almgren}).
  Then the limit $\gamma:=\lim_{r\to 0^+}\mathcal N(r)$ exists, is
  finite and
\begin{equation}\label{eq:59-2d}
\gamma\geq 0 .
\end{equation}
\end{Lemma}
 \begin{pf}
   By Lemma \ref{l:hprime} and Lemma \ref{l:dprime-2d} we have that
   $\mathcal N\in W^{1,1}_{{\rm loc}}(0,\tilde r_1)$ and
$\mathcal N(r)\geq -1$ for all $r\in (0,\tilde r_1)$.
 As explained in the proof of Lemma \ref{l:limitN} it is not
 restrictive to assume that $\Sigma\neq \emptyset$ and that $0$ is a
 limit point of $\Sigma$ since otherwise the convergence of $\mathcal
 N$ as $r\to 0^+$ is immediate.  Therefore by Lemma
 \ref{l:if_Nprime_neg-2d} and \eqref{eq:35} we obtain
 \begin{align}\label{eq:48-2d}
  \mathcal N'(r)&= \frac{2r\Big(\int_{S_r} \frac{|\tildA \nabla
    w\cdot\nu|^2}\mu\,d\sigma\Big)\Big(\int_{S_r}\mu w^2 d\sigma\Big)
+\Big(\int_{S_r}\mu w^2 d\sigma\Big)\Big(
\int_{\Gamma_r}
  \frac{|\nabla w|^2(\tildA\tilde\nu\cdot\tilde\nu)(\tildA
  x\cdot\tilde\nu)}{\mu}\,d\sigma\Big)}{
\Big(\int_{S_r}\mu w^2 d\sigma\Big)^2
}\\
\notag&\quad
-\frac r2\frac{(H'(r))^2}{H^2(r)}+\frac{H'(r)}{H(r)}O(r^{\delta})+
  O\Big(r^{-1+\tilde \delta}+g(r)\Big)\Big(\mathcal
  N(r)+1 \Big)
\end{align}
as $r\to 0^+$, $r\in\Sigma$.  Proceeding as in the proof of Lemma
\ref{l:limitN} we arrive to
$$
\mathcal N'(r)\geq
O\Big(r^{-1+\tilde \delta}+g(r)\Big)\bigg(\mathcal
N(r)+\frac{1}2\bigg) .
$$
By Lemma \ref{l:boundary_non_2d} and integration it follows that
$\mathcal N$ is bounded also from above and, in turn, that $\mathcal
N'$ is the sum of a nonnegative function and of a $L^1$-integrable
function in a neighborhood of $0$.  Therefore $\mathcal N$ has a
limit as $r\to 0^+$.  Finally, (\ref{eq:59-2d}) follows immediately
from Lemma \ref{welld-2} (ii).
\end{pf}

\noindent We conclude this subsection with the following estimates on
the function $H$.

\begin{Lemma}\label{l:uppb-2d}
There exists $r_1\in (0,\tilde r_1)$ and  $K_1>0$ such that
\begin{equation} \label{1stest-2d}
H(r)\leq K_1 r^{2\gamma}  \quad \text{for all } r\in (0,r_1)
\end{equation}
and
\begin{equation} \label{eq:doubling-2d}
H(2r)\leq K_1 H(r) \quad \text{for all } r\in (0,r_1/2).
\end{equation}
On the other hand for any $\sigma>0$ there exists a constant
$K_2(\sigma)>0$ depending on $\sigma$ such that
\begin{equation} \label{2ndest-1}
 H(r)\geq K_2(\sigma)\,
  r^{2\gamma+\sigma} \quad \text{for all } r\in (0, r_1).
\end{equation}
\end{Lemma}
\begin{pf} It follows by Lemma \ref{l:limitN-2d} by proceeding exactly
  as in the proof of Lemma \ref{l:uppb}.
\end{pf}

\section{The blow-up argument}\label{sec:blow-up-argument}
Throughout this section, we let $\tildA,\tildb,\tilde f,\widetilde h$
be as in (\ref{eq:tildA}--\ref{eq:tildh}) with $A,\bi,\Psi,f,h,V$ as
in assumptions (\ref{eq:matrix1}--\ref{eq:F_assumption2}),
(\ref{eq:def_psi}). Let $\tildO$ be as in (\ref{eq:tildeomega}) with
$\Omega$ satisfying (\ref{eq:omega}) and
(\ref{eq:phi1}--\ref{eq:phi4}).  Let $w\in H^1(\tildO)\setminus\{0\}$
be a non-trivial weak solution to (\ref{eq:tildequation}).

\begin{Lemma}\label{l:blowup} Let $\gamma$ be as in Lemmas
  \ref{l:limitN}, \ref{l:limitN-2d}
  respectively for $N\ge 3$ and $N=2$. Then
\begin{enumerate}
\item[(i)] there exists $k_0\in \N$ such that
  $\gamma=-\frac{N-2}2+\sqrt{\left(\frac{N-2}2
    \right)^{\!2}+\mu_{k_0}(V)}$;
\item[(ii)] for every sequence $\lambda_n\to 0^+$ there exists a
  subsequence $\lambda_{n_k}$ and $\psi\in H^1_0(C)\subset H^1(\SN)$
  eigenfunction of the operator ${\mathcal L}_V=-\Delta_{\SN}-V$
  associated to the eigenvalue $\mu_{k_0}(V)$ such that
  $\|\psi\|_{L^2(\SN)}=1$ and
$$
\frac{w(\lambda_{n_k}x)}{\sqrt{H(\lambda_{n_k})}}\to |x|^{\gamma}
\psi\Big(\frac{x}{|x|}\Big)
$$
weakly in $H^1(B_1)$, strongly in $C^{1,\alpha}_{{\rm loc}}(\mathcal
C\cap B_1)$
and in $C^{0,\alpha}_{{\rm loc}}(B_1\setminus\{0\})$
for any $\alpha\in (0,1)$, strongly in $H^1(B_r)$ for all
$r\in(0,1)$, and strongly in $L^2(\partial B_1)$, where $w$ is meant to be
trivially extended outside $\tildO$.
\end{enumerate}

\end{Lemma}
\begin{pf}
Let us set
\begin{equation}\label{eq:w_lambda}
w^\lambda(x)=\frac{w(\lambda x)}{\sqrt{H(\lambda)}}.
\end{equation}
We notice that
\begin{equation}\label{eq:normalization}
\int_{C_\lambda}\mu(\lambda \theta)(w^\lambda(\theta))^2 d\sigma=1
\end{equation}
where $C_\lambda$ is defined in (\ref{eq:Cr}).
If $N\ge 3$, by Lemma \ref{welld} we have that, for all $\lambda\in(0,r_0)$,
\begin{align}\label{eq:51}
  \bigg(\mathcal N(\lambda)+\frac{N-2}{2}\bigg)&\geq
  \frac{\overline{C}}{\lambda^{N-2}H(\lambda)}
\int_{\Omega_\lambda}\bigg(|\nabla
    w(x)|^2+ \frac{w^2(x)}{|x|^2}\bigg)\,dx\\
\notag&=  \overline{C}
\int_{\Omega_\lambda/\lambda}\bigg(|\nabla
    w^\lambda(x)|^2+ \frac{(w^\lambda(x))^2}{|x|^2}\bigg)\,dx.
  \end{align}
  Similarly, if $N=2$, by Lemma \ref{welld-2} we have that, for all
  $\lambda\in(0,\tilde r_1)$,
\begin{align}\label{eq:5}
  \big(\mathcal N(\lambda)+1\big)&\geq \frac{\overline{C}_1}{H(\lambda)}
  \bigg(\int_{\Omega_\lambda}|\nabla w(x)|^2dx
  +\frac1\lambda\int_{S_\lambda}w^2\bigg)\\
  &\notag= \overline{C}_1\bigg( \int_{\Omega_\lambda/\lambda}|\nabla
  w^\lambda(x)|^2dx +\int_{C_\lambda}(w^\lambda)^2d\sigma\bigg).
  \end{align}
  From (\ref{eq:45}), Lemma \ref{l:limitN-2d},
  (\ref{eq:51}), and (\ref{eq:5}), we deduce that the trivial extension
$$
\widetilde w^\lambda(x)=\begin{cases}
w^\lambda(x),&\text{if }x\in \Omega_\lambda/\lambda,\\
0,&\text{if }x\in B_1\setminus \Omega_\lambda/\lambda,
\end{cases}
$$
is bounded in $H^1(B_1)$ uniformly with respect to $\lambda\in(0,r_1)$
with $r_1$ as in Lemmas \ref{l:uppb} and \ref{l:uppb-2d}.  Therefore,
for any given sequence $\lambda_n\to 0^+$, there exists a subsequence
$\lambda_{n_k}\to0^+$ such that $\widetilde w^{\lambda_{n_k}}\weakly
\widetilde w$ weakly in $H^1(B_1)$ and a.e. in $B_1$ for some
$\widetilde w\in H^1(B_1)$.  Due to compactness of the trace embedding
$H^1(B_1)\hookrightarrow L^2(\partial B_1)$, we obtain that
$\widetilde w^{\lambda_{n_k}}\to \widetilde w$ in $L^2(\partial B_1)$
and consequently from (\ref{eq:normalization}) $\int_{\partial
  B_1}|\widetilde w|^2d\sigma=1$. In particular $\widetilde
w\not\equiv 0$. Moreover $\widetilde w=0$ a.e. in
$B_1\setminus\mathcal C$ where $\mathcal C$ is defined in
(\ref{eq:cilindro}), as it easily follows from the definition of
$\widetilde w^{\lambda}$, a.e. convergence of $\widetilde
w^{\lambda_{n_k}}\to\widetilde w$ and the fact
\begin{align}\label{eq:conv_comp_lambda}
  \text{for every $x\in B_1\setminus\overline{\mathcal C}$ there
    exists $\lambda_x>0$ such that, for all $\lambda
    \in(0,\lambda_x)$, $x\not\in \Omega_\lambda/\lambda$}.
\end{align}
To prove (\ref{eq:conv_comp_lambda}), it is enough to observe that if
$x=(x',x_N)\in B_1\setminus\overline{\mathcal C}$, then
$x_N<\varphi_0(x')$ and hence from (\ref{eq:23}) $\lambda
x_N<|x'|\widetilde\varphi\big({\lambda x'}/{|x'|}\big)$ for $\lambda$
sufficiently small. From (\ref{eq:9}), we deduce that $\lambda
x/|x'|\not\in \tildO$ for $\lambda$ small which in particular yields
$x\not\in \Omega_\lambda/\lambda$ for $\lambda$ small. This proves claim
(\ref{eq:conv_comp_lambda}).

By scaling of equation (\ref{eq:tildequation}), we have that
$w^\lambda$  weakly solves
\begin{equation}\label{eq:tildequation_lambda}
\begin{cases}
  -\dive(\tildA(\lambda x)\nabla w^\lambda(x))+\lambda \tildb(\lambda
  x)\cdot\nabla w^\lambda(x)-
  \dfrac{V\big(\frac{x}{|x|}\big)}{|x|^2}w^\lambda(x)\\[6pt]
\hskip2cm  =\lambda^2\widetilde h(\lambda x)w^\lambda(x)
  +\frac{\lambda^2}{\sqrt{H(\lambda)}}\tilde f(\lambda
  x,\sqrt{H(\lambda)}w^\lambda (x)),
  &\text{in }\Omega_{\lambda}/\lambda,\\[5pt]
  w^\lambda=0,&\text{on }\partial(\Omega_{\lambda}/\lambda)
\cap B_{1}.
\end{cases}
\end{equation}
In order to pass to the limit in (\ref{eq:tildequation_lambda}), we
observe that
\begin{equation}\label{eq:52}
  \text{if }\xi\in C^{\infty}_{\rm c}(\mathcal C \cap B_1)\quad \text{then}
  \quad \xi\in C^{\infty}_{\rm c}(\Omega_\lambda/\lambda)
  \text{ for sufficiently small }\lambda.
\end{equation}
Indeed, let us consider $\xi\in C^{\infty}_{\rm c}(\mathcal C \cap B_1)$ and
denote $K=\mathop{\rm supp} \xi$. Since $K$ is compact, we have that
$$
\tau=\min_{(x',x_N)\in K}(x_N-\varphi_0(x'))>0.
$$
From (\ref{eq:23}), there exists $t_0$ such that
\begin{equation*}
\Big|
  \frac{\widetilde\varphi(t\nu)}{t}-g(\nu)\Big|<\tau
\quad\text{for all $t\in
(0,t_0)$ and for all $\nu\in {\mathbb S}^{N-2}$}.
\end{equation*}
Then for all $\lambda\in (0,t_0)$ and $(x',x_N)\in K$  we have that
$$
\lambda x_N-\widetilde\varphi(\lambda x')
=\lambda(x_N-\varphi_0(x'))+\lambda |x'|\bigg(g(x'/|x'|)
-\frac{\widetilde\varphi(\lambda|x'|
\frac{x'}{|x'|})}{\lambda |x'|}\bigg)>0
$$
and hence, by (\ref{eq:9}), $K\subset \Omega_\lambda/\lambda$ for all
$\lambda \in (0,t_0)$, thus proving claim (\ref{eq:52}).
Hence we can test equation (\ref{eq:tildequation_lambda}) with  every
$\xi\in C^{\infty}_{\rm c}(\mathcal C \cap B_1)$. From (\ref{eq:2}) we have that
\begin{equation}\label{eq:54}
\int_{\Omega_{\lambda_{n_k}}/\lambda_{n_k}}
\tildA(\lambda_{n_k} x)\nabla w^{\lambda_{n_k}}(x)\cdot\nabla \xi(x)\,dx=
\int_{\mathcal C\cap B_1}
\nabla \widetilde w(x)\cdot\nabla \xi(x)\,dx+o(1)\quad\text{ as }k\to+\infty.
\end{equation}
From (\ref{eq:3}) and (\ref{eq:4})
\begin{equation}\label{eq:55}
\lambda_{n_k}\int_{\Omega_{\lambda_{n_k}}/\lambda_{n_k}}
 \tildb(\lambda_{n_k}
  x)\cdot\nabla w^{\lambda_{n_k}}(x)\xi(x)\,dx
-\lambda_{n_k}^2\int_{\Omega_{\lambda_{n_k}}/\lambda_{n_k}}
\widetilde h(\lambda_{n_k} x)w^{\lambda_{n_k}}(x)\xi(x)\,dx
=o(1)
\end{equation}
as $k\to+\infty$. From (\ref{eq:5_2}) and H\"older and Sobolev
inequalities, we have that, denoting $\tilde p=2^*$ if $N\ge 3$ and
$\tilde p=p$ with $p$ as in (\ref{eq:F_assumption2}) if $N=2$,
\begin{multline}\label{eq:53}
  \frac{\lambda_{n_k}^2}{\sqrt{H(\lambda_{n_k})}}
  \left|\int_{\Omega_{\lambda_{n_k}}/\lambda_{n_k}}\tilde
    f(\lambda_{n_k}
    x,\sqrt{H(\lambda_{n_k})}w^{\lambda_{n_k}}(x))\xi(x)\,dx\right|\\
  \leq C_{\tilde
    f}\lambda_{n_k}^2\int_{\Omega_{\lambda_{n_k}}/\lambda_{n_k}}
  |w^{\lambda_{n_k}}(x)||\xi(x)|\,dx +C_{\tilde
    f}\lambda_{n_k}^2\int_{\Omega_{\lambda_{n_k}}/\lambda_{n_k}}|w(\lambda_{n_k}
  x)|^{\tilde p-2}|w^{\lambda_{n_k}}(x)||\xi(x)|\,dx
  \\
  \leq C_{\tilde f}\lambda_{n_k}^2\|\widetilde
  w^{\lambda_{n_k}}\|_{H^1(B_1)}\|\xi\|_{H^1(B_1)} +C_{\tilde
    f}\lambda_{n_k}^{2-N+\frac{2N}{\tilde p}} \|\widetilde
  w^{\lambda_{n_k}}\|_{L^{\tilde p}(B_1)} \|\xi\|_{L^{\tilde p}(B_1)}
  \|w\|_{L^{\tilde p}(\Omega_{\lambda_{n_k}})}^{\tilde p-2}
  =o(1)\end{multline} as $k\to +\infty$.  Testing equation
(\ref{eq:tildequation_lambda}) with $\xi\in C^{\infty}_{\rm
  c}(\mathcal C \cap B_1)$, letting $k\to+\infty$, and using
(\ref{eq:54}--\ref{eq:53}), we obtain that $\widetilde w$ is a weak
solution to
\begin{equation}\label{eq:eq_limite}
\begin{cases}
  -\Delta \widetilde w-
  \dfrac{V\big(\frac{x}{|x|}\big)}{|x|^2} \widetilde w=0,
\quad&\text{in }\mathcal C\cap B_1,\\
\widetilde w=0, \quad&\text{on }\partial \mathcal C\cap B_1.
\end{cases}
\end{equation}
For
$\lambda \in (0,r_1)$ we define
$$
\Psi_\lambda:\{(y',y_N)\in \R^{N-1}\times \R:|y'|\leq 1\}\to
\{(y',y_N)\in \R^{N-1}\times \R:|y'|\leq 1\}
$$
as
$$
\Psi_\lambda(y',y_N):=\bigg(y', y_N+\frac{\widetilde\varphi(\lambda
  y')}{\lambda}\bigg).
$$
We notice that $\Psi_\lambda$ is invertible and
$\Psi_\lambda^{-1}(x',x_N):= \big(x',
x_N-\frac{\widetilde\varphi(\lambda x')}{\lambda}\big)$.  Let us fix
$r\in (0,1)$, $s_1,s_2,\rho_1,\rho_2,R_1,R_2$ such that
$0<R_1<\rho_1<s_1<r<s_2<\rho_2<R_2<1$, and denote
\begin{align*}
&\mathcal A_{R_1,R_2}:=\Big\{(y',y_N)\in \R^{N-1}\times (0,+\infty):
R_1< \sqrt{
|y'|^2+(y_N+\varphi_0(y'))^2}<R_2\Big\},\\
&\mathcal A_{\rho_1,\rho_2}:=\Big\{(y',y_N)\in \R^{N-1}\times (0,+\infty):
\rho_1< \sqrt{
|y'|^2+(y_N+\varphi_0(y'))^2}<\rho_2\Big\}.
\end{align*}
Using
(\ref{eq:23}) it is easy to verify that there exist
$\lambda_0\in(0,r_1)$ and $c_0>0$ such that for all $\lambda\in (0,\lambda_0)$
\begin{equation*}
|\Psi_\lambda(y)|\geq c_0 \quad\text{for every }y\in  \mathcal A_{R_1,R_2}
\end{equation*}
and
\begin{multline}\label{eq:64}
\Big\{(y',y_N)\in \R^{N-1}\times (0,+\infty):\sqrt{
|y'|^2+\Big(y_N+\lambda^{-1}\widetilde\varphi(\lambda y')\Big)^2}=r\Big\}\\
\subseteq
\Big\{(y',y_N)\in \R^{N-1}\times (0,+\infty):s_1<\sqrt{
|y'|^2+\Big(y_N+\lambda^{-1}\widetilde\varphi(\lambda y')\Big)^2}<s_2\Big\}
\subseteq \mathcal A_{\rho_1,\rho_2}
\\
\subseteq \mathcal A_{R_1,R_2}
\subseteq
\Big\{(y',y_N)\in \R^{N-1}\times (0,+\infty):
\sqrt{
|y'|^2+\Big(y_N+\lambda^{-1}\widetilde\varphi(\lambda y')\Big)^2}
<1\Big\}
\end{multline}
namely
\begin{align}\label{eq:57}
  \Psi_{\lambda}^{-1}\bigg(\frac{\Omega_\lambda}{\lambda}\cap \partial
  B_r\bigg) \subseteq
  \Psi_{\lambda}^{-1}\bigg(\frac{\Omega_\lambda}{\lambda}\cap (
  B_{s_2}\setminus \overline{B}_{s_1})\bigg) \subseteq \mathcal
  A_{\rho_1,\rho_2} \subseteq \mathcal A_{R_1,R_2}\subseteq
  \Psi_{\lambda}^{-1}\bigg(\frac{\Omega_\lambda}{\lambda}\bigg)
\end{align}
for all $\lambda\in(0,\lambda_0)$.
From (\ref{eq:tildequation_lambda}) and (\ref{eq:57}), the functions
$v^\lambda(y):=w^{\lambda}(\Psi_{\lambda}(y))$ satisfy
\begin{align*}
\left\{ \hskip-10pt
\begin{array}{ll}
&  -\dive(\tildA^\lambda(y)\nabla v^\lambda(y))+\ 
  \tildb^\lambda(y)\cdot\nabla v^\lambda(y)-
  \frac{V\big(\frac{\Psi_\lambda(y)}{|\Psi_\lambda(y)|}\big)}
  {|\Psi_\lambda(y)|^2}v^\lambda(y)\\[6pt]
&\qquad = \widetilde
  h^\lambda(y)v^\lambda(y) +\tilde f^\lambda(y,v^\lambda (y)),
  \hskip1.5cm\text{ in }\mathcal A_{R_1,R_2},\\[5pt]
  &v^\lambda=0,\hskip5.5cm\text{ on }\partial \mathcal
  A_{R_1,R_2}\cap\{(y',y_N): y_N=0\},
\end{array}
\right.
\end{align*}
for all $\lambda\in(0,\lambda_0)$, where
\begin{align*}
  &\tildA^\lambda(y)= (\mathop{\rm Jac}\Psi_{\lambda}(y))^{-1}
  \tildA(\lambda \Psi_{\lambda}(y)) ((\mathop{\rm
    Jac}\Psi_{\lambda}(y))^T)^{-1}, \quad \tildb^{\lambda}(y)= \lambda
  \tildb(\lambda\Psi_\lambda(y))
  ((\mathop{\rm Jac}\Psi_\lambda(y))^T)^{-1},\\[2pt]
  & \tilde f^\lambda(y,s) =\frac{\lambda^2}{\sqrt{H(\lambda)}} \tilde
  f(\lambda\Psi_\lambda(y),\sqrt{H(\lambda)}s),\quad \widetilde
  h^\lambda(y)=\lambda^2 \widetilde h(\lambda\Psi_\lambda(y)).
\end{align*}
From (\ref{eq:58}) $\|\tildA^\lambda\|_{W^{1,\infty}(\mathcal
  A_{R_1,R_2})}$ is bounded uniformly with respect to $\lambda\in
(0,\lambda_0)$. From (\ref{eq:3}) and (\ref{eq:4})
$\|\tildb^{\lambda}\|_{L^{\infty}(\mathcal A_{R_1,R_2})}, \|\widetilde
h^\lambda\|_{L^{\infty}(\mathcal A_{R_1,R_2})}$ are bounded uniformly
with respect to $\lambda\in (0,\lambda_0)$.  From (\ref{eq:5_2}),
Lemmas \ref{l:uppb} and \ref{l:uppb-2d}, (\ref{eq:59}), and (\ref{eq:59-2d}),
we have that, denoting again $\tilde p=2^*$ if $N\ge
3$ and $\tilde p=p$ with $p$ as in (\ref{eq:F_assumption2}) if $N=2$,
\begin{align}\label{eq:7}
  \bigg|\frac{\tilde f^\lambda(y,v^\lambda (y))}{v^\lambda (y)}\bigg|
  &\leq C_{\tilde f} \lambda^2 \Big(1+|H(\lambda)|^{(\tilde
    p-2)/2}|v^\lambda
  (y)|^{\tilde p-2}\Big)\\
  &\notag\leq{\rm const\,}\Big(\lambda^2 +\lambda^{(\tilde
    p-2)(\frac2{\tilde p-2}+\gamma)} |v^\lambda (y)|^{\tilde
    p-2}\Big)\leq{\rm const\,}\Big(1+ |v^\lambda (y)|^{\tilde
    p-2}\Big).
\end{align}
Hence, if we define $s=q/(\tilde p-2)>N/2$ with $q$ as in Proposition \ref{SMETS} if 
$N\ge 3$ and $q>\tilde p-2$ if $N=2$, then by \eqref{eq:69} and 
two changes of variables, we obtain
\begin{align*}
\left\|\frac{\tilde f^\lambda(y,v^\lambda (y))}{v^\lambda (y)}\right\|
_{L^s(\mathcal A_{R_1,R_2})} & \leq {\rm const} \left[1+\lambda^2 
|H(\lambda)|^{\frac{\tilde p-2}2}\left(\int_{\mathcal A_{R_1,R_2}} |v^{\lambda}(y)|
^{(\tilde p-2)s} dy \right)^{1/s}  \right]   \\
& \leq
{\rm const} \left[1+\lambda^2 
|H(\lambda)|^{\frac{\tilde p-2}2}\left(\int_{B_1} |\widetilde w^{\lambda}(x)|^q
|{\rm det} \, {\rm Jac}\Psi_{\lambda}^{-1}(x)|\, dx \right)^{1/s}  \right] \\
& \leq {\rm const} \left[1+\lambda^{2-\frac{N}s}  
\left(\int_{\Omega_\lambda} |w(x)|^q dx \right)^{1/s}
\right]=O(1) \qquad \text{as } \lambda\to 0^+ \ .
\end{align*}

Furthermore, up to shrinking
$\lambda_0$, it is easy to verify that $\{v^\lambda\}_{\lambda\in
  (0,\lambda_0)}$ is bounded in $H^1(\mathcal A_{R_1,R_2})$ uniformly
with respect to $\lambda$ and that
$$
\inf_{\substack{\lambda\in (0,\lambda_0)\\
y\in \mathcal A_{R_1,R_2}}}\inf_{\xi\in\R^N\setminus\{0\}}
\frac{\tildA^\lambda(y)\xi\cdot\xi}{|\xi|^2}>0.
$$
Therefore, using classical iterative estimates of Brezis-Kato
\cite{BrezisKato} type (see also Proposition \ref{SMETS}), standard
bootstrap, elliptic regularity theory, (\ref{eq:57}), (\ref{eq:7}),  we first
deduce that
\begin{equation}\label{eq:61}
\{v^\lambda\}_{\lambda\in
  (0,\lambda_0)} \text{ is bounded in $C^{1,\alpha}(\mathcal A_{\rho_1,\rho_2})$ uniformly
with respect to $\lambda$ for all $\alpha\in(0,1)$}.
\end{equation}
From (\ref{eq:61}) and local Lipschitz continuity of
$\widetilde\varphi$, it follows that, for all $x,z\in
\Psi_{\lambda}(\mathcal A_{\rho_1,\rho_2})$,
\begin{align*}
  |w^\lambda(x)-w^\lambda(z)|&=
  |v^\lambda(\Psi_\lambda^{-1}(x))-v^\lambda(\Psi_\lambda^{-1}(z))|
  \leq \|v^\lambda\|_{C^{0,\alpha}(\mathcal A_{\rho_1,\rho_2})}
  |\Psi_\lambda^{-1}(x)-\Psi_\lambda^{-1}(z)|^{\alpha}\\
  &\notag \leq  \|v^\lambda\|_{C^{0,\alpha}(\mathcal
    A_{\rho_1,\rho_2})} \bigg(|x-z|+\frac{|\widetilde\varphi(\lambda
    x')-\widetilde\varphi(\lambda z')|}
  {\lambda}\bigg)^{\alpha}\\
  &\notag \leq {\rm const\,}\|v^\lambda\|_{C^{0,\alpha}(\mathcal
    A_{\rho_1,\rho_2})} |x-z|^{\alpha},
\end{align*}
while from (\ref{eq:61}) and (\ref{eq:58}) we deduce
\begin{align*}
  &|\nabla w^\lambda(x)-\nabla w^\lambda(z)|\leq \big| \big(\nabla
  v^\lambda(\Psi_\lambda^{-1}(x))-\nabla
  v^\lambda(\Psi_\lambda^{-1}(z))\big) \mathop{\rm
    Jac}\Psi^{-1}_\lambda(x)\big|
  \\
  &\notag\hskip3cm+ \big| \nabla
  v^\lambda(\Psi_\lambda^{-1}(z))\big(\mathop{\rm
    Jac}\Psi^{-1}_\lambda(x) -\mathop{\rm
    Jac}\Psi^{-1}_\lambda(z)\big)
  \big|\\
  &\notag \leq {\rm const\,}\Big(\|\nabla
  v^\lambda\|_{C^{0,\alpha}(\mathcal A_{\rho_1,\rho_2})}
  |\Psi_\lambda^{-1}(x)-\Psi_\lambda^{-1}(z)|^{\alpha} +\|\nabla
  v^\lambda\|_{L^{\infty}(\mathcal A_{\rho_1,\rho_2})}
  |\nabla\widetilde\varphi(\lambda x')-\nabla\widetilde\varphi(\lambda
  z')|\Big)
  \\
  &\notag \leq {\rm const\,}|x-z|^{\alpha}.
\end{align*}
In particular, the above estimates yield that
$\|w^\lambda\|_{C^{1,\alpha}(\Psi_{\lambda} (\mathcal
  A_{\rho_1,\rho_2}))}$ is bounded uniformly with respect to
$\lambda\in (0,\lambda_0)$ for all $\alpha\in(0,1)$, and hence, taking
into account (\ref{eq:57}),
\begin{align}\label{eq:65}
  \|w^\lambda\|_{C^{1,\alpha}\big( \frac{\Omega_\lambda}\lambda\cap (
    B_{s_2}\setminus \overline{B}_{s_1})\big)} \text{ is bounded
    uniformly with respect to $\lambda\in (0,\lambda_0)$}
\end{align}
for all $\alpha\in(0,1)$.
(\ref{eq:65}) implies that   $\|\widetilde w^\lambda\|_{C^{0,\alpha}(
    B_{s_2}\setminus \overline{B}_{s_1})}$ is bounded
    uniformly with respect to $\lambda\in (0,\lambda_0)$, and hence
\begin{align}\label{eq:93}
\{\widetilde w^{\lambda}\}_{\lambda\in (0,\lambda_0)}\text{ is relatively compact in }
C^{0,\alpha}(
    B_{s_2}\setminus \overline{B}_{s_1}),
\end{align}
for all $\alpha\in(0,1)$.  From (\ref{eq:64}) and (\ref{eq:65}) it
follows that
\begin{equation} \label{eq:compactness} \{\nabla \widetilde
  w^{\lambda}\}_{\lambda\in (0,\lambda_0)}\text{ is relatively compact
    in } L^{2}(\partial B_r).
\end{equation}
and consequently, from weak convergence $\widetilde
w^{\lambda_{n_k}}\weakly \widetilde w$ in $H^1(B_1)$ we deduce that
for every $r\in (0,1)$
\begin{align}\label{eq:66}
\nabla \widetilde w^{\lambda_{n_k}}\to\nabla \widetilde w\text{ in }
L^{2}(\partial B_r) .
\end{align}
From (\ref{eq:93}), (\ref{eq:65}), and compact embedding of H\"older
spaces, reasoning as in the proof of (\ref{eq:52}), we also obtain
that for all $\alpha\in(0,1)$
\begin{align}\label{eq:67}
\widetilde w^{\lambda_{n_k}}\to\widetilde w\text{ in }
C^{0,\alpha}_{\rm loc}(B_1\setminus\{0\})
\quad\text{and}\quad
w^{\lambda_{n_k}}\to\widetilde w\text{ in }
C^{1,\alpha}_{\rm loc}(\mathcal C \cap B_1) .
\end{align}
Testing equation (\ref{eq:tildequation_lambda}) with $w^\lambda$ and
integrating over $(\Omega_\lambda/\lambda)\cap B_r$ with $r\in (0,1)$,
we obtain
\begin{multline} \label{pl-id}
\int_{(\Omega_\lambda/\lambda)\cap B_r}
\tildA(\lambda x)\nabla w^\lambda(x)\cdot \nabla w^\lambda(x)\,dx
+\lambda \int_{(\Omega_\lambda/\lambda)\cap B_r}
\tildb(\lambda
  x)\cdot\nabla w^\lambda(x)w^\lambda(x)\,dx\\
-
\int_{(\Omega_\lambda/\lambda)\cap B_r}
  \dfrac{V\big(\frac{x}{|x|}\big)}{|x|^2}|w^\lambda(x)|^2\,dx =
\lambda^2\int_{(\Omega_\lambda/\lambda)\cap B_r}
\widetilde h(\lambda x)|w^\lambda(x)|^2\,dx\\
  +\frac{\lambda^2}{\sqrt{H(\lambda)}}\int_{(\Omega_\lambda/\lambda)\cap B_r}
\tilde f(\lambda
  x,\sqrt{H(\lambda)}w^\lambda (x))w^\lambda (x)\,dx
+\int_{(\Omega_\lambda/\lambda)\cap \partial B_r}
\tildA(\lambda x)\nabla w^\lambda(x)\cdot\nu(x)  w^\lambda(x)\,d\sigma(x).
\end{multline}
From (\ref{eq:2}) and boundedness of $\{\widetilde
w^\lambda\}_{\lambda\in (0,\lambda_0)}$ in $H^1(B_1)$ we have that
\begin{align} \label{approx-A}
\int_{(\Omega_\lambda/\lambda)\cap B_r}
\tildA(\lambda x)\nabla w^\lambda(x)\cdot \nabla w^\lambda(x)\,dx
=\int_{(\Omega_\lambda/\lambda)\cap B_r}
|\nabla w^\lambda(x)|^2\,dx+O(\lambda^\delta)
\quad\text{as }\lambda\to 0^+.
\end{align}
From (\ref{eq:3}) and (\ref{eq:4}) we have that
\begin{equation} \label{pl-1} \lambda
  \int_{(\Omega_\lambda/\lambda)\cap B_r} \tildb(\lambda x)\cdot\nabla
  w^\lambda(x)w^\lambda(x)\,dx=O(\lambda^\delta) \quad\text{as
  }\lambda\to 0^+
\end{equation}
and
\begin{equation} \label{pl-2}
\lambda^2\int_{(\Omega_\lambda/\lambda))\cap B_r} \widetilde
h(\lambda x)|w^\lambda(x)|^2\,dx=O(\lambda^\delta) \quad\text{as }\lambda\to
0^+ .
\end{equation}
Proceeding as in \eqref{eq:53} we can prove that
\begin{equation} \label{pl-3}
\frac{\lambda^2}{\sqrt{H(\lambda)}}\int_{(\Omega_\lambda/\lambda)\cap
B_r} \tilde f(\lambda
  x,\sqrt{H(\lambda)}w^\lambda (x))w^\lambda (x)\,dx=o(1)
 \quad\text{as }\lambda\to 0^+ .
\end{equation}
By (\ref{pl-id}--\ref{pl-3}), (\ref{eq:2}), and (\ref{eq:65}), we obtain
\begin{align*}
\int_{(\Omega_\lambda/\lambda)\cap B_r} \bigg(|\nabla
w^\lambda(x)|^2
-  \dfrac{V\big(\frac{x}{|x|}\big)}{|x|^2}|w^\lambda(x)|^2\bigg)\,dx
  = \int_{(\Omega_\lambda/\lambda)\cap \partial B_r}
\frac{\partial w^\lambda }{\partial \nu}w^\lambda\,d\sigma+o(1)
\end{align*}
as $\lambda\to 0^+$,
so that, along the sequence $\lambda_{n_k}$, by \eqref{eq:66}, the
strong convergence $\widetilde w^{\lambda_{n_k}}\to \widetilde w$ in
$L^2(\partial B_r)$, and \eqref{eq:eq_limite}, we obtain
for any positive constant $\widetilde C$
\begin{align*}
  & \lim_{k\to +\infty} \Bigg(
  \int_{(\Omega_{\lambda_{n_k}/\lambda_{n_k}})\cap B_r}\!\!\!\bigg(
  |\nabla w^{\lambda_{n_k}}(x)|^2\,dx-
  \dfrac{V\big(\frac{x}{|x|}\big)}{|x|^2}|w^{\lambda_{n_k}}(x)|^2\bigg)\,dx
  \\ & \notag \qquad \qquad +\frac{\widetilde C}{r}
  \int_{(\Omega_{\lambda_{n_k}}/\lambda_{n_k})\cap \partial B_r}
  |w^{\lambda_{n_k}}(x)|^2 \, d\sigma(x)\Bigg) \\
  & \notag =\int_{{\mathcal C}\cap \partial B_r} \frac{\partial
    \widetilde w}{\partial \nu}(x) \widetilde
  w(x)\,d\sigma(x)+\frac{\widetilde C}{r} \int_{{\mathcal C}\cap
  \partial B_r}
|\widetilde w(x)|^2 \, d\sigma(x) \\
& \notag =\int_{{\mathcal C}\cap B_r} |\nabla \widetilde w(x)|^2
dx-\int_{{\mathcal C}\cap B_r} \frac{V\big(\frac{x}{|x|}\big)}{|x|^2}
|\widetilde w(x)|^2 dx +\frac{\widetilde C}{r} \int_{{\mathcal C}\cap
\partial B_r}
  |\widetilde w(x)|^2 \, d\sigma(x)
\end{align*}
as $k\to +\infty$ and consequently
\begin{multline*}
\lim_{k\to +\infty} \bigg(
\int_{(\Omega_{\lambda_{n_k}/\lambda_{n_k}}\cup\mathcal C)\cap B_r}\!\!\!
\bigg( |\nabla
(\widetilde w^{\lambda_{n_k}}-\widetilde w)
(x)|^2\,dx-\dfrac{V\big(\frac{x}{|x|}\big)}{|x|^2}
(\widetilde w^{\lambda_{n_k}}-\widetilde w)^2(x)\bigg)\,dx\\
+\frac{\widetilde C}{r}
\int_{(\Omega_{\lambda_{n_k}}/\lambda_{n_k}\cup \mathcal C)\cap \partial B_r}
  (\widetilde w^{\lambda_{n_k}}-\widetilde w)^2(x) d\sigma(x)\bigg)=0 ,
\end{multline*}
which, in view of Corollary \ref{c:suppC} with $\lambda r$ in place of
$r$ and through the change of variable $y=\lambda x$ yields the strong
convergence
\begin{equation} \label{H^1-conv}
\widetilde w^{\lambda_{n_k}}\to \widetilde w  \qquad \text{in }
H^1(B_r) .
\end{equation}
According to \eqref{eq:tildequation_lambda} we define the
functions
\begin{align*}
  &D_k(r)\\
  &\notag=\frac1{r^{N-2}}\int_{(\Omega_{\lambda_{n_k}}/\lambda_{n_k})\cap
    B_r }\!\!\! \bigg( \tildA(\lambda_{n_k}y)\nabla
  w^{\lambda_{n_k}}(y)\cdot \nabla
  w^{\lambda_{n_k}}(y)+\lambda_{n_k}\tildb(\lambda_{n_k}y)\cdot\nabla
  w^{\lambda_{n_k}}(y)\,w^{\lambda_{n_k}}(y)\bigg) dy\\
  &\notag\qquad
  -\frac1{r^{N-2}}\int_{(\Omega_{\lambda_{n_k}}/\lambda_{n_k})\cap
    B_r}
  \bigg(\dfrac{V\big(\frac{y}{|y|}\big)}{|y|^2}|w^{\lambda_{n_k}}(y)|^2+
  \lambda_{n_k}^2\widetilde h(\lambda_{n_k}y) |w^{\lambda_{n_k}}(y)|^2\bigg)dy\\
  &\notag
  \qquad-\frac{\lambda_{n_k}^2}{r^{N-2}\sqrt{H(\lambda_{n_k})}}
  \int_{(\Omega_{\lambda_{n_k}}/\lambda_{n_k})\cap B_r} \tilde
  f(\lambda_{n_k}
  y,\sqrt{H(\lambda_{n_k})}w^{\lambda_{n_k}} (y))w^{\lambda_{n_k}} (y)dy,\\
  &H_k(r)=\frac1{r^{N-1}}\int_{(\Omega_{\lambda_{n_k}}/\lambda_{n_k})\cap
\partial B_r
}\mu(\lambda_{n_k}y)|w^{\lambda_{n_k}}(y)|^2\,d\sigma(y) .
\end{align*}
By (\ref{approx-A}--\ref{H^1-conv}) we infer that, for any $r\in
(0,1)$,
\begin{equation} \label{convergences}
D_k(r)\to D_{\widetilde w}(r)  \quad \text{and} \quad H_k(r)\to
H_{\widetilde w}(r)
\end{equation}
as $k\to +\infty$ where
\begin{equation*}
D_{\widetilde w}(r):= \frac{1}{r^{N-2}} \int_{\mathcal C\cap
B_r}\left( |\nabla \widetilde
w(y)|^2-\frac{V\big(\frac{y}{|y|}\big)}{|y|^2} \widetilde
w^2(y)\right) \, dy \quad \text{and} \quad H_{\widetilde
w}(r):=\frac{1}{r^{N-1}} \int_{\mathcal C\cap
\partial B_r} \widetilde w^2 d\sigma .
\end{equation*}
By (\ref{eq:eq_limite}) and \eqref{eq:V_assumption} we have that
$H_{\widetilde w}(r)>0$ for all $r\in (0,1)$.
Therefore the function
\begin{equation}\label{eq:56}
\mathcal N_{\widetilde w}(r):=\frac{D_{\widetilde
w}(r)}{H_{\widetilde w}(r)} \qquad \text{for any }
r\in(0,1)
\end{equation}
is well defined.
Moreover by direct computation one verifies that
\begin{align}\label{eq:62}
\frac{D_k(r)}{H_k(r)}=\mathcal N(\lambda_{n_k} r)
\end{align}
for all $r\in (0,1)$.  By (\ref{convergences}--\ref{eq:62}), Lemmas
\ref{l:limitN} and \ref{l:limitN-2d}, letting $k\to+\infty$, we obtain
\begin{equation} \label{eq:w=gamma}
\mathcal N_{\widetilde w}(r)=\gamma \quad\text{for all }r\in(0,1),
\end{equation}
where $\gamma$ is as in Lemmas \ref{l:limitN} and \ref{l:limitN-2d}.

Proceeding as in Propositions \ref{p:poho} and \ref{Poho-2d}, Lemmas
\ref{l:dprime}, \ref{l:dprime-2d} and \ref{l:HprimeD}, and taking into
account that $\widetilde w$ solves problem \eqref{eq:eq_limite} in the
domain $\mathcal C\cap B_1$, we deduce that $D_{\widetilde w},
H_{\widetilde w}, \mathcal N_{\widetilde w} \in W^{1,1}_{{\rm
    loc}}(0,1)$ and
\begin{align*}
  D'_{\widetilde w}(r)&=2r^{2-N} \int_{\mathcal C\cap \partial B_r}
  \bigg|\frac{\partial \widetilde w}{\partial \nu} \bigg|^2 d\sigma
=2r^{2-N} \int_{\partial B_r}
  \bigg|\frac{\partial \widetilde w}{\partial \nu} \bigg|^2 d\sigma\\
  H'_{\widetilde w}(r)&=2r^{1-N} \int_{\mathcal C\cap \partial B_r}
  \frac{\partial \widetilde w}{\partial \nu}\, \widetilde w \, d\sigma
  =2r^{1-N} \int_{\partial B_r}
  \frac{\partial \widetilde w}{\partial \nu}\, \widetilde w \, d\sigma
  \\
  D_{\widetilde w}(r)&=\frac{r}2 \,H'_{\widetilde w} (r) \ ,
  \\
  \mathcal N'_{\widetilde w}(r)&= \frac{2r\Big[ \Big(\int_{\partial B_r}
\big|\frac{\partial \widetilde w} {\partial
      \nu}\big|^2 d\sigma \Big)\Big(\int_{\partial B_r}
    \widetilde w^2 d\sigma\Big) -\left(\int_{\partial
        B_r} \frac{\partial \widetilde w}{\partial \nu}\, \widetilde w
      \,d\sigma \right)^{\!2}\Big]} {\Big(\int_{\partial
      B_r} \widetilde w^2 d\sigma\Big)^{\!2}}
\end{align*}
for a.e. $r\in (0,1)$.  On the other hand, by (\ref{eq:w=gamma}),
$\mathcal N_{\widetilde w}$ is a constant function thus implying
\begin{equation*}
 \bigg(\int_{\partial B_r} \bigg|\frac{\partial \widetilde w}
{\partial \nu}\bigg|^2 d\sigma \bigg)\bigg(\int_{\partial B_r}
\widetilde w^2 d\sigma\bigg)
-\bigg(\int_{\partial B_r} \frac{\partial \widetilde w}{\partial \nu}\,
\widetilde w\, d\sigma \bigg)^{\!\!2}=0.
\end{equation*}
The above identity shows that the functions $\frac{\partial\widetilde
  w}{\partial \nu}$ and $\widetilde w$ have the same direction as
vectors in $L^2(\partial B_r)$ and hence there exists a function
$\eta=\eta(r)$ such that
\begin{equation} \label{eq:ode-w} \frac{\partial \widetilde
    w}{\partial \nu}(r,\theta)=\eta(r)\widetilde w(r,\theta) \qquad
  \text{for a.e. } r\in(0,1),\ \theta\in \SN  .
\end{equation}
Since necessarily  $\eta(r)=\frac{H'_{\widetilde
    w}(r)}{2H_{\widetilde w}(r)}$, then $\eta\in L^1_{{\rm
    loc}}(0,1)$.  Integration of \eqref{eq:ode-w} yields
\begin{equation*}
  \widetilde w(r,\theta)=e^{\int_1^r \eta(s)\, ds} \widetilde w(1,\theta)=
  \varphi(r)\psi(\theta)
  \qquad \text{for all } r\in (0,1),\  \theta\in \SN \ ,
\end{equation*}
where $\varphi(r)=e^{\int_1^r \eta(s)\, ds}$ and
$\psi(\theta)=\widetilde w(1,\theta)$.  We notice that $\psi\in
H^1_0(C)$ and \eqref{eq:eq_limite} may be written in polar
coordinates as
\begin{equation*}
\left(-\varphi''(r)-\frac{N-1}r \, \varphi'(r)\right) \psi(\theta)
-\frac{\varphi(r)}{r^2}\,  {\mathcal L}_V \psi(\theta)=0
\quad\text{in }(0,1)\times C.
\end{equation*}
Taking $r$ fixed we may observe that $\psi$ has to be necessarily an
eigenfunction of the operator ${\mathcal L}_V$ on $C\subset\SN$ with
homogeneous Dirichlet boundary conditions. Hence, if we denote by $\mu_{k_0}(V)$
the corresponding eigenvalue, it follows that $\varphi$ solves the
equation
\begin{equation*}
  -\varphi''(r)-\frac{N-1}r \, \varphi'(r)+\frac{\mu_{k_0}(V)}{r^2} \,
  \varphi(r)=0 \quad\text{in }(0,1) .
\end{equation*}
The general solution of the above equation is given by
\begin{equation*}
\varphi(r)=c_1 r^{\sigma_{k_0}^+}+c_2 r^{\sigma_{k_0}^-}, \quad 
c_1,c_2\in\R,
\end{equation*}
where
$\sigma_{k_0}^{\pm}=-\frac{N-2}2\pm\sqrt{\left(\frac{N-2}2\right)^2+\mu_{k_0}(V)}$.
We observe that the function $|x|^{\sigma_{k_0}^-}
\psi(\frac{x}{|x|})\notin H^1(B_1)$ and hence $c_2=0$. Moreover
$\varphi(1)=1$ implies $c_1=1$, so that $\widetilde w$ takes the form
$\widetilde w(r,\theta)=r^{\sigma_{k_0}^+} \psi(\theta)$.
Finally, inserting this representation of $\widetilde w$ in $\mathcal
N_{\widetilde w}$ and taking into account that
$$
\mathcal N_{\widetilde w}(r)=\frac{D_{\widetilde w}(r)}{H_{\widetilde
    w}(r)} =\frac{r}2 \frac{H'_{\widetilde w}(r)}{H_{\widetilde w}(r)}
= \frac{r\int_{\mathcal C\cap \partial B_r} \frac{\partial
    \widetilde w}{\partial \nu} \widetilde w\, d\sigma}{\int_{\mathcal
    C\cap \partial B_r} \widetilde w^2 \, d\sigma },
$$
from \eqref{eq:w=gamma} it follows that $\sigma_{k_0}^+=\gamma$.
The proof is thereby complete.
\end{pf}

\noindent We  now study the behavior of $H(\lambda)$ as $\lambda\to0^+$.

\begin{Lemma} \label{l:lim-exist}
Let $\gamma$ as in Lemmas \ref{l:limitN} and \ref{l:limitN-2d}. Then
\begin{equation*}
\lim_{\lambda\to 0^+} \lambda^{-2\gamma} H(\lambda)
\end{equation*}
exists and is finite.
\end{Lemma}

\begin{pf} In view of (\ref{1stest}), \eqref{1stest-2d}, Lemma
  \ref{welld}(i) and Lemma \ref{welld-2}(i), it is sufficient to
  prove that the limit exists. By  (\ref{eq:35}) and
  Lemmas \ref{l:limitN}, \ref{l:limitN-2d} we have
\begin{align}\label{eq:d/dr}
\frac{d}{dr} \frac{H(r)}{r^{2\gamma}}& =-2\gamma r^{-2\gamma-1}
H(r)+r^{-2\gamma} H'(r)\\
\notag& =2r^{-2\gamma-1} (D(r)-\gamma
H(r)+H(r)O(r^\delta))=2r^{-2\gamma-1} H(r) \bigg(\int_0^r {\mathcal N}'(s) ds
+O(r^\delta)\bigg).
\end{align}
Let us define the functions
\begin{align*}
&\nu_1(r):=\frac{
\int_{\Gamma_r}
  \frac{|\nabla w|^2(\tildA\tilde\nu\cdot\tilde\nu)(\tildA
  x\cdot\tilde\nu)}{\mu}\,d\sigma}
{\int_{S_r}\mu w^2 d\sigma}\\
&\qquad\qquad\qquad+\frac{2r \Big[
\Big(\int_{S_r} \frac{|\tildA \nabla
    w\cdot\nu|^2}\mu\,d\sigma\Big)\Big(\int_{S_r}\mu w^2 d\sigma\Big)
-\Big( \int_{S_r} (\tildA\nabla
    w\cdot\nu)w\,d\sigma(y)\Big)^{\!2}\Big]}{
\Big(\int_{S_r}\mu w^2 d\sigma\Big)^2} \\
&\nu_2(r):=\mathcal N'(r)-\nu_1(r)   .
\end{align*}
By Lemma \ref{l:tildAytildnu} and Schwartz inequality we have that
$\nu_1\ge 0$. On the other hand from the proofs of Lemmas
\ref{l:limitN}, \ref{l:limitN-2d}, \ref{l:boundary_non},
\ref{l:boundary_non_2d}, we infer that the function $\nu_2$ is
$L^1$-integrable in a right neighborhood of zero and moreover
\begin{equation} \label{eq:int-nu2}
\int_0^r \nu_2(s)\, ds=O(r^{\hat\delta}), \qquad \text{as } r\to 0^+
\end{equation}
where
\begin{align}\label{eq:tilde_delta}
\hat\delta=
\begin{cases}
\min\{\delta,
(q-2^*)/q\},&\text{if }N\geq 3,\\
\min\left\{\delta,\frac 4p,\frac{p-2}{2p} \right\},&\text{if }N=2,
\end{cases}
\end{align}
with $p$ as in (\ref{eq:F_assumption2}) and $q$ as in (\ref{eq:50}).
After integration of \eqref{eq:d/dr} over the interval $(r,r_1)$ we obtain
\begin{align}\label{inte}
\frac{H(r_1)}{r_1^{2\gamma}}-\frac{H(r)}{r^{2\gamma}}&=
\int_r^{r_1} 2s^{-2\gamma-1}
  H(s) \left( \int_0^s \nu_1(t) dt \right) ds +\int_r^{r_1} 2s^{-2\gamma-1}
  H(s) \left( \int_0^s \nu_2(t) dt \right) ds  \\
&\notag\quad +O\bigg( \int_r^{r_1} s^{-2\gamma-1+\delta}
  H(s)\,ds\bigg).
\end{align}
From the nonnegativity of $\nu_1$ the limit
$\lim_{r\to 0^+} \int_r^{r_1} 2s^{-2\gamma-1} H(s) \left( \int_0^s
  \nu_1(t) dt \right) ds$
exists.
On the other hand, by \eqref{eq:int-nu2}, \eqref{1stest}, \eqref{1stest-2d}
\begin{equation*}
  \left| s^{-2\gamma-1} H(s) \left( \int_0^s \nu_2(t) dt \right)
  \right|\leq
  O(1) s^{-1} \left|\int_0^s \nu_2(t) \, dt\right|=O(s^{-1+\hat\delta})
  \quad \text{as } s\to 0^+
\end{equation*}
which proves that $s^{-2\gamma-1} H(s) \left( \int_0^s \nu_2(t) dt
\right)$ is $L^1$-integrable in a right neighborhood of the origin.
From \eqref{1stest} and \eqref{1stest-2d} it follows that
$s^{-2\gamma-1+\delta} H(s)=O(s^{-1+\delta})$ as $s\to 0^+$ and hence
$s^{-2\gamma-1+\delta} H(s)$ is $L^1$-integrable in a right
neighborhood of the origin.  We may therefore conclude that all terms
in the right hand side of (\ref{inte}) admit a limit as $r\to 0^+$
thus completing the proof of the lemma.~\end{pf}

\section{Straightening the domain}\label{sec:straightening-domain}

\begin{Lemma}\label{l:straightening}
There exists $\widehat R\in(0,\widetilde R)$ such that the function
\begin{align}\label{eq:70}
  &\Xi:
\tildO\cap B_{\widehat R}
\to
\mathcal C\cap B_{\widehat R},
\\
\notag&\Xi(y)=\Xi(y',y_N)= \frac{ (y',
    y_N-\widetilde\varphi(y')+\varphi_0(y'))}{ \sqrt{ 1+\dfrac{
        (\varphi_0(y')-\widetilde\varphi(y'))^2
        +2y_N(\varphi_0(y')-\widetilde\varphi(y'))} { |y'|^2+y_N^2 }
    }}
\end{align}
is invertible. Furthermore, putting $\Phi=\Xi^{-1}$, we have
\begin{align}
  &\label{eq:76}\Phi\in C^1(\mathcal C\cap B_{\widehat R},\tildO\cap
  B_{\widehat R}), \quad
  \Phi^{-1}\in C^1(\tildO\cap B_{\widehat R},\mathcal C\cap B_{\widehat R}),\\
  &\label{eq:77}\Phi(\mathcal C\cap\partial B_r)=\tildO\cap\partial
  B_r
  \quad\text{for all }r\in(0,\widehat R),\\
  &\label{eq:78}\Phi(x)=x+O(|x|^{1+\delta})\quad\text{and}\quad
  \mathop{\rm
    Jac}\Phi(x)={\rm Id}_N+O(|x|^{\delta}) \quad \text{as }|x|\to0,\\
  &\label{eq:79}\Phi^{-1}(y)=y+O(|y|^{1+\delta}) \quad\text{and}\quad
  \mathop{\rm Jac}\Phi^{-1}(y)={\rm Id}_N+O(|y|^{\delta}) \quad
  \text{as
  }|y|\to0,\\
  &\label{eq:80}\mathop{\rm det}\mathop{\rm
    Jac}\Phi(x)=1+O(|x|^\delta)\quad \text{as }|x|\to0.
\end{align}
\end{Lemma}
\begin{pf}
  It follows from the Local Inversion Theorem, 
  (\ref{eq:68}--\ref{eq:69}), and direct calculations.~\end{pf}

\noindent Let $u\in H^1(\Omega)$ be a weak solution to
(\ref{eq:equation}), so that $w=u\circ \Psi\in H^1(\tildO)$  weakly
solves (\ref{eq:tildequation}).  Then
\begin{equation}\label{eq:74}
v=w\circ \Phi=u\circ\Psi\circ\Phi \in
H^1(\mathcal C\cap B_{\widehat R})
\end{equation}
is a weak solution to 
\begin{equation}\label{eq:hatquation}
\begin{cases}
  -\dive(\widehat A(x)\nabla v(x))+\widehat {\mathbf b}(x)\cdot\nabla
  v(x)- \dfrac{V\big(\frac{x}{|x|}\big)}{|x|^2}v(x) = \widehat
  h(x)v(x)+\widehat f(x,v(x)),&\text{in }
  \mathcal C\cap B_{\widehat R},\\[5pt]
  v=0,&\text{on }\partial\mathcal C\cap B_{\widehat R},
\end{cases}
\end{equation}
where
\begin{align*}
  &\widehat A(x)= |\mathop{\rm det}\mathop{\rm Jac}\Phi(x)|
  (\mathop{\rm Jac}\Phi(x))^{-1}\tildA(\Phi(x)) ((\mathop{\rm
    Jac}\Phi(x))^T)^{-1}, \\[2pt]
  &\widehat {\mathbf b}(x)=|\mathop{\rm det}\mathop{\rm Jac}\Phi(x)|
  \tildb(\Phi(x))((\mathop{\rm Jac}\Phi(x))^T)^{-1},\\[2pt]
  & \widehat f(x,s)=|\mathop{\rm det}\mathop{\rm Jac}\Phi(x)|
  \tilde f(\Phi(x),s),\\[-5pt]
  & \widehat h(x)=|\!\mathop{\rm det}\!\mathop{\rm Jac}\Phi(x)| \widetilde
  h(\Phi(x))+ |\!\mathop{\rm det}\!\mathop{\rm Jac}\Phi(x)|
  \bigg(\tfrac{V(\frac{\Phi(x)}{|\Phi(x)|})}{|\Phi(x)|^2}-
  \tfrac{V(\tfrac{x}{|x|})}{|x|^2}\bigg)+ \big( |\!\mathop{\rm
    det}\!\mathop{\rm
    Jac}\Phi(x)|-1\big)\tfrac{V(\frac{x}{|x|})}{|x|^2}.
\end{align*}
By Lemmas \ref{l:straightening}, \ref{l:new_coeff}, and direct
calculations, we obtain
\begin{align}
  \label{eq:71}&
\widehat A(x)={\rm Id}_N+O(|x|^{\delta})\quad\text{as }|x|\to0,\\
\label{eq:72}&\widehat \bi\in L^\infty_{\rm loc}(\mathcal C\cap B_{\widehat R},\R^N),
\quad
|\widehat \bi(x)|=O(|x|^{-1+\delta})\quad\text{as }|x|\to0,\\
\label{eq:73}
&\widehat f\in C^0((\mathcal C\cap B_{\widehat R})\times \R)
\quad\text{and}\quad |\widehat f(x,s)s|\leq
\begin{cases}
C_{\widehat f}\,(|s|^2+|s|^{2^*}),&\text{if }N\geq 3,\\
C_{\widehat f}\,(|s|^2+|s|^{p}),&\text{if }N=2,
\end{cases}\\
\label{eq:75}&\widehat h \in L^\infty_{\rm loc}(\mathcal C\cap
B_{\widehat R}),\quad \widehat h(x)=O(|x|^{-2+\delta})\quad\text{as
}|x|\to0.
\end{align}
\begin{Lemma}\label{l:Hstraight}
Let $H$ be as in (\ref{eq:H(r)}) and $v=w\circ \Phi$ as in (\ref{eq:74}). Then
\begin{align}
\label{eq:81}
&H(\lambda)=\big(1+O(\lambda^\delta)\big)
\int_Cv^2(\lambda\theta)\,d\sigma(\theta),\\
\label{eq:82}
&\frac{\int_{\mathcal C\cap B_1} |\nabla \hat
  v^\lambda(x)|^2dx}{{H(\lambda)}}
=(1+O(\lambda^\delta))\int_{\Omega_\lambda/\lambda}|\nabla
w^\lambda(y)|^2\,dy =O(1),
\end{align}
as $\lambda\to 0^+$, where $w^\lambda$ is defined in (\ref{eq:w_lambda}) and
$\hat v^\lambda(x):=v(\lambda x)$.
\end{Lemma}
\begin{pf}
From (\ref{eq:77}), by a change of variable
$$
H(\lambda)=\int_C\mu(\Phi(\lambda\theta))v^2(\lambda\theta)
|\mathop{\rm det}\mathop{\rm Jac}\Phi(\lambda\theta)|\,d\sigma(\theta)
$$
and
$$
\frac{\int_{\mathcal C\cap B_1} |\nabla \hat
  v^\lambda(x)|^2dx}{{H(\lambda)}}
=\int_{\Omega_\lambda/\lambda}|\nabla w^\lambda(y) \mathop{\rm
  Jac}\Phi(\Phi^{-1}(\lambda y))|^2 |\mathop{\rm det}\mathop{\rm
  Jac}\Phi^{-1}(\lambda y)|\,dy
$$
for all $\lambda\in (0,\widehat R)$. We conclude from
(\ref{eq:11}), (\ref{eq:78}--\ref{eq:80}), and
$H^1$-boundedness of $\{w^\lambda\}$ (see the proof of Lemma
\ref{l:blowup}).
\end{pf}

\begin{Lemma}\label{l:blow_v}
  Let $v=w\circ \Phi$ be as in (\ref{eq:74}) and let $k_0$ and
  $\gamma$ as in Lemma \ref{l:blowup}(i). Then for every sequence
  $\lambda_n\to 0^+$ there exist a subsequence $\lambda_{n_k}$ and
  $\psi\in H^1_0(C)\subset H^1(\SN)$ eigenfunction of the operator
  ${\mathcal L}_V=-\Delta_{\SN}-V$ associated to the eigenvalue
  $\mu_{k_0}(V)$ such that $\|\psi\|_{L^2(\SN)}=1$, the convergences
  of $\frac{w(\lambda_{n_k}x)}{\sqrt{H(\lambda_{n_k})}}$ to
  $|x|^{\gamma} \psi\big(\frac{x}{|x|}\big)$ stated in part ii) of
  Lemma \ref{l:blowup} hold, and
$$
\frac{v(\lambda_{n_k}\cdot)}{\sqrt{
\int_Cv^2(\lambda_{n_k}\theta)\,d\sigma(\theta)}}\to \psi
$$
strongly in $L^2(C)$.
\end{Lemma}
\begin{pf}
  From Lemma \ref{l:blowup}, there exist a subsequence $\lambda_{n_k}$
  and $\psi\in H^1_0(C)\subset H^1(\SN)$ eigenfunction of the operator
  ${\mathcal L}_V=-\Delta_{\SN}-V$ associated to the eigenvalue
  $\mu_{k_0}(V)$ such that $\|\psi\|_{L^2(\SN)}=1$ and
  $(H(\lambda_{n_k}))^{-1/2}w(\lambda_{n_k}x)\to |x|^{\gamma}
  \psi\big(\frac{x}{|x|}\big)$ in senses claimed in part ii) of
  Lemma~\ref{l:blowup}, in particular strongly in $L^2(\SN)$ and a.e.
  on $\SN$. Moreover from $H^1$-boundedness of $w^\lambda$ (see the
  proof of Lemma \ref{l:blowup}) and (\ref{eq:82}) it follows that
  $\{\hat v^\lambda/\sqrt{H(\lambda)}\}_{\lambda}$ is bounded in
  $H^1(\mathcal C\cap B_1)$ and relatively compact in $L^2(C)$. Hence
  from (\ref{eq:81}) there exists $\widetilde\psi\in L^2(C)$ such
  that, up to a further subsequence,
$$
\frac{v(\lambda_{n_k}\cdot)}
{\sqrt{\int_Cv^2(\lambda_{n_k}\theta)\,d\sigma(\theta)}}\to
\widetilde\psi\quad\text{in } L^2(C)\text{ and a.e.}
$$
From (\ref{eq:65}) together with (\ref{eq:doubling}) and
(\ref{eq:doubling-2d}) which allow extending estimate (\ref{eq:65}) up
to $\partial B_1$, we have that, in view of (\ref{eq:78}) and
(\ref{eq:81}), for a.e. $\theta\in \SN,$
\begin{align*}
  \frac{v(\lambda_{n_k}\theta)}
  {\sqrt{\int_Cv^2(\lambda_{n_k}\theta)\,d\sigma(\theta)}}&=
  (1+O(\lambda_{n_k}^\delta))\frac{w(\Phi(\lambda_{n_k}\theta))}
  {\sqrt{H(\lambda_{n_k})}}\\
  &=(1+O(\lambda_{n_k}^\delta))\bigg[\frac{w(\lambda_{n_k}\theta)}
  {\sqrt{H(\lambda_{n_k})}}
  +\bigg(w^{\lambda_{n_k}}\Big(\frac{\Phi(\lambda_{n_k}\theta)}{\lambda_{n_k}}\Big)-
  w^{\lambda_{n_k}}(\theta)\bigg)\bigg]\\
  &=\frac{w(\lambda_{n_k}\theta)}{\sqrt{H(\lambda_{n_k})}}
  +O(\lambda_{n_k}^\delta)\to \psi (\theta)\quad\text{as }k\to+\infty.
\end{align*}
Then $\widetilde\psi=\psi$ and the lemma is proved.
\end{pf}

In the sequel we denote by $\psi_i$ a $L^2$-normalized eigenfunction
of the operator ${\mathcal L}_V=-\Delta_{\SN}-V$
on the spherical cap
$C\subset\mathbb S^{N-1}$ under null Dirichlet boundary conditions
associated to the $i$-th eigenvalue
$\mu_i(V)$, i.e.
\begin{equation} \label{eq:2rad}
\begin{cases}
{\mathcal L}_V\psi_i(\theta)
=\mu_i(V) \,\psi_i(\theta),&\text{in }C,\\[3pt]
\psi_i=0,&\text{on }\partial C,\\[3pt]
\int_{{\mathbb S}^{N-1}}|\psi_i(\theta)|^2\,d\sigma(\theta)=1.
\end{cases}
\end{equation}
Moreover, we  choose the $\psi_i$'s in such a way that
 the set $\{\psi_i\}_{i\in \N\setminus\{0\}}$
forms an orthonormal basis of $L^2(C)$.
For all $i\in\N,i\geq1$, and $\lambda\in(0,\widehat R)$, we also define
\begin{equation} \label{78}
  \varphi_i(\lambda):=\int_{C}v(\lambda\,\theta)
  \psi_i(\theta)\,d\sigma(\theta).
\end{equation}
From Lemma \ref{l:blowup}, there exist $j_0,m\in\N$, $j_0,m\geq 1$
such that $m$ is the multiplicity of the eigenvalue
$\mu_{j_0}(V)=\mu_{j_0+1}(V)=\cdots=\mu_{j_0+m-1}(V)$ and
\begin{equation}\label{eq:83}
  \gamma=\lim_{r\rightarrow 0^+} {\mathcal N}(r)=
-\frac{N-2}{2}+\sqrt{\bigg(\frac{N-2}
    {2}\bigg)^{\!\!2}+\mu_{i}(V)},
  \quad i=j_0,\dots,j_0+m-1.
\end{equation}
Let ${\mathcal E}_0$ be the eigenspace
of the operator $\mathcal L_{V}$ associated to the eigenvalue
$\mu_{j_0}(V)$, so that the set $\{\psi_i\}_{i=j_0,\dots,j_0+m-1}$ is an
 orthonormal basis of ${\mathcal E}_0$.

 \begin{Lemma}\label{l:phii}
   Let $v=w\circ \Phi$ be as in (\ref{eq:74}), $j_0$ and $m$ as in
   (\ref{eq:83}) and $\varphi_i$ as in (\ref{78}). Then for all
   $i\in \{j_0,\dots,j_0+m-1\}$ and $R\in (0,\widehat R)$
   \begin{multline} \label{eq:id} \varphi_i(\lambda)= \lambda^{\gamma}
     \bigg( R^{-\gamma}\varphi_i(R)+ \frac{2-N-\gamma}{2-N-2\gamma}
     \int_\lambda^Rs^{-N+1-\gamma}\Upsilon_i(s)ds \\
     -\frac{\gamma R^{-N+2-2\gamma}}{2-N-2\gamma} \int_0^R
     s^{\gamma-1} \Upsilon_i(s)\,ds \bigg)+O(\lambda^{\gamma+\hat\delta})
\end{multline}
as $\lambda\to 0^+$
with $\hat\delta$ as in (\ref{eq:tilde_delta})
and $\Upsilon_i\in L^1(0,\widehat R)$ defined as
\begin{multline} \label{Upsilon}
  \Upsilon_i(\lambda)=-\int_{\mathcal C\cap B_\lambda} (\widehat
  A(x)-{\rm Id}_N)\nabla
  v(x)\cdot\frac{\nabla_{\SN}\psi_i(x/|x|)}{|x|}\,dx\\
+\int_{\mathcal C\cap B_\lambda}
\Big(
 - \widehat {\mathbf
    b}(x)\cdot\nabla v(x) + \widehat h(x)
  v(x)+\widehat f(x,v(x))\Big)\psi_i(x/|x|)\,dx\\
  +\int_{\mathcal C\cap \partial B_\lambda} (\widehat A(x)-{\rm Id}_N)
\nabla v(x)\cdot \frac{x}{|x|}\psi_i(x/|x|)\,d\sigma(x).
\end{multline}
\end{Lemma}

\begin{pf}
For any $\lambda\in (0,\widehat R)$, we
expand $\theta\mapsto v(\lambda\theta)\in L^2(C)$
 in  Fourier series with respect to the orthonormal
basis $\{\psi_i\}$ of $L^2({\mathbb S}^{N-1})$ defined in \eqref{eq:2rad}, i.e.
\begin{equation}\label{77}
v(\lambda\,\theta)=\sum_{i=1}^\infty\varphi_i(\lambda)\psi_i(\theta)
\quad \text{in }L^2(C),
\end{equation}
with $\varphi_i$ is defined in (\ref{78}).
For all $i$, we consider the distribution $\zeta_i$ on $(0,\widehat R)$
defined as
\begin{multline*}
  {}_{\mathcal D'(0,\widehat R)}\langle \zeta_i,\omega\rangle_{\mathcal D(0,\widehat R)}\\
  = {}_{H^{-1}(B_{\widehat R}\cap \mathcal C)}\Big\langle \dive\big((\widehat A-{\rm
    Id}_N)\nabla v\big) - \widehat {\mathbf b}\cdot\nabla v+ \widehat h
  v+\widehat f(x,v),\frac{\omega(|x|)}{|x|^{N-1}}\psi_i(x/|x|)
  \Big\rangle_{H^{1}_0(B_{\widehat R}\cap \mathcal C)}
\end{multline*}
for all $\omega\in \mathcal D(0,\widehat R)$.  Letting $\Upsilon_i$ as
in (\ref{Upsilon}), by direct calculations we have that
\begin{align}\label{eq:84}
  \Upsilon_i'(\lambda)= \lambda^{N-1}\zeta_i(\lambda)\quad\text{ in
  }\mathcal D'(0,\widehat R).
\end{align}
On the other hand, from the definition of $\zeta_i$ and the fact that
$v$ solves (\ref{eq:hatquation}), it follows that, for all $i$, the
function $\varphi_i$ defined in (\ref{78}) solves
\begin{equation*}
  -\varphi_i''(\lambda)-\frac{N-1}{\lambda}\varphi_i^\prime(\lambda)+
  \frac{\mu_i(V)}{\lambda^2}\varphi_i(\lambda)=
  \zeta_i(\lambda)\quad\text{in the sense of distributions in }(0,\widehat R),
\end{equation*}
which can be also written as
\begin{equation*}
  -\left(
    \lambda^{N-1+2\sigma_i}\big(\lambda^{-\sigma_i}\varphi_i(\lambda)\big)'\right)'
  =\lambda^{N-1+\sigma_i}
  \zeta_i(\lambda)\quad\text{in the sense of distributions in }(0,\widehat  R),
\end{equation*}
where
\begin{equation} \label{sigma_i}
  \sigma_i=-\frac{N-2}{2}+\sqrt{\bigg(\frac{N-2}
    {2}\bigg)^{\!\!2}+\mu_i(V)}.
\end{equation}
Let us fix $R\in (0,\widehat R)$. Integrating by parts the right hand side
and taking into account (\ref{eq:84}), we obtain that there exists
$c_i \in \R$ (depending on $R$) such that
\begin{equation*}
\big(
    \lambda^{-\sigma_i}\varphi_i(\lambda)\big)'
  =-\lambda^{-N+1-\sigma_i}\Upsilon_i(\lambda)-
\sigma_i\lambda^{-N+1-2\sigma_i}\bigg(
c_i+\int_\lambda^Rs^{\sigma_i-1}\Upsilon_i(s)ds\bigg)
\end{equation*}
in the sense of distributions in $(0,R)$. In particular
$\varphi_i\in W^{1,1}_{\rm loc}(0,\widehat R)$. A further integration yields
\begin{align}\label{eq:42-bis}
  \varphi_i(\lambda)&=\lambda^{\sigma_i} \bigg(
  R^{-\sigma_i}\varphi_i(R)+
  \int_\lambda^Rs^{-N+1-\sigma_i}\Upsilon_i(s)ds\bigg)\\[5pt]
  \notag&\qquad\qquad+\sigma_i\lambda^{\sigma_i}\int_\lambda^R
  s^{-N+1-2\sigma_i}\bigg(
  c_i+\int_s^Rt^{\sigma_i-1}\Upsilon_i(t)dt\bigg)ds\\[5pt]
  \notag&=\lambda^{\sigma_i} \bigg( R^{-\sigma_i}\varphi_i(R)+
  \frac{2-N-\sigma_i}{2-N-2\sigma_i}
  \int_\lambda^Rs^{-N+1-\sigma_i}\Upsilon_i(s)ds
  +\frac{\sigma_ic_iR^{-N+2-2\sigma_i}}{2-N-2\sigma_i} \bigg)\\[5pt]
  \notag&\qquad\qquad
+  \frac{\sigma_i\lambda^{-N+2-\sigma_i}}{N-2+2\sigma_i}
\bigg(c_i+\int_\lambda^R
  t^{\sigma_i-1}\Upsilon_i(t)\,dt\bigg).
\end{align}
Let $j_0,m\in\N$ be as in (\ref{eq:83}), so that the eigenvalue
$\mu_{j_0}(V)=\mu_{j_0+1}(V)=\cdots=\mu_{j_0+m-1}(V)$ has multiplicity
$m$ and
\begin{equation}\label{eq:85}
  \gamma=\lim_{r\rightarrow 0^+} {\mathcal N}(r)=\sigma_{i},
  \quad i=j_0,\dots,j_0+m-1,
\end{equation}
see  Lemma \ref{l:blowup}.
Estimate (\ref{eq:81}) and the Parseval identity yield
\begin{equation}\label{eq:86}
  H(\lambda)=(1+O(\lambda^\delta))\int_{C}|v(\lambda\,\theta)|^2\,d\sigma(\theta)=
  (1+O(\lambda^\delta))\sum_{i=1}^{\infty}|\varphi_i(\lambda)|^2,
  \quad\text{for all }0<\lambda\leq R.
\end{equation}
We claim that
\begin{equation}\label{eq:87}
\Upsilon_i(\lambda)=O(\lambda^{N-2+\hat\delta+\sigma_i})\quad
\text{for every }i\in\{j_0,\dots,j_0+m-1\}\quad
\text{as }\lambda\to0^+,
\end{equation}
with $\hat\delta$ defined in (\ref{eq:tilde_delta}).
Let us prove \eqref{eq:87}. By \eqref{1stest}, \eqref{1stest-2d},
\eqref{eq:71}, \eqref{eq:82}, H\"older inequality and a change of
variable we obtain
\begin{align} \label{eq:st-1} & \left| \int_{\mathcal C\cap B_\lambda}
    (\widehat A(x)-{\rm Id}_N)\nabla
    v(x)\cdot\frac{\nabla_{\SN}\psi_i(x/|x|)}{|x|}\,dx \right|
  \\
  & \notag \qquad \leq O(\lambda^{N-2+\delta+\sigma_i}) \int_{\mathcal
    C\cap B_1} |x|^\delta \frac{|\nabla\hat
    v^\lambda(x)|}{\sqrt{H(\lambda)}}
  \, \frac{|\nabla_{\SN}\psi_i(x/|x|)|}{|x|}\, dx \\
  & \notag \qquad \leq O(\lambda^{N-2+\delta+\sigma_i}) \left(
    \frac{\int_{\mathcal C\cap B_1} |\nabla \hat v^\lambda(x)|^2
      dx}{H(\lambda)} \right)^{\!\!1/2} \left( \int_{\mathcal C\cap
      B_1} |x|^{2\delta-2} |\nabla_{\SN}\psi_i(x/|x|)|^2 dx
  \right)^{\!\!1/2}\\
  & \notag \qquad =O(\lambda^{N-2+\delta+\sigma_i})\quad\text{as
  }\lambda\to 0^+.
\end{align}
Similarly, by a change of variable, \eqref{eq:3}, \eqref{eq:4},
\eqref{1stest}, \eqref{1stest-2d}, and boundedness in $H^1(B_1)$ of the
set $\{\widetilde w^\lambda\}_{\lambda}$ (see the proof of Lemma
\ref{l:blowup}), we obtain
\begin{align} \label{eq:st-2} \int_{\mathcal C\cap B_\lambda} \Big( -
  \widehat {\mathbf b}(x)\cdot\nabla v(x) + \widehat h(x)
  v(x)\Big)\psi_i(x/|x|)\,dx
  =O(\lambda^{N-2+\delta+\sigma_i})\quad\text{as }\lambda\to
  0^+.
\end{align}
Moreover, \eqref{eq:5_2}, Proposition
\ref{SMETS}, (\ref{1stest}), (\ref{1stest-2d}),
 and boundedness of 
$\{\widetilde w^\lambda\}_{\lambda}$ in $H^1(B_1)$ imply that 
\begin{align}\label{eq:st-4} 
  \bigg|&\int_{\mathcal C\cap B_\lambda} \widehat
  f(x,v(x))\psi_i(x/|x|)\,dx\bigg| \leq {\rm const\,}
  \bigg|\int_{\Omega_\lambda} \tilde f(y,w(y))\psi_i\big({\textstyle{
      \frac{\Phi^{-1}(y)}{|\Phi^{-1}(y)|}}}\big)\,dy\bigg|\\
  &\notag\leq {\rm const\,}\int_{\Omega_\lambda}(|w(y)|+|w(y)|^{\tilde
    p-1})\,dy\\
  &\notag \leq {\rm const\,}\lambda^N\sqrt{H(\lambda)}\bigg(
  \int_{\Omega_\lambda/\lambda}|w^\lambda(x)|\,dx+
  \int_{\Omega_\lambda/\lambda}|w(\lambda x)|^{\tilde
    p-2}|w^\lambda(x)|\,dx\bigg)\\
  &\notag \leq {\rm
    const\,}\lambda^{N+\sigma_i}\bigg(\|w^\lambda\|_{H^1(\Omega_\lambda/\lambda)}
  +\bigg(\int_{\Omega_{\lambda/\lambda}}|w(\lambda x)|^{\tilde
    q}dx\bigg)^{\!\!\frac{\tilde p-2}{\tilde
      q}}\|w^\lambda\|_{L^{\tilde
      p}(\Omega_\lambda/\lambda)}\bigg)\\
  &\notag \leq {\rm const\,}\lambda^{N+\sigma_i}(1+\lambda^{-N(\tilde
    p-2)/\tilde q}) =O(\lambda^{N-2+\sigma_i+(2-N\frac{\tilde
      p-2}{\tilde q})}) =O(\lambda^{N-2+\hat\delta+\sigma_i})
  \quad\text{as }\lambda\to 0^+,
\end{align}
where $\tilde p=2^*$ if $N\ge 3$ and $\tilde p=p$ with $p$ as in
(\ref{eq:F_assumption2}) if $N=2$, while $\tilde q=q$ with $q$ is as
in (\ref{eq:50}) if $N\ge 3$ and $\tilde q=2p$ if $N=2$, so that
$2-N\frac{\tilde p-2}{\tilde q}\geq \hat\delta$.  In order to estimate
the boundary term in \eqref{Upsilon}, we perform the change of
variables $x=\Phi^{-1}(y)$ and $y=\lambda\theta$ to obtain
\begin{align*}
  &\int_{\mathcal C\cap \partial B_\lambda} (\widehat A(x)-{\rm Id}_N)
  \nabla v(x)\cdot \frac{x}{|x|}\psi_i(x/|x|)\,d\sigma(x)\\
  & \notag = \lambda^{N-1}\int_{C_\lambda}(\widehat
  A(\Phi^{-1}(\lambda\theta)))-{\rm Id}_N) \big(\nabla
  w(\lambda\theta) {\rm Jac}\Phi(\Phi^{-1}(\lambda\theta))\big)^{T}
  \cdot \tfrac{\Phi^{-1}(\lambda\theta)}{|\Phi^{-1}(\lambda\theta)|}
  \times
  \\
  & \notag \qquad \qquad \times
  \psi_i\left(\tfrac{\Phi^{-1}(\lambda\theta)}{|\Phi^{-1}(\lambda\theta)|}
  \right) |{\rm det \, Jac}\Phi^{-1}(\lambda \theta)|\,
  d\sigma(\theta)
\end{align*}
and from this, using \eqref{eq:78}, \eqref{eq:79}, \eqref{eq:80},
\eqref{eq:71}, \eqref{1stest}, \eqref{1stest-2d}, (\ref{eq:doubling}),
(\ref{eq:doubling-2d}), and \eqref{eq:compactness}, we arrive to
\begin{align}\label{eq:st-3}
  \bigg|\int_{\mathcal C\cap \partial B_\lambda} &(\widehat A(x)-{\rm
    Id}_N)
  \nabla v(x)\cdot \frac{x}{|x|}\psi_i(x/|x|)\,d\sigma(x)\bigg|\\
  & \notag \leq O(\lambda^{N-2+\delta+\sigma_i})
  \left(\int_{C_\lambda} |\nabla w^\lambda (\theta)|^2
    d\sigma(\theta)\right)^{1/2} \left(\int_{C_\lambda} \left|
      \psi_i\left(\tfrac{\Phi^{-1}
          (\lambda\theta)}{|\Phi^{-1}(\lambda\theta)|}\right)
    \right|^2 d\sigma(\theta)
  \right)^{1/2}\\
  & \notag =O(\lambda^{N-2+\delta+\sigma_i}).
\end{align}
Inserting \eqref{eq:st-1}, \eqref{eq:st-2}, \eqref{eq:st-4},
\eqref{eq:st-3} into \eqref{Upsilon}, the proof of \eqref{eq:87}
follows.

In the rest of the proof it is not restrictive to assume that
$\sigma_i\neq 0$, since otherwise the proof of the lemma follows
immediately from \eqref{eq:42-bis}.  From \eqref{eq:87} we deduce that
the map
\begin{equation} \label{eq:int}
s\mapsto s^{-N+1-\sigma_i} \Upsilon_i(s)\in L^1(0,R)
\end{equation}
so that
\begin{multline}\label{eq:smo}
  \lambda^{\sigma_i} \bigg( R^{-\sigma_i}\varphi_i(R)+
  \frac{2-N-\sigma_i}{2-N-2\sigma_i}
  \int_\lambda^Rs^{-N+1-\sigma_i}\Upsilon_i(s)ds
  +\frac{\sigma_ic_iR^{-N+2-2\sigma_i}}{2-N-2\sigma_i} \bigg)
  \\
  =O(\lambda^{\sigma_i}) =o(\lambda^{-N+2-\sigma_i}) \qquad \text{as }
  \lambda\to 0^+ .
\end{multline}
On the other hand, by \eqref{eq:87} we also have that $t\mapsto
t^{\sigma_i-1} \Upsilon_i(t)\in L^1(0,R)$.  We now claim that
\begin{equation} \label{eq:claim2} c_i+\int_0^R t^{\sigma_i-1}
  \Upsilon_i(t)\, dt=0 .
\end{equation}
Suppose by contradiction that \eqref{eq:claim2} is not true. Then, by
\eqref{eq:42-bis} and \eqref{eq:smo} we infer
\begin{equation} \label{eq:contradiction} \varphi_i(\lambda) \sim
  \frac{\sigma_i}{N-2+2\sigma_i} \bigg(c_i+\int_0^R
  t^{\sigma_i-1}\Upsilon_i(t)\,dt\bigg)\lambda^{-N+2-\sigma_i}
  \quad\text{as }\lambda\to 0^+ .
\end{equation}
If $N\ge 3$, Hardy inequality and the fact that $v\in H^1(\mathcal
C\cap B_R)$ imply
\begin{equation*}
  \int_0^R \lambda^{N-3} |\varphi_i(\lambda)|^2 d\lambda\leq
  \int_0^R \lambda^{N-3} \left(\int_C |v(\lambda\theta)|^2
d\sigma(\theta) \right)  d\lambda
  =\int_{\mathcal C\cap B_R} \frac{|v(x)|^2}{|x|^2} \, dx<+\infty
\end{equation*}
thus contradicting \eqref{eq:contradiction}.
 If $N=2$, 
\eqref{eq:contradiction}, \eqref{sigma_i}, and the fact we are
assuming $\sigma_i\neq 0$ would imply
$$
\lim_{\lambda\to 0^+} |\varphi_i(\lambda)|=+\infty
$$
and, in turn, by \eqref{eq:86} we would have
$$
\lim_{\lambda\to 0^+} H(\lambda)=+\infty
$$
in contradiction with \eqref{1stest-2d}.
Claim \eqref{eq:claim2} is thereby proved.
By \eqref{eq:42-bis} and \eqref{eq:claim2} we then obtain
\begin{align} \label{eq:be}
&
\varphi_i(\lambda)=\lambda^{\sigma_i} \bigg( R^{-\sigma_i}\varphi_i(R)+
  \frac{2-N-\sigma_i}{2-N-2\sigma_i}
  \int_\lambda^Rs^{-N+1-\sigma_i}\Upsilon_i(s)ds
  +\frac{\sigma_ic_iR^{-N+2-2\sigma_i}}{2-N-2\sigma_i} \bigg)\\[5pt]
  \notag&\qquad\qquad
- \frac{\sigma_i\lambda^{-N+2-\sigma_i}}{N-2+2\sigma_i}
\int_0^\lambda t^{\sigma_i-1}\Upsilon_i(t)\,dt .
\end{align}
The proof the lemma follows inserting \eqref{eq:claim2} into
\eqref{eq:be} and observing that by \eqref{eq:87}
$$
\frac{\sigma_i\lambda^{-N+2-\sigma_i}}{N-2+2\sigma_i} \int_0^\lambda
t^{\sigma_i-1}\Upsilon_i(t)\,dt=O(\lambda^{\sigma_i+\hat\delta})
$$
as $\lambda\to 0^+$.
\end{pf}

\noindent The asymptotic behavior of $H(\lambda)$ as $\lambda\to0^+$
is evaluated in the following lemma.

\begin{Lemma}\label{l:limitHgammapos}
Let $H$ be as in (\ref{eq:H(r)}) and let $\gamma$ be as
  in Lemmas \ref{l:limitN} and \ref{l:limitN-2d} respectively in the cases
  $N\ge 3$ and $N=2$. Then
\begin{equation} \label{eq:limit-H}
\lim_{\lambda\to 0^+} \lambda^{-2\gamma} H(\lambda)>0 .
\end{equation}
\end{Lemma}

\begin{pf} The fact that the limit in \eqref{eq:limit-H} exists and is
  finite was proved in Lemma \ref{l:lim-exist} and hence we may proceed by
  contradiction by supposing that $\lim_{\lambda\to 0^+}
  \lambda^{-2\gamma} H(\lambda)=0$.
Let $j_0$ and $m$ be as in
   (\ref{eq:83}) and $\varphi_i$ as in (\ref{78}).
From \eqref{eq:86} we
  deduce that for any $i\in \{j_0,\dots,j_0+m-1\}$
  \begin{equation*}
    \lim_{\lambda\to 0^+} \lambda^{-\gamma} \varphi_i(\lambda)=0 .
\end{equation*}
Therefore by Lemma \ref{l:phii}, \eqref{eq:85}, and \eqref{eq:int} we obtain
\begin{equation*}
  R^{-\gamma}\varphi_i(R) -\frac{\gamma R^{-N+2-2\gamma}}{2-N-2\gamma}
  \int_0^R s^{\gamma-1} \Upsilon_i(s)\,ds
  =-\frac{2-N-\gamma}{2-N-2\gamma}
  \int_0^R s^{-N+1-\gamma}\Upsilon_i(s)ds
\end{equation*}
which replaced in \eqref{eq:id} yields
\begin{equation*}
  \varphi_i(\lambda)=
  -\frac{2-N-\gamma}{2-N-2\gamma} \, \lambda^{\gamma}
  \int_0^\lambda s^{-N+1-\gamma}\Upsilon_i(s) \, ds +O(\lambda^{\gamma+\hat\delta})
  \quad \text{as } \lambda\to 0^+.
\end{equation*}
Inserting in the last estimate \eqref{eq:87} we conclude that for any
$i\in\{j_0,\dots,j_0+m-1\}$
\begin{equation*}
\varphi_i(\lambda)=O(\lambda^{\gamma+\hat\delta}) \qquad \text{as } \lambda\to 0^+ .
\end{equation*}
This, together with \eqref{78}, implies
\begin{equation} \label{eq:gamma+e} \int_C
  v(\lambda\theta)\phi(\theta)\,
  d\sigma(\theta)=O(\lambda^{\gamma+\hat\delta}) \quad \text{as }
  \lambda\to 0^+
\end{equation}
for any function $\phi$ in the eigenspace $\mathcal E_0$. By
\eqref{eq:gamma+e}, \eqref{eq:81}, \eqref{2ndest} in the case $N\geq
3$, and \eqref{2ndest-1} with $\sigma=\hat\delta$ in the case $N=2$, we
obtain
\begin{equation}\label{eq:o(1)}
  \int_C \frac{v(\lambda\theta)}{\|\hat
    v^\lambda\|_{L^2(C)}} \, \phi(\theta) \,
  d\sigma(\theta)=O(\lambda^{\hat\delta/2})=o(1) \qquad
  \text{as } \lambda\to 0^+
\end{equation}
for any $\phi\in \mathcal E_0$.  On the other hand, Lemma
\ref{l:blow_v} states that for any sequence $\lambda_n\to 0^+$ there
exists a subsequence $\lambda_{n_k}$ and a function $\psi\in\mathcal
E_0$ with $\|\psi\|_{L^2(C)}=1$ such that
\begin{equation*}
\frac{v(\lambda_{n_k}\cdot)}{\|\hat v^{\lambda_{n_k}}\|_{L^2(C)}}\to \psi
\end{equation*}
strongly in $L^2(C)$.
Therefore, taking $\phi=\psi$ in \eqref{eq:o(1)} we conclude that
$$
0=\lim_{k\to+\infty}\Big(
\frac{\hat v^{\lambda_{n_k}}}{\|\hat v^{\lambda_{n_k}}\|_{L^2(C)}},\psi\Big)_{L^2(C)}
=\|\psi\|^2_{L^2(C)}=1
$$
thus giving rise to a contradiction.
\end{pf}

\noindent 
The following theorem is a more precise and complete version of Theorem
\ref{t:main-w2}.

\begin{Theorem}\label{t:main-w}
  Let $\tildA,\tildb,\tilde f,\widetilde h$ be as in
  (\ref{eq:tildA}--\ref{eq:tildh}) with $A,\bi,\Psi,f,h,V$ as in
  assumptions (\ref{eq:matrix1}--\ref{eq:F_assumption2}),
  (\ref{eq:def_psi}). Let $\tildO$ be as in (\ref{eq:tildeomega}) with
  $\Omega$ satisfying (\ref{eq:omega}) and
  (\ref{eq:phi1}--\ref{eq:phi4}).  Let $w\in
  H^1(\tildO)\setminus\{0\}$ be a non-trivial weak solution to
  (\ref{eq:tildequation}).  Then, letting $\mathcal N(r)$ as in
  (\ref{eq:alm_fun}) and (\ref{Almgren}), there exists $k_0\in \N$, $k_0\geq 1$, such that
\begin{equation}\label{lim-N}
  \lim_{r\to 0^+} \mathcal
  N(r)=-\frac{N-2}2+\sqrt{\left(\frac{N-2}2\right)^{\!\!2}+\mu_{k_0}(V)}.
\end{equation}
Furthermore, if $\gamma$ denotes the limit in (\ref{lim-N}), $m\geq 1$
is the multiplicity of the eigenvalue $\mu_{k_0}(V)$ and
$\{\psi_i:j_0\leq i\leq j_0+m-1\}$ ($j_0\leq k_0\leq j_0+m-1$) is an
$L^2(C)$-orthonormal basis for the eigenspace associated to
$\mu_{k_0}(V)$, then, denoting again as $w$ its trivial extension
outside $\tildO$,
\begin{equation} \label{convergence} \lambda^{-\gamma} w(\lambda x)\to
  |x|^\gamma\sum_{i=j_0}^{j_0+m-1} \beta_i
  \psi_i\bigg(\frac{x}{|x|}\bigg) \quad \text{as } \lambda\to 0^+
\end{equation}
in $H^1(B_1)$ and  in $C^{1,\alpha}_{{\rm loc}}(\mathcal
C\cap B_1)$ for any $\alpha\in (0,1)$,  where
$$
(\beta_{j_0},\beta_{j_0+1},\dots,\beta_{j_0+m-1})\neq(0,0,\dots,0)
$$
and 
\begin{align}\label{eq:88}
  \beta_i= \int_{C} R^{-\gamma}w(\Phi(R\theta))\psi_i(\theta)\,d\sigma(\theta)
+\frac{1}{2-N-2\gamma}
\int_0^R\bigg(\frac{2-N-\gamma}{s^{N-1+\gamma}}-\gamma
\frac{s^{\gamma-1}}{R^{N-2+2\gamma}}
\bigg)\Upsilon_i(s)\,ds,
 \end{align}
 for all $R\in (0,\widehat R)$ for some $\widehat R>0$, 
$\Upsilon_i$ being defined in (\ref{Upsilon}).
\end{Theorem}

\begin{pf}
  Identity \eqref{lim-N} follows immediately from Lemma
  \ref{l:blowup}.  As in the statement of the theorem, let $m$ be the
  multiplicity of the eigenvalue $\mu_{k_0}(V)$ found in Lemma
  \ref{l:blowup},  $j_0\in \N\setminus\{0\}$,  such that
$j_0\leq k_0\leq j_0+m-1$,
  $\mu_{j_0}(V)=\mu_{j_0+1}(V)=\dots=\mu_{j_0+m-1}(V)$, and
  $\gamma=\lim_{r\to 0^+} \mathcal N(r)$.

  In order to prove \eqref{convergence}, let $\{\lambda_n\}_{n\in
    \N}\subset (0,\infty)$ be a sequence such that $\lambda_n\to 0^+$
  as $n\to +\infty$. By Lemmas \ref{l:blowup},
  \ref{l:lim-exist}, \ref{l:blow_v},  \ref{l:limitHgammapos},
and (\ref{eq:81}),
there
  exist a subsequence $\lambda_{n_j}$ and
  $\beta_{j_0},\dots,\beta_{j_0+m-1}\in \R$ such that
  $(\beta_{j_0},\beta_{j_0+1},\dots,\beta_{j_0+m-1})\neq(0,0,\dots,0)$,
\begin{equation}\label{eq:23-bis}
  \lambda_{n_j}^{-\gamma}w(\lambda_{n_j}x)\to
  |x|^\gamma\sum_{i=j_0}^{j_0+m-1} \beta_i\psi_{i}\bigg(\frac{x}{|x|}\bigg)
\quad\text{ in $H^1(B_1)$ and $C^{1,\alpha}_{{\rm loc}}(\mathcal
C\cap B_1)$ for any $\alpha\in (0,1)$}
\end{equation}
 (with $w$ meant to be
trivially extended outside $\tildO$), and
\begin{equation}\label{eq:89}
\lambda_{n_j}^{-\gamma}v(\lambda_{n_j}\cdot)
\to
\sum_{i=j_0}^{j_0+m-1} \beta_i\psi_{i}\quad \text{in }
L^2(C) \quad \text{as }j\to+\infty.
\end{equation}
We now prove that the $\beta_i$'s depend neither on the sequence
$\{\lambda_n\}_{n\in\N}$ nor on its subsequence
$\{\lambda_{n_j}\}_{j\in\N}$.  Let us fix $R\in(0,\widehat R)$ with
$\widehat R$ as in Lemma \ref{l:straightening}.  Defining $\varphi_i$
as in (\ref{78}), from (\ref{eq:89}) it follows that, for any
$i=j_0,\dots, j_0+m-1$,
\begin{equation}\label{eq:25-bis}
  \lambda_{n_j}^{-\gamma}\varphi_i(\lambda_{n_j}) =\int_{C}
  \frac{v(\lambda_{n_j}\theta)}{\lambda_{n_j}^{\gamma}}
  \psi_i(\theta)\,d\sigma(\theta)
  \to\sum_{\ell=j_0}^{j_0+m-1} \beta_\ell\int_{{\mathbb
      S}^{N-1}}\psi_{\ell}(\theta)\psi_i(\theta)\,d\sigma(\theta)=\beta_i
\end{equation}
as $j\to+\infty$.  On the other hand, from Lemma \ref{l:phii},
it follows that,  for any $i=j_0,\dots, j_0+m-1$,
\begin{align*}
  \lambda^{-\gamma}\varphi_i(\lambda)\ \to\
  &
 R^{-\gamma}\varphi_i(R)+ \frac{2-N-\gamma}{2-N-2\gamma}
  \int_0^Rs^{-N+1-\gamma}\Upsilon_i(s)\,ds- \frac{\gamma
    R^{-N+2-2\gamma}}{2-N-2\gamma} \int_0^R s^{\gamma-1}\Upsilon_i(s) \,ds
\end{align*}
as $\lambda\to0^+$
and therefore from (\ref{eq:25-bis}) we deduce that
\begin{align*}
  \beta_i= R^{-\gamma}\varphi_i(R)+ \frac{2-N-\gamma}{2-N-2\gamma}
  \int_0^Rs^{-N+1-\gamma}\Upsilon_i(s)\,ds- \frac{\gamma
    R^{-N+2-2\gamma}}{2-N-2\gamma} \int_0^R s^{\gamma-1}\Upsilon_i(s) \,ds,
\end{align*}
 for any $i=j_0,\dots, j_0+m-1$.
In particular the $\beta_i$'s depend neither on the sequence
$\{\lambda_n\}_{n\in\N}$ nor on its subsequence
$\{\lambda_{n_k}\}_{k\in\N}$, thus implying that the convergence in
(\ref{eq:23-bis}) actually holds as $\lambda\to 0^+$
and proving the theorem.
\end{pf}

\begin{pfn}{Theorem \ref{t:main-u}}
Let us first observe that the family of functions
\begin{equation*}
u^\lambda(x)=\frac{u(\lambda x)}{\sqrt{H(\lambda)}}
\end{equation*}
is bounded in $H^1(B_1)$, where $u$ is meant to be trivially extended
outside $\Omega$. Indeed, by the change of variable $\lambda
x=\Psi(\lambda y)$
\begin{align*}
  &\int_{B_1}|\nabla u^\lambda(x)|^2dx=
  \int_{\frac{\Psi^{-1}(B_\lambda)}{\lambda}} \Big|\nabla \widetilde
  w^\lambda(y)(\mathop{\rm
    Jac}\Psi(\lambda y))^{-1}\Big|^2 |{\rm det \, Jac}\Psi(\lambda y)|\,dy,\\
  &\int_{B_1}|u^\lambda(x)|^2dx=
  \int_{\frac{\Psi^{-1}(B_\lambda)}{\lambda}} |\widetilde
  w^\lambda(y)|^2 |{\rm det \, Jac}\Psi(\lambda y)|\,dy,
\end{align*}
and hence from (\ref{eq:def_psi}), (\ref{eq:doubling}),
(\ref{eq:doubling-2d}), and boundedness in $H^1(B_1)$ of the set
$\{\widetilde w^\lambda\}_{\lambda}$ it follows that
$\{u^\lambda\}_\lambda$ is bounded in $H^1(B_1)$ uniformly with
respect to $\lambda$. Hence we can repeat for $u^\lambda$ the same
arguments performed in the proof of Lemma \ref{l:blowup} for
$w^\lambda$ to obtain that $\{u^{\lambda}\}_\lambda$ is relatively
compact in $C^{1,\alpha}_{\rm loc}(\mathcal C \cap B_1)$, in
$C^{0,\alpha}_{{\rm loc}}(B_1\setminus\{0\})$, and in $H^1(B_1)$ and
hence, by Lemma \ref{l:limitHgammapos},
\begin{align}\label{eq:90}
\{\lambda^{-\gamma}u(\lambda\cdot)\}_\lambda\text{ is relatively compact in }
C^{1,\alpha}_{\rm loc}(\mathcal C \cap B_1),
\text{ in }C^{0,\alpha}_{{\rm loc}}(B_1\setminus\{0\}),
\text{ and in }H^1(B_1),
\end{align}
with $\gamma$ as in Theorem \ref{t:main-w}.
Furthermore, from
Lemma \ref{l:limitHgammapos}, (\ref{eq:def_psi}), and (\ref{eq:65})
\begin{align*}
  |\lambda^{-\gamma}u(\lambda x)-\lambda^{-\gamma}w(\lambda x)|=
  \frac{\sqrt{H(\lambda)}}{\lambda^\gamma}\bigg(w^\lambda\Big(\frac{\Psi^{-1}(\lambda
    x)}{\lambda}\Big)-w^\lambda (x)\bigg)\to 0
\end{align*}
for all $x\in \mathcal C\cap B_1$. From the above limit and Theorem
\ref{t:main-w} we deduce that, for all $x\in \mathcal C\cap B_1$,
\begin{align}\label{eq:91}
  \lambda^{-\gamma}u(\lambda x)
\to |x|^\gamma\sum_{i=j_0}^{j_0+m-1}
\beta_i \psi_i\bigg(\frac{x}{|x|}\bigg) \quad
\text{as } \lambda\to 0^+
\end{align}
with $\beta_i$ and $\psi_i$ as in  Theorem \ref{t:main-w}.
Combining (\ref{eq:90}) and (\ref{eq:91}) we deduce that
\begin{align*}
  \lambda^{-\gamma}u(\lambda x)
\to |x|^\gamma\sum_{i=j_0}^{j_0+m-1}
\beta_i \psi_i\bigg(\frac{x}{|x|}\bigg) \quad
\text{as } \lambda\to 0^+
\end{align*}
in $C^{1,\alpha}_{\rm loc}(\mathcal C \cap B_1)$,
 in $C^{0,\alpha}_{{\rm loc}}(B_1\setminus\{0\})$
for all
$\alpha\in(0,1)$, and in $H^1(B_1)$, thus completing the~proof.\end{pfn}

\section{An example}\label{sec:examples}
In this section we show that the presence of a logarithmic term in the asymptotic 
expansion cannot be excluded without assuming conditions \eqref{eq:phi2} and
\eqref{eq:63}. 

Let us consider a domain $\Omega\subset \R^2$ admitting a local
representation in a neighborhood of the origin as in
\eqref{eq:omega} where the corresponding function $\varphi$
satisfies \eqref{eq:phi1},
$$
\sup_{\nu\in {\mathbb S}^{N-2}}\Big|\frac{\varphi(t\nu)}{t}-g(\nu)\Big|=o(1)
\quad\text{as }t\to0^+,
$$
with $g$ as in (\ref{eq:g_continuous}), but not
(\ref{eq:phi2}--\ref{eq:63}). To this purpose let us define in
Gauss plane the sets
\begin{equation*}
A_1:=\C\setminus \{z=ix_2\in\C:x_2\leq 0\}, \qquad
A_2:=\C\setminus\{x_1\in \R\subset \C:x_1\leq 0\}
\end{equation*}
and the holomorphic functions $\eta_1:A_1\to\C$, $\eta_2:A_2\to\C$
defined as follows:
\begin{align*}
&\eta_1(z):=\log r+i\theta \quad \text{for any $z=re^{i\theta}\in A_1$, 
 $r>0$,
$\theta\in\left(-\frac{\pi}2,\frac{3\pi}2\right)$},\\
&\eta_2(z):=\log r+i\theta \quad \text{for any $z=re^{i\theta}\in A_2$,
 $r>0$,
$\theta\in\left(-\pi,\pi\right)$}.
\end{align*}
Given $\alpha\in (0,2)$, we are going to define $\Omega$ in such a way
that $\partial\Omega$ admits at $0$ a corner with amplitude
$\alpha\pi$.  We distinguish the cases $\alpha\in(0,1)$, and
$\alpha\in [1,2)$.

{\bf The case $\alpha\in (0,1)$.}  \ Let us consider the
holomorphic function
$$
v_1(z):=e^{\frac{2}{\alpha}\eta_1(-iz)}\eta_1(-iz) \qquad \text{for
  any } z\in\{w\in \C:\Im w>0\}
$$
and the set
\begin{equation} \label{Zeta}
\mathcal Z_1:=\{z: \Im z>0 \text{ and } \Im(v_1(z))=0\} .
\end{equation}
If $z=re^{i\theta}$ with $r>0$, $\theta\in (0,\pi)\setminus
\{\frac{\pi}2\}$, then $z\in \mathcal Z_1$ if and only
\begin{equation} \label{rho_1} r=\rho_1(\theta):=
  \exp\left[-\bigg(\theta-\frac{\pi}2\bigg)\cot\bigg(\frac{2}{\alpha}
    \bigg(\theta-\frac{\pi}2\bigg) \bigg) \right] .
\end{equation}
For some fixed $\sigma\in \left(0,\frac\pi2(1-\alpha)\right)$, 
we define the curves $\gamma_{\sigma}^+\subset
\mathcal Z_1$ and $\gamma_{\sigma}^-\subset \mathcal Z_1$
respectively parametrized by
\begin{equation} \label{gamma+}
\gamma_{\sigma}^+:
\begin{cases}
x_1(\theta)=\rho_1(\theta)\cos\theta \\
x_2(\theta)=\rho_1(\theta)\sin\theta
\end{cases}
\qquad \theta\in
\left(\frac{\pi}2-\frac{\alpha\pi}2-\sigma,\frac{\pi}2-\frac{\alpha\pi}2\right)
\end{equation}
and
\begin{equation} \label{gamma-}
\gamma_{\sigma}^-:
\begin{cases}
x_1(\theta)=\rho_1(\theta)\cos\theta \\
x_2(\theta)=\rho_1(\theta)\sin\theta
\end{cases}
\qquad \theta\in
\left(\frac{\pi}2+\frac{\alpha\pi}2,\frac{\pi}2+\frac{\alpha\pi}2+\sigma\right)
  .
\end{equation}
If we choose $\sigma>0$ sufficiently small then the union of these
two curves is the graph of a function $\varphi$ defined in a
neighborhood $U$ of $0$. Moreover $\varphi$ is a Lipschitz function in
$U$, $\varphi\in C^1(U\setminus\{0\})$ and
\begin{equation} \label{derivatives}
\lim_{t\to 0^-}
\frac{\varphi(t)}{t}=\tan\left(\frac{\pi}2+\frac{\alpha\pi}2\right)
\ , \qquad \lim_{t\to 0^+}
\frac{\varphi(t)}{t}=\tan\left(\frac{\pi}2-\frac{\alpha\pi}2\right)
.
\end{equation}
At this point it is possible to construct a bounded domain
$\Omega\subset \{z\in \C:\Im z>0\}$ satisfying \eqref{eq:omega}
for some $R>0$ sufficiently small. Then we define the harmonic function
\begin{equation*}
u(x_1,x_2):=\Im(v_1(z)) \qquad \text{ for any } z=x_1+ix_2\in
\Omega .
\end{equation*}
In polar coordinates the function $u$ reads
\begin{equation} \label{polar}
u(r,\theta)=r^{\frac{2}{\alpha}} \left[(\log r)
\sin\left(\frac{2}{\alpha}
\left(\theta-\frac{\pi}2\right)\right)+\left(\theta-\frac{\pi}2\right)
\cos\left(\frac{2}{\alpha}
\left(\theta-\frac{\pi}2\right)\right)\right] .
\end{equation}
Since $\Omega$ is bounded, then $u\in H^1(\Omega)$. From
(\ref{Zeta}--\ref{rho_1}) and \eqref{polar} we deduce that $u$
vanishes on $\gamma^+_\sigma\cup \gamma^-_\sigma$ and in particular on
$\partial\Omega\cap B_R$.

Next we show that $\varphi$ does not satisfy \eqref{eq:phi3}
for any $C_0>0$. Since by
(\ref{rho_1}--\ref{gamma-}) $\varphi$ is an even
function, it is sufficient to study the behavior of $\varphi(x_1)-x_1
\varphi'(x_1)$ in a right neighborhood of zero.  By \eqref{gamma+} and
the fact that $\alpha\in (0,1)$ we may assume that $\theta\in
\big(0,\frac{\pi}2\big)$ and hence, if $x_1$ belongs to a
sufficiently small right neighborhood of $0$, by \eqref{rho_1} we have
\begin{equation} \label{identity} \frac{1}2 \log
  \big(x_1^2+(\varphi(x_1))^2\big)
  \tan\left[\frac{2}{\alpha}\left(\arctan\left(\frac{\varphi(x_1)}{x_1}\right)
      -\frac{\pi}2\right)\right]
  +\arctan\left(\frac{\varphi(x_1)}{x_1}\right)-\frac{\pi}2=0 .
\end{equation}
By \eqref{derivatives} and \eqref{identity} we have that, as $x_1\to 0^+$,
\begin{align} \label{tangent} \tan &
  \left[\frac{2}{\alpha}\left(\arctan\Big(\frac{\varphi(x_1)}{x_1}\Big)
      -\frac{\pi}2\right)\right]
  =-\frac{2\arctan\big(\frac{\varphi(x_1)}{x_1}\big)-\pi}{\log
    \big(x_1^2+(\varphi(x_1))^2\big)}
   =\frac{\alpha\pi}2 \frac{1}{\log
    x_1}+o\left(\frac{1}{\log x_1}\right)   .
\end{align}
Differentiating both sides of \eqref{identity} and multiplying by
$x_1^2+(\varphi(x_1))^2$ we obtain the identity
\begin{multline*}
  \big(x_1+\varphi(x_1)\varphi'(x_1)\big) \, \tan
  \Big[\tfrac{2}{\alpha}\Big(\arctan\big(\tfrac{\varphi(x_1)}{x_1}\big)
  -\tfrac{\pi}2\Big)\Big]\\
  +\left\{ 1+\frac{\log \big(x_1^2+(\varphi(x_1))^2\big)}
    {\alpha\cos^2\big[\frac{2}{\alpha}
      \big(\arctan\big(\frac{\varphi(x_1)}{x_1}\big)
      -\frac{\pi}2\big)\big]} \right\}
  \big(x_1\varphi'(x_1)-\varphi(x_1) \big)=0
\end{multline*}
and hence \eqref{derivatives} and \eqref{tangent} yield
\begin{equation} \label{stima} x_1\varphi'(x_1)-\varphi(x_1)\sim
  -\frac{\alpha^2 \pi\left[1+\tan^2
      \left(\frac{\pi}2-\frac{\alpha\pi}2\right)\right]}4 \,
  \frac{x_1}{\log^2 x_1} \qquad \text{as } x_1\to 0^+ .
\end{equation}
This shows that $\varphi$ does not satisfy condition \eqref{eq:phi3}.

{\bf The case $\alpha\in [1,2)$.}  \  Let us consider the
holomorphic function
$$
v_2(z):=e^{\frac{2}{\alpha}\eta_2(-iz)}\eta_2(-iz) \qquad \text{for
  any } z\in \C\setminus \{iy: y\le 0\}
$$
and the set
\begin{equation*}
\mathcal Z_2:=\{z\in \C\setminus \{iy: y\le 0\}:\Im(v_2(z))=0\} .
\end{equation*}
Similarly to the previous case one may define two curves
$\gamma_\sigma^+$, $\gamma_\sigma^-$ and a corresponding function
$\varphi$ whose graph coincides with $\gamma_\sigma^+\cup
\gamma_\sigma^-$.  Next one may also construct a bounded domain
$\Omega$ satisfying \eqref{eq:omega} for a suitable choice of $R>0$
and define a harmonic function $u$ as
$$
u(x_1,x_2):=\Im(v_2(z)) \qquad \text{for any } z=x_1+ix_2\in \Omega. 
$$
Finally one can prove that the new function $\varphi$ satisfies
\eqref{stima}.

\bigskip\noindent 
 \textbf{Aknowledgements.} The authors would like to thank
 Prof. Susanna Terracini  for  helpful comments and discussions.


\begin{thebibliography}{99}
\bibitem{ae} V. Adolfsson, L. Escauriaza, {\it
$C^{1,\alpha}$ domains and unique continuation at the boundary},
Comm. Pure Appl. Math. 50 (1997), no. 10, 935--969.

\bibitem{almgren} F. J. Jr. Almgren, {\it $Q$ valued functions
    minimizing Dirichlet's integral and the regularity of area
    minimizing rectifiable currents up to codimension two}, Bull.
  Amer. Math. Soc. 8 (1983), no. 2, 327--328.

\bibitem{nistor} C. B{\u{a}}cu{\c{t}}{\u{a}}, A. L. Mazzucato, V.
  Nistor, L.  Zikatanov, {\it Interface and mixed boundary value
    problems on $n$-dimensional polyhedral domains}, Doc. Math. 15
  (2010), 687--745.

\bibitem{BrezisKato} H.~Br{\'e}zis, T.~Kato, {\it Remarks on the
    {S}chr\"odinger operator with singular complex potentials}, J.
  Math. Pures Appl. (9) 58 (1979), no.~2, 137--151.

\bibitem{cd} M. Costabel, M. Dauge, {\it Singularit\'es d'ar\^etes
    pour les probl\`emes aux limites elliptiques}, Journ\'ees
  ``\'Equations aux D\'eriv\'ees Partielles'' (Saint-Jean-de-Monts,
  1992), Exp. No. IV, 12 pp., \'Ecole Polytech., Palaiseau, 1992.


\bibitem{dauge} M. Dauge, {\it Elliptic boundary value problems on
    corner domains.  Smoothness and asymptotics of solutions}, Lecture
  Notes in Mathematics, 1341.  Springer-Verlag, Berlin, 1988.

\bibitem{FFT} V. Felli, A. Ferrero, S. Terracini, {\it Asymptotic
    behavior of solutions to Schr\"odinger equations near an isolated
    singularity of the electromagnetic potential},
  J. Eur. Math. Soc. (JEMS) 13 (2011), no. 1, 119--174.

\bibitem{FFT2} V. Felli, A. Ferrero, S. Terracini, {\it On the
    behavior at collisions of solutions to Schr\"odinger equations
    with many-particle and cylindrical potentials}, Preprint 2010.

\bibitem{FFT3} V. Felli, A. Ferrero, S. Terracini, {\it A note on
    local asymptotics of solutions to singular elliptic equations via
    monotonicity methods}, Preprint 2011.

\bibitem{GL} N. Garofalo, F.-H.  Lin, {\it Monotonicity properties of
    variational integrals, $A\sb p$ weights and unique continuation},
  Indiana Univ. Math. J.  35 (1986), no. 2, 245--268.

\bibitem{Grisvard} P.  Grisvard, {\it Behavior of the solutions of an
    elliptic boundary value problem in a polygonal or polyhedral
    domain}, Numerical solution of partial differential equations, III
  (Proc. Third Sympos. (SYNSPADE), Univ. Maryland, College Park, Md.,
  1975), pp. 207--274. Academic Press, New York, 1976.

\bibitem{kawohl} B. Kawohl, {\it On nonlinear mixed boundary value
    problems for second order elliptic differential equations on
    domains with corners}, Proc. Roy. Soc. Edinburgh Sect. A 87
  (1980/81), no. 1--2, 35--51.

  \bibitem{ko} V. A. Kondrat'ev, O. A. Ole\u{\i}nik, {\it Boundary
      value problems for partial differential equations in nonsmooth
      domains}, Uspekhi Mat. Nauk 38 (1983), no. 2(230), 3--76.

  \bibitem{lehman} R. S. Lehman, {\it Development of the mapping
      function at an analytic corner}, Pacific J. Math. 7 (1957),
    1437--1449.

\bibitem{llm} V. Liskevich, S. Lyakhova, V. Moroz, {\it Positive
    solutions to singular semilinear elliptic equations with critical
    potential on cone-like domains}, Adv. Differential Equations 11
  (2006), no. 4, 361--398.

\bibitem{Mazja} V. G. Maz'ja, {\it Sobolev Spaces}, Springer Series in
  Soviet Mathematics, Springer-Verlag, Berlin, 1985.

\bibitem{MNP1} V. Maz'ya, S. Nazarov, B. Plamenevskij, {\it Asymptotic
    theory of elliptic boundary value problems in singularly perturbed
    domains. Vol. I}, Operator Theory: Advances and Applications, 111.
  Birkh\"auser Verlag, Basel, 2000.

\bibitem{tz} X. Tao, S. Zhang, {\it Weighted doubling properties and
    unique continuation theorems for the degenerate Schr\"odinger
    equations with singular potentials}, J. Math. Anal. Appl. 339
  (2008), no. 1, 70--84.

  \bibitem{terracini96} S. Terracini, {\it On positive entire solutions
      to a class of equations with singular coefficient and critical
      exponent}, Adv. Diff. Equa. 1 (1996), no. 2, 241--264.

  \bibitem{tolksdorf} P. Tolksdorf, {\it On the Dirichlet problem for
      quasilinear equations in domains with conical boundary points},
    Comm. Partial Differential Equations 8 (1983), no. 7, 773--817.

  \bibitem{wz} Z.-Q. Wang, M. Zhu, {\it {H}ardy inequalities with
      boundary terms}, Electron. J. Differential Equations 2003, No.
    43.

  \bibitem{wasow} W. Wasow, {\it Asymptotic development of the
      solution of Dirichlet's problem at analytic corners}, Duke Math.
    J. 24 (1957), 47--56.

  \bibitem{wigley} N. M. Wigley, {\it Corner behavior of solutions of
      semilinear Dirichlet problems}, Canad. J. Math. 37 (1985), no.
    6, 1025--1046.

  \bibitem{wolff} T. H. Wolff, {\it A property of measures in $\R^ N$
      and an application to unique continuation}, Geom. Funct. Anal.
    2 (1992), no. 2, 225--284.


\end{thebibliography}
\end{document}